\documentclass[a4paper,11pt,twoside]{article}
\usepackage[T1]{fontenc}
\usepackage[english]{babel}
\usepackage{amsmath, amssymb, amsthm, amsfonts, mathtools, mathrsfs}	
\usepackage{enumerate}
\usepackage[shortlabels]{enumitem}
\usepackage{booktabs, caption, graphicx, subfig, wrapfig} 	
\usepackage{geometry}
\usepackage{algorithm}
\usepackage{algpseudocode} 
\usepackage{epstopdf} 
\usepackage{soul}

\usepackage{xcolor}

\usepackage{cite} 

\usepackage{fancyhdr} 

\usepackage{upgreek}
\usepackage{bbm}

\newcommand{\R}{\mathbb{R}}
\newcommand{\N}{\mathbb{N}}

\newcommand{\C}{\mathbb{C}}

\numberwithin{equation}{section}

\theoremstyle{plain}
\newtheorem{theorem}{Theorem}[section] 

\newtheorem{corollary}[theorem]{Corollary}

\theoremstyle{definition}

\newtheorem{remark}[theorem]{Remark}

\DeclarePairedDelimiter{\abs}{\lvert}{\rvert}
\DeclarePairedDelimiter{\norm}{\lVert}{\rVert}
\DeclarePairedDelimiter{\duality}{\langle}{\rangle}

\DeclareMathOperator*{\argmin}{arg\,min}

\newcommand{\eps}{\varepsilon}
\renewcommand{\phi}{\varphi}
\renewcommand{\bar}{\overline}
\renewcommand{\vec}{\boldsymbol}

\def\genspazio #1#2#3#4#5{#1^{#2}(#5,#4;#3)}
\def\spazio #1#2#3{\genspazio {#1}{#2}{#3}T0}

\def\LT {\spazio L}
\def\HT {\spazio H}
\def\C #1#2{\mathcal{C}^{#1}([0,T];#2)}

\def\Lx #1{L^{#1}(\Omega)}
\def\Lt #1{L^{#1}(0,T)}
\def\Lqt #1{L^{#1}(Q_T)}
\def\Hx #1{H^{#1}(\Omega)}
\def\Wx #1{W^{#1}(\Omega)}

\def\Cqt #1{\mathcal{C}^{#1}(\bar{Q_T})}

\def\Accorpa #1#2 #3 {\gdef #1{\eqref{#2}--\eqref{#3}}%
	\wlog{}\wlog{\string #1 -> #2 - #3}\wlog{}}%
	
\def\ls{<}
\def\gs{>}
\def\mezzo {\frac{1}{2}}
\def\de {\mathrm{d}}
\def\D {\mathrm{D}}
\def\ddt {\frac{\de}{\de t}}

\def\weakstar {\stackrel{\ast}{\rightharpoonup}}
\def\weak {\rightharpoonup}

\def\n {\vec{n}}
\def\hh {\mathbbm{h}}
\def\phib {\bar{\phi}}
\def\mub {\bar{\mu}}
\def\sigmab {\bar{\sigma}}
\def\phiob {\bar{\phi_0}}
\def\sigmaob {\bar{\sigma_0}}
\def\Uad {\mathcal{U}_{\text{ad}}}
\def\Vad {\mathcal{V}_{\text{ad}}}
\def\phimeas {\phi_{\text{meas}}}
\def\sigmameas {\sigma_{\text{meas}}}
\def\Rcal {\mathcal{R}}
\def\HH {\mathbb{H}}
\def\VV {\mathbb{V}}
\def\WW {\mathbb{W}}
\def\AA {\mathbb{A}}
\def\vpsi {\vec{\uppsi}}
\def\Jcal {\mathcal{J}}
\def\Jtil {\widetilde{\mathcal{J}}}

\usepackage{hyperref}

\begin{document}
	
	\begin{center}
		
		{\LARGE \textbf{Identifying early tumour states in a \\ 
        Cahn--Hilliard-reaction-diffusion model}}

		\vskip0.5cm
		
		{\large\textsc{Abramo Agosti$^1$}} \\
		{\normalsize e-mail: \texttt{abramo.agosti@unipv.it}} \\
		\vskip0.35cm
		
		{\large\textsc{Elena Beretta$^2$}} \\
		{\normalsize e-mail: \texttt{eb147@nyu.edu}} \\
		\vskip0.35cm
		
		{\large\textsc{Cecilia Cavaterra$^3$}} \\
		{\normalsize e-mail: \texttt{cecilia.cavaterra@unimi.it}} \\
		\vskip0.35cm
		
		{\large \textsc{Matteo Fornoni$^1$}} \\
		{\normalsize e-mail: \texttt{matteo.fornoni01@universitadipavia.it}} \\
		\vskip0.35cm
		
		{\large\textsc{Elisabetta Rocca$^4$}} \\
		{\normalsize e-mail: \texttt{elisabetta.rocca@unipv.it}} \\
		\vskip0.35cm
		
		{\footnotesize $^1$Department of Mathematics ``F. Casorati'', University of Pavia, 27100 Pavia, Italy}
		\vskip0.1cm
		
		{\footnotesize $^2$Division of Science, New York University Abu Dhabi, Abu Dhabi, United Arab Emirates}
		\vskip0.1cm
		
		{\footnotesize $^3$Department of Mathematics ``F. Enriques'', University of Milan, 20133 Milan, Italy \\ \& IMATI-C.N.R., 27100 Pavia, Italy}
		\vskip0.1cm

        {\footnotesize$^4$Department of Mathematics ``F. Casorati'', University of Pavia \& IMATI-C.N.R., 27100 Pavia, Italy}
        \vskip0.5cm
		
	\end{center}

	\begin{abstract}\noindent
        In this paper we tackle the problem of reconstructing earlier tumour configurations starting from a single spatial measurement at a later time. 
        We describe the tumour evolution through a diffuse interface model coupling a Cahn--Hilliard-type equation for the tumour phase field to a reaction-diffusion equation for a key nutrient proportion, also accounting for chemotaxis effects.
        We stress that the ability to reconstruct earlier tumour states is crucial for calibrating the model used to predict the tumour dynamics and also to identify the areas where the tumour initially began to develop.
        However, backward-in-time inverse problems are well-known to be severely ill-posed, even for linear parabolic equations. 
        Moreover, we also face additional challenges due to the complexity of a non-linear fourth-order parabolic system. 
        Nonetheless, we can establish uniqueness by using logarithmic convexity methods under suitable a priori assumptions.
        To further address the ill-posedness of the inverse problem, we propose a Tikhonov regularisation approach that approximates the solution through a family of constrained minimisation problems. For such problems, we analytically derive the first-order necessary optimality conditions. Finally, we develop a computationally efficient numerical approximation of the optimisation problems by employing standard $C^0$-conforming first-order finite elements. We conduct numerical experiments on several pertinent test cases and observe that the proposed algorithm consistently meets expectations, delivering accurate reconstructions of the original ground truth. 
		
		\vskip3mm
		
		
		\vskip3mm
		
		
	\end{abstract}


\pagestyle{fancy}
\fancyhf{}	
\fancyhead[EL]{\thepage}
\fancyhead[ER]{\textsc{Agosti -- Beretta -- Cavaterra -- Fornoni -- Rocca}} 
\fancyhead[OL]{\textsc{Identifying early states in a Cahn--Hilliard tumour model}} 
\fancyhead[OR]{\thepage}

\renewcommand{\headrulewidth}{0pt}
\setlength{\headheight}{5mm}

\thispagestyle{empty} 

\section{Introduction}
In the last few decades, several efforts have been spent to develop mathematical models for tumour growth, in order to predict the tumour progression at its different stages of evolution by appropriately modelling the dynamics of multiple tumour cell species and their interaction with the tumour microenvironment (e.g. extracellular matrix and water). See for example the review papers \cite{Lowengrub1,Lowengrub2, BLM}. These efforts aim to create computational tools to assist clinicians in enhancing diagnosis and optimising therapy schedules. The modelling approaches fall into two main categories: discrete modelling, to describe cell-based mechanisms at the molecular scales, and continuum models, to describe the multiple tumour cell species interactions and their invasive dynamics in the microenvironment at the tissue scale.

Since tumours consist of multiple interacting phases, and due to the fact that the biomechanics of invasion of tumour cells in the microenvironment is characterised by non-sharp interfaces between the tumour and the host tissues, where tumour cells infiltrate into the parenchyma deteriorating the extracellular matrix, the natural modelling framework to describe tumour growth at the tissue level is the multiphase approach based on diffuse-interface mixture theory. In the latter approach, sharp interfaces are replaced by narrow transition layers where cell species can mix and interact due to their differential adhesive forces. The solid tumour and its microenvironment are described as a saturated medium composed of different tumour and healthy cell species and liquid phases, with possible additional phases to model growth and invasion mechanisms at different stages of evolution, like nutrients and neovascularisation (at late stages) components. The governing equations consist of mass and momentum balance equations for each phase, which, due to the diffuse interface description of interfaces between the phases, take the form of Cahn--Hilliard type system of equations. The mass and momentum equations are complemented with mass and momentum exchange terms between the phases, and appropriate constitutive laws to close the model equations. Recent theoretical and numerical studies of Cahn--Hilliard type tumour growth models, which describe tumour progression in avascular early stages, can be found e.g.~in \cite{Lowengrub3,CGH2015,FGR2015_TumGrowth,ACGAC2018,Agosti1,GLSS2016,EG2019,CGSS2023,GLS2021,GKT2023,F2024_maxreg}, while complete models which describe tumour progression in vascular late stages can be found e.g.~in \cite{Lowengrub4,Agosti2}. 
In this paper, we consider the following tumour growth model, based on the one introduced in \cite{HZO2012}:
\begin{alignat}{2}
	& \partial_t \phi - \Delta \mu
	= P(\phi) \left(\sigma + \chi (1-\phi) - \mu \right) - c(x,t) \hh (\phi)
	\qquad && \hbox{in $Q_T$,} \label{eq:phi} \\
	& \mu 
	= - \Delta \phi + F'(\phi) - \chi \sigma 
	\qquad && \hbox{in $Q_T$,} \label{eq:mu} \\
	& \partial_t \sigma - \Delta \sigma + \chi \Delta \phi 
	= - P(\phi) \left(\sigma + \chi (1-\phi) - \mu \right) + (1 - \sigma)
	\qquad && \hbox{in $Q_T$,} \label{eq:sigma} \\
	& \partial_{\n} \phi = \partial_{\n} \mu = \partial_{\n} \sigma = 0 
	\qquad && \hbox{in $\Sigma_T$,} \label{bc} \\
	& \phi(0) = \phi_0, \quad \sigma(0) = \sigma_0 
	\qquad && \hbox{in $\Omega$,} \label{ic}
\end{alignat}
where $\Omega \subset \mathbb{R}^N$, $N=2,3$, is a smooth open and bounded domain with outward unit normal $\vec{n}$, $T>0$ is the final time, $Q_T := \Omega \times (0,T)$ and $\Sigma_T = \partial \Omega \times (0,T)$.
Here, $\phi$ is a phase-field variable representing the tumour volume fraction, where $\phi \equiv 0$ is the healthy phase and $\phi \equiv 1$ is the tumour one, and satisfies a Cahn--Hilliard-type equation, being $\mu$ the corresponding chemical potential. 
Conversely, $\sigma$ represents the concentration of a key nutrient driving the tumour proliferation and satisfies a reaction-diffusion equation. 
We also include chemotaxis effects through the parameter $\chi \ge 0$.
More details about the derivation of the model and its structure, as well as precise information on all the model parameters, can be found in Section \ref{sec:model}.

The problem we want to address is the identification of the initial data $(\phi_0, \sigma_0)$ starting from a single measurement of $(\phi(T), \sigma(T))$ at the final time. 
This step is of paramount importance in the development of patient-specific tumour forecasts, which can be then leveraged to better inform clinical decision-making \cite{Lorenzo2022_review}. 
Indeed, depending on the type of cancer and the available screening methods, it may not be possible to obtain tumour measurements until the tumour has developed sufficiently.
Moreover, standard monitoring strategies may not enable the assessment of the tumour status often enough to capture its growth or its treatment response in detail.
Consequently, the possibility of reconstructing earlier states of the tumour evolution could be crucial in having better-calibrated models to be later used to predict and manage the tumour dynamics.  
Furthermore, the proposed backward problem could also be useful to locate with more precision the areas where tumours started growing. 
For instance, this is especially useful in the case of brain tumours, where the presence of cancerous cells in different areas of the brain can affect the patient in a multitude of different ways. 
Hence, locating precisely the earlier position of the tumour could give valuable suggestions to the therapy decision-making. 
For these reasons, inverse problems of this kind are becoming increasingly interesting in the mathematical oncology community, but we can still count just a few contributions up to now \cite{JBJA2019, SSMB2020, subramanian2022ensemble, BCFLR2024, BCFLR2024_num, FLS2021}.

More precisely, we call $\mathcal{R}: (\phi_0, \sigma_0) \mapsto (\phi(T), \sigma(T))$ the solution operator which associates to any initial data $(\phi_0, \sigma_0) \in \Hx1 \times \Lx2$ the value of the solution at the final time. 
Given some measurements $(\phimeas, \sigmameas)$ of the variables at the terminal time, we can formulate our inverse problem in the following way:
\begin{equation}
    \label{inv:prob}
    \hbox{find initial data $(\phi_0, \sigma_0)$ such that $(\phi(T), \sigma(T)) = (\phimeas, \sigmameas)$,}
\end{equation}
where $(\phi, \sigma)$ is the solution to \eqref{eq:phiad}--\eqref{icad}.
First of all, by recalling well-posedness results on the forward system \cite{FGR2015_TumGrowth, GY2020} and by improving some continuous dependence estimates, we show that the forward operator $\Rcal$ is Lipschitz continuous under minimal assumptions on the initial data. 
However, backward inverse problems of this kind are well-known to be severely ill-posed, even in the case of linear parabolic equations. 
Thus, to obtain uniqueness or stability results, one usually regularises the problem by assuming some additional \emph{a priori} assumptions on the initial data to be reconstructed. 
In our case, if we consider the reconstruction of more regular initial data $(\phi_0, \sigma_0) \in \Hx2 \times \Hx1$, we are able to show that the operator $\Rcal$ is injective (Theorem \ref{thm:backuniq}). 
Hence, the solution to the inverse problem is unique in this more regular class. 
To do this, we first show the existence of strong solutions to \eqref{eq:phiad}--\eqref{icad} and then prove a backward uniqueness result by employing the logarithmic convexity approach \cite{AN1967}.
We additionally mention that the high regularity of the strong solutions to our system would also allow us to prove that $\Rcal$ is actually continuously Fr\'echet differentiable and that its derivative is injective at any point (see Remark \ref{rmk:lipstab}). 
These properties would then open the possibility of getting Lipschitz stability estimates by assuming the unknown initial data to be lying in a finite-dimensional subspace \cite{B2013}.
However, we do not outline this procedure in full detail, because one of our future aims would be to prove a quantitative Lipschitz stability estimate as done in \cite{BCFLR2024}. 
Indeed, in the cited paper, the authors consider a very similar inverse problem on a prostate cancer model and, due to the specific structure of their model, they could extract more information on the stability of the backward problem. 
The same kind of analysis seems presently out of reach for our system, due to the fourth-order structure of the Cahn--Hilliard equation, but it would be also crucial in establishing convergence properties of Landweber-like iterative reconstruction algorithms \cite{HNS1995, kaltenbacher:neubauer:scherzer, BCFLR2024_num}.

In practice, to reconstruct the solution to the ill-posed inverse problem, we propose a Tikhonov regularisation approach to approximate the previous problem with a family of constrained minimisation ones.
Tikhonov regularisation is a common method in approximating solutions to inverse problems \cite{engl:hanke:neubauer, EKN1989}. 
It generally consists of replacing the ill-posed inverse problem with a one-parameter family of well-posed constrained minimisation ones. 
The key to this strategy is the addition of a regularising term in the minimisation functional. 
In our case, due to the inherent structure of the phase field system, we use the Ginzburg--Landau energy as a regularisation term for $\phi_0$ and the $\Lx 2$-norm for $\sigma_0$.
More precisely, we consider the following problem:
\begin{equation}
\label{oc:problem}
	\begin{split}
	(\phiob, \sigmaob) = \argmin_{(\phi_0, \sigma_0) \, \in \, \Uad \times \Vad} & \Bigg( \frac{\lambda_1}{2} \norm{ \phi(T) - \phimeas }^2_{\Lx 2} + \frac{\lambda_2}{2} \norm{ \sigma(T) - \sigmameas }^2_{\Lx 2} \\
	& \quad + \alpha_1 \int_{\Omega} \Big( F(\phi_0) + \mezzo |\nabla \phi_0|^2 \Big) \,\de x + \frac{\alpha_2}{2} \norm{\sigma_0}^2_{\Lx 2} \Bigg),
	\end{split}
\end{equation} 
where $( \phi(T), \sigma(T) ) = \Rcal ( \phi_0, \sigma_0 )$ and
\begin{gather*}
	\Uad \times \Vad = \left\{ \phi_0 \in H^1(\Omega) \mid 0 \le \phi_0 \le 1 \text{ a.e. in } \Omega \right\} \times \left\{ \sigma_0 \in L^2(\Omega) \mid 0 \le \sigma_0 \le 1 \text{ a.e. in } \Omega \right\} 
\end{gather*}

\noindent
Here, $\lambda_1, \lambda_2 \ge 0$ are two non-negative parameters that can be used to calibrate the relative orders of magnitude of the variables, while $\alpha_1, \alpha_2 \ge 0$ are the Tikhonov-regularisation parameters. 
We assume to be reconstructing unknown initial data $(\phi_0, \sigma_0)$ belonging to the admissible set $\Uad \times \Vad$. Hence, they are assumed to have minimal regularity for the existence of weak solutions, as well as some physical bounds representing the fact that they are respectively a phase field and a concentration of nutrient. 
Even if we can prove uniqueness in the more regular class $\Hx2 \times \Hx1$, we choose to reconstruct the initial data in the broader class introduced above. 
The motivation behind such a choice is twofold.
First of all, for application purposes, it may be useful to allow the reconstruction of less regular initial data in some realistic situations. 
Secondly, in view of the numerical applications, we wanted to avoid putting higher-order regularisation terms in the minimisation functional, because they would then impose the employment of computationally demanding discretisation schemes.
Nevertheless, regarding the constrained minimisation problem \eqref{oc:problem}, we first show the existence of minimisers (Theorem \ref{thm:control_existence}) and then characterise such minimisers through first-order necessary optimality conditions (Theorem \ref{thm:optcond}). 
The key ingredient in proving these kinds of results is the Fr\'echet differentiability of the forward map $\Rcal: \Hx1 \times \Lx2 \to \Lx2 \times \Lx2$, together with the study of the linearised and adjoint systems. 
In particular, since we just assume the initial data to be in $\Hx1 \times \Lx2$, we mention that in Theorem \ref{thm:frechet} we prove the Fr\'echet differentiability of the solution operator starting only from weak solutions to \eqref{eq:phiad}--\eqref{icad}.
To do this, we need some stricter hypotheses on the non-linearities in the system, which, however, still encompass the physically relevant cases \eqref{f:example} and \eqref{p:example}.
This fact can be seen as a new technical achievement of this paper, as in optimal control applications of this kind of model the Fr\'echet-differentiability of the forward map was always proved starting from strong solutions \cite{CGRS2017, CSS2021_secondorder, CRW2021}. 
We also emphasize that, although Gate\^aux differentiability might suffice to establish optimality conditions for the Tikhonov-regularised problem, we endeavour to prove the Fr\'echet differentiability of the forward map. 
This property is valuable as it provides additional insights into the inverse problem and is crucial in developing convergent and robust reconstruction algorithms \cite{HNS1995, EKN1989}.

We finally mention that a similar Tikhonov-regularised problem was considered in \cite{FLS2021} for the reconstruction of the initial data for a non-local version of the tumour growth model \eqref{eq:phiad}--\eqref{icad}. 
In comparison, here we not only derive the optimality conditions, but also prove some theoretical results on the inverse problem and perform a numerical study to explore the practical feasibility of the reconstruction.
Indeed, the last sections of this paper are devoted to the numerical discretisation of the model \eqref{eq:phiad}--\eqref{icad} and the Tikhonov regularised problem \eqref{oc:problem}. 
In order to design a computationally efficient numerical approximation of the optimisation problem, we employ standard $\mathcal{C}^0$-conforming first-order finite elements to discretise the forward and adjoint problem in space and a backward Euler scheme in time, together with a mesh refinement technique to better track the diffuse interface \cite{GHK2016}. 
We observe that the use of $\mathcal{C}^0$-conforming finite elements does not allow us to consider initial data and regularisation terms of class $H^2(\Omega)$ at the discrete level. This would require $\mathcal{C}^1$-conformity, which for standard finite elements entails the definition of a large number of degrees of freedom per element (see e.g. \cite{Argyris}), making the numerical approximation computationally expensive. We leave the design of discrete optimisation problems associated to the Cahn--Hilliard equation through non-standard $\mathcal{C}^1$-conforming finite elements (see e.g. \cite{Antonietti}) for future developments. 
Then, we derive the discrete version of the optimality conditions, employing a ``first discretise then optimise'' strategy (see e.g. \cite[Chapter 3]{Hinzebook}), and approximate them through a projected gradient method with a learning rate chosen according to a line search algorithm. We emphasise that in the present work, by proving Fr\'echet differentiability within the framework of weak solutions, we derive the optimality conditions for $\phi_0$ in $H^1(\Omega)$ at both the continuous and discrete levels.
Consequently, we conduct numerical experiments in some relevant test scenarios, getting highly encouraging results on the performance of the reconstruction algorithm.
We must observe that, while the tumour distribution at time $T$ can be obtained by segmentation of conventional MRI data, the distribution of nutrients, e.g. oxygen, or metabolites associated with tumour growth at time $T$ are more difficult to identify and may be inferred through non-conventional MR data, like perfusion MRI and MR spectroscopy, which are scarcely available from retrospective datasets acquired in clinical protocols. Also, the reconstruction of the nutrient distribution from the blood volume maps obtained from the perfusion MRI is difficult due to the presence of edema and blood-brain-barrier leakage induced by tumour growth. Hence, in our numerical tests, after considering a few test cases where we assume to have measurements of both $\phi(T)$ and $\sigma(T)$, we will further investigate the behaviour of the numerical algorithm (like convergence, dependence from the regularization parameter and the final time $T$) assuming to know only $\phi(T)$ (i.e.~by setting $\lambda_2=\alpha_2=0$ in \eqref{oc:problem}), in the perspective of applying our numerical algorithm to real test cases. 
In our numerical experiments, we explore the behaviour of the algorithm by varying in different ways all the following parameters: the Tikhonov regularisation parameter, the final time horizon, the position of the initial guess and the spatial dimension ($2D$ or $3D$).
We also consider the case in which the end-time measurement is affected by some noise, confirming that the Tikhonov-regularised problem acts as a regularisation method \cite{engl:hanke:neubauer}. 
In all situations, the proposed algorithm performed within the expectations, showing faithful reconstructions of the original ground truth. 
Hence, we believe this is a first promising step in the direction of applying these techniques to more sophisticated models and eventually to real medical data. 

The paper is organised as follows. 
First, in Section \ref{sec:model} we introduce our tumour growth model \eqref{eq:phi}--\eqref{ic}, by explaining in detail all the modelling choices and parameters.
In Section \ref{sec:prelim} we collect some useful preliminary notions, introduce our main hypotheses on the parameters of the model and state the well-posedness results on the forward system (Theorems \ref{thm:weaksols} and \ref{thm:contdep}). 
Section \ref{sec:uniq} is dedicated to the uniqueness result for the inverse problem, in case the unknown initial data is supposed to be more regular. 
This is done by showing first the existence of strong solutions to the forward system and then by proving a backward uniqueness result. 
Section \ref{sec:tikhonov}, instead, is devoted to the analysis of the Tikhonov-regularised minimisation problem in a broader regularity context. 
Here, we show our Fr\'echet-differentiability result and then derive the first-order optimality conditions.
Finally, in Section \ref{sec:numerics} we introduce our numerical setup, as well as our reconstruction algorithm, and in Section \ref{sec:simulations} we show representative numerical experiments to validate our approximation strategy.

\section{Mathematical model}
\label{sec:model}

In the present paper, we consider the tumour growth model introduced in \cite{HZO2012}, describing avascular tumour growth and based on a four-phase continuum mixture composed by a tumour cell phase, with volume fraction $\phi$, and a healthy cell phase, with volume fraction $h$, together with nutrient-rich and nutrient-poor extracellular water phases, with volume fractions $\sigma$ and $w$ respectively. The mixture is closed and saturated, and as a further assumption the total volume concentration of cells is assumed to be constant, i.e. $\phi+h=C$, and hence $\sigma+w=1-C$. The constants $C$ and $1-C$ can be rescaled to be unity, so that all the variables range between $0$ and $1$. This reduces to two the number of independent variables in the model, which takes the form of the following coupled system of PDEs: 
\begin{alignat}{2}
	& \partial_t \phi - M_{\phi}\Delta \mu_{\phi}
	= \gamma_{\phi}
	\qquad && \hbox{in $Q_T$,} \label{phihd1} \\
	& \mu_{\phi} 
	= - \eps^2 \Delta \phi + \Gamma F'(\phi) - \chi \sigma 
	\qquad && \hbox{in $Q_T$,} \label{muhd1} \\
	& \partial_t \sigma - M_{\sigma}\Delta \mu_{\sigma} 
	= \gamma_{\sigma}
	\qquad && \hbox{in $Q_T$,} \label{sigmahd1} \\
        & \mu_{\sigma} 
	= \frac{1}{\delta}\sigma + \chi (1-\phi) 
	\qquad && \hbox{in $Q_T$,} \label{mushd1} \\
	& \partial_{\n} \phi = \partial_{\n} \mu = \partial_{\n} \sigma = 0 
	\qquad && \hbox{in $\Sigma_T$,} \label{bchd1} \\
	& \phi(0) = \phi_0, \quad \sigma(0) = \sigma_0 
	\qquad && \hbox{in $\Omega$,} \label{ichd1}
\end{alignat}
\Accorpa\eqnhd {phihd1} {ichd1} 
where $\Omega \subset \mathbb{R}^N$, $N=2,3$, is a smooth open and bounded domain, $T>0$ is the final time, $Q_T := \Omega \times (0,T)$ and $\Sigma_T = \partial \Omega \times (0,T)$. The dynamics of the tumour cells and the nutrient-rich water component in {\eqnhd} are described by mass action laws, with associated chemical potentials given as functional derivatives of the Helmholtz free energy of the isothermal mixture,
\begin{equation}
    \label{hfe}
    E:=\int_{\Omega}\left(\Gamma F(\phi)+\frac{\eps^2}{2}|\nabla \phi|^2+\chi \sigma(1-\phi)+\frac{1}{2\delta}\sigma^2\right)\,\de x.
\end{equation}
The first term in \eqref{hfe} accounts for the free energy per unit volume due to cell-cell and cell-matrix adhesion, whereas the
second term represents the non-local intermixing that generates a surface tension between the tumour and
the host tissue, across a diffuse interface of thickness proportional to $\eps$. The third term in \eqref{hfe} accounts for the chemotactic coupling between $\phi$ and $\sigma$; indeed, for high values of $\sigma$ this term is minimised when $\phi$ grows to one. Finally, the last term is a mass contribution for the nutrient associated to a pressure term in the nutrient-rich water phase. As a diffuse interface model, $\phi$ can vary continuously between the tumour phase $\phi \approx 1$ and the healthy phase $\phi \approx 0$. 
The function $F$ is generally a double-well potential with equal minima on the pure phases, such as
\begin{equation}
\label{f:example}
    F(\phi) = \frac{1}{4} \phi^2 (\phi - 1)^2.
\end{equation}
The system, without source terms, can be seen as an $H^{-1}$-gradient flow of the free energy \eqref{hfe}. 
As one can see, the energy penalises frequent transitions on the interface and drives the evolution of $\phi$ towards the minima of $F$. 
Moreover, the chemotactic term related to the coefficient $\chi \ge 0$ acts as a cross-diffusion term, which models the fact that tumours are naturally attracted towards regions with higher concentrations of nutrients. 

A further assumption in \cite{HZO2012} is that the mass exchange terms in \eqref{phihd1} and \eqref{sigmahd1} are given as chemical kinetic rates of the form
\begin{equation}
    \label{kr}
    \gamma_{\phi}=\delta P_0P(\phi)\left(\mu_{\sigma}-\mu_{\phi}\right),\quad \gamma_{\sigma}=-\gamma_{\phi},
\end{equation}
where $P(\phi)=\max\{0,\phi\}$. Note that the model equations {\eqnhd} do not guarantee that $\phi\in [0,1]$.
With assumption \eqref{kr} the mass of $\phi+\sigma$ is constant and the dynamics is fully dissipative, i.e. the free energy of the system is decreasing along solutions. The latter property is peculiar to system \eqnhd, while in other models for tumour growth introduced in the literature (see e.g.~\cite{Lowengrub3,Agosti1,Lowengrub4,Agosti2}) the source terms break the Lyapunov-type dissipative structure of the dynamics, which makes its analytical study based on energy-type estimates more involved. 

In \cite{HZO2012} the parameter $\delta$ is taken as a small positive value, and $P_0\geq 0$, hence the leading contribution in the source term $\gamma_{\phi}$ is the proliferation term $P_0\max\{0,\phi\}\sigma$, which is linear in the nutrient concentration and null for negative values of the tumour cells concentration, while the additional terms $\delta P_0 \max\{0,\phi\}(\chi(1-\phi)-\mu_{\phi})$ are small and serve the only purpose of maintaining the dissipative structure of the system. At the same time, the leading contribution in the source term $\gamma_{\sigma}$ is $-P_0\max(0,\phi)\sigma$, modelling the consumption of nutrient by tumour cells. We observe that this picture is lacking in the description of apoptosis mechanisms, i.e. cell death, and in the inclusion of source terms for the nutrients from the healthy tissues, which both influence tumour growth dynamics at all stages. For this reason, we modify the expression of the source terms in the following way 
\begin{equation}
    \label{kr2}
    \gamma_{\phi}=\delta P_0 P(\phi)\left(\mu_{\sigma}-\mu_{\phi}\right)-c(x,t)\hh(\phi),\quad \gamma_{\sigma}=-\delta P_0 P(\phi)\left(\mu_{\sigma}-\mu_{\phi}\right)+\kappa (1 - \sigma),
\end{equation}
where $c(x,t):=D_0+D_T(x,t)$ is the sum of the apoptosis rate $D_0$ with a possible profile for a therapy death rate $D_T(x,t)$, while the term $\kappa(1-\sigma)$ describes the release of nutrients from the healthy vasculature in the nutrient-rich liquid phase at a rate $\kappa$, up to the saturation level. A possible expression for the function $P$, which avoids the introduction of an unrealistic proliferation rate when $\phi>1$, is the following
\begin{equation}
	\label{p:example}
	P(s) = \begin{cases}
		0 & \text{if } s \le 0, \\
		s & \text{if } 0 \ls s \ls 1, \\
		1 & \text{if } s \ge 1.
	\end{cases}
\end{equation}
The function $\hh(\phi)$ may also have the same expression of the function $P$.

In \eqnhd, \eqref{hfe} and \eqref{kr2} we have explicitly introduced the following physical parameters of the problem (with corresponding units in three space dimensions): the mobility parameters $M_{\phi},M_{\sigma}$ (units of $mm^2Pa^{-1}\text{day}^{-1}$), the proliferation, apoptosis and therapy rates $P_0,D_0,D_T(x,t)$ (units of $\text{day}^{-1}$), the Young modulus of the tumour tissue $\Gamma$ and the chemotactic coefficient $\chi$ (units of $Pa$), the interface thickness of the tumour cells phase $\eps$ (units of $\sqrt{Pa}\,mm$), the inverse of the pressure contribution from the nutrient-rich water phase $\delta$ (units of $Pa^{-1}$) and the release rate of nutrients from the healthy tissues $\kappa$ (units of $\text{day}^{-1}$). Note that, substituting \eqref{muhd1} in \eqref{phihd1}, we would obtain a chemotactic flux of the form $M_{\phi}\chi \Delta \sigma$, where $M_{\phi}\chi$ has units of $mm^2\text{day}^{-1}$, which are the standard units for a chemotactic coefficient. To obtain an a-dimensionalized version of \eqnhd, we introduce the nutrient penetration length 
\[
l_n:=\sqrt{\frac{M_{\sigma}\Gamma}{P_0}},
\]
and we make the change of variables
\[
\tilde{t}=tP_0, \quad \tilde{x}=\frac{x}{l_n}.
\]
We also introduce the transformed variables and parameters
\begin{align*}
&\tilde{\delta}=\delta \Gamma, \;\; \tilde{c}(x,t)=\frac{c(x,t)}{P_0}, \;\; \tilde{\kappa}=\frac{\kappa}{P_0},\\
&\tilde{\mu}=\frac{\mu}{\Gamma}, \;\; \tilde{\eps}=\frac{\eps}{\sqrt{\Gamma}l_n}, \;\; \tilde{\chi}=\frac{\chi}{\Gamma}.
\end{align*}
We then obtain the following a-dimensionalized version of \eqnhd, which we write without reporting the tilde superscript on the independent variables for ease of notation and substituting \eqref{mushd1} in \eqref{phihd1} and \eqref{sigmahd1}:
\begin{alignat}{2}
	& \partial_t \phi - \frac{M_{\phi}}{M_{\sigma}}\Delta \mu
	= \tilde{\delta} P(\phi) \left(\frac{\sigma}{\tilde{\delta}} + \tilde{\chi} (1-\phi) - \mu\right) - \tilde{c}(x,t) \hh (\phi)
	\qquad && \hbox{in $Q_T$,} \label{eq:phiad} \\
	& \mu 
	= -\tilde{\eps}^2\Delta \phi + F'(\phi) - \tilde{\chi} \sigma 
	\qquad && \hbox{in $Q_T$,} \label{eq:muad} \\
	& \partial_t \sigma - \Delta \sigma + \tilde{\chi} \Delta \phi 
	= - \tilde{\delta} P(\phi) \left(\frac{\sigma}{\tilde{\delta}} + \tilde{\chi} (1-\phi) - \mu\right) + \tilde{\kappa} (1 - \sigma)
	\qquad && \hbox{in $Q_T$,} \label{eq:sigmaad} \\
	& \partial_{\n} \phi = \partial_{\n} \mu = \partial_{\n} \sigma = 0 
	\qquad && \hbox{in $\Sigma_T$,} \label{bcad} \\
	& \phi(0) = \phi_0, \quad \sigma(0) = \sigma_0 
	\qquad && \hbox{in $\Omega$.} \label{icad}
\end{alignat}
\Accorpa\eqnhdad {eq:phiad} {icad}
Since in the following analysis all the a-dimensional combinations of parameters in \eqnhdad, except $\tilde{\chi}$, do not play a significant role, we can take them, without loss of generality and for ease of notation, as equal to one by choosing $M_{\phi}=M_{\sigma}$, $\tilde{\eps}=1$ (i.e. $\eps=\sqrt{\Gamma}l_n$), $\tilde{\kappa}=1$ (i.e. $\kappa=P_0$). On the contrary, we will keep track of the specific value of the parameter $\tilde{\chi}$ throughout the analysis, renaming it as $\chi$, without reporting the tilde superscript, for ease of notation.
With the above-mentioned choices, the tumour growth system then takes the form \eqref{eq:phi}--\eqref{ic}.  

\section{Preliminaries and weak well-posedness}
\label{sec:prelim}

We first introduce some notation that will be used throughout the paper.
Let $\Omega \subset \R^N$, $N = 2,3$, be an open and bounded domain with $\mathcal{C}^2$ boundary and outward unit normal vector $\vec{n}$.
Let $T > 0$ be the final time and denote $Q_t = \Omega \times (0, t)$ and $\Sigma_t = \partial \Omega \times (0, t)$, for any $t \in (0, T]$.
We recall the usual conventions regarding the Hilbertian triplet used in this context. 
If we define
\[ H = L^2(\Omega), \quad V = H^1(\Omega), \quad W = \{ u \in H^2(\Omega) \mid \partial_{\n} u = 0 \text{ on } \partial \Omega \}, \]
then we have the continuous and dense embeddings:
\[ W \hookrightarrow V \hookrightarrow H \cong H^* \hookrightarrow V^* \hookrightarrow W^*. \]
We denote by $\duality{\cdot, \cdot}_V$ the duality pairing between $V$ and $V^*$ and by $(\cdot, \cdot)_H$ the scalar product in $H$.
Regarding Lebesgue and Sobolev spaces, we will use the notation $\norm{\cdot}_{\Lx{p}}$ for the $\Lx p$-norm and $\norm{\cdot}_{\Wx{k,p}}$ for the $\Wx {k,p}$-norm, with $k \in \N$ and $1 \le p \le \infty$. 
Moreover, by standard elliptic regularity estimates, we equip $W$ with the equivalent norm
\[ \norm{u}^2_W := \norm{u}^2_H + \norm{\Delta u}^2_H. \]
Finally, we recall the Riesz isomorphism $\mathcal{N}: V \to V^*$:
\[ \duality{\mathcal{N}u, v}_V := \int_{\Omega} \left( \nabla u \cdot \nabla v + uv \right) \, \de x \quad \forall u,v \in V. \]
It is well-known that for $u\in W$ we have $\mathcal{N}u = - \Delta u + u \in H$ and that the restriction of $\mathcal{N}$ to $W$ is an isomorphism from $W$ to $H$. 
Additionally, by the spectral theorem, there exists a sequence of eigenvalues $0 < \lambda_1 \le \lambda_2 \le \dots$, with $\lambda_j \to +\infty$, and a family of eigenfunctions $w_j \in W$ such that $\mathcal{N}w_j = \lambda_j w_j$, which forms an orthonormal Schauder basis in $H$ and an orthogonal Schauder basis in $V$. 
In particular, $w_1$ is constant.
We also recall some useful inequalities that will be used throughout the paper: 
\begin{itemize}
	\item \emph{Gagliardo--Nirenberg inequality}. Let $\Omega \subset \R^N$ bounded Lipschitz, $m\in\N$, $1 \le r,q \le \infty$, $j\in\N$ with $0\le j \le m$ and $j/m \le \alpha \le 1$ such that
	\[ \frac{1}{p} = \frac{j}{N} + \left( \frac{1}{r} - \frac{m}{N} \right) \alpha + \frac{1-\alpha}{q}, \]
	then 
    \begin{equation}
        \label{gn:ineq}
        \norm{D^j f}_{L^p(\Omega)} \le C \, \norm{f}^\alpha_{W^{m,r}(\Omega)} \norm{f}^{1-\alpha}_{L^q(\Omega)}.
    \end{equation}
	In particular, we recall the following versions with $N=2,3$:
	\begin{equation*}
		\begin{split}
			& \norm{f}_{\Lx4} \le C \norm{f}^{1/2}_{\Hx1} \norm{f}^{1/2}_{\Lx2} \quad \text{ if } N = 2, \\
			& \norm{f}_{\Lx3} \le C \norm{f}^{1/2}_{\Hx1} \norm{f}^{1/2}_{\Lx2} \quad \text{ if } N = 3.
		\end{split}
	\end{equation*}
\end{itemize}

As mentioned in the introduction, we consider for simplicity the a-dimensionalised tumour growth model with all constants, aside from $\chi$, equal to $1$. Namely, we examine the equations \eqref{eq:phi}--\eqref{ic}, which we recall here for the reader's convenience:
\begin{alignat*}{2}
	& \partial_t \phi - \Delta \mu
	= P(\phi) \left(\sigma + \chi (1-\phi) - \mu \right) - c(x,t) \hh (\phi)
	\qquad && \hbox{in $Q_T$,} \\
	& \mu 
	= - \Delta \phi + F'(\phi) - \chi \sigma 
	\qquad && \hbox{in $Q_T$,} \\
	& \partial_t \sigma - \Delta \sigma + \chi \Delta \phi 
	= - P(\phi) \left(\sigma + \chi (1-\phi) - \mu \right) + (1 - \sigma)
	\qquad && \hbox{in $Q_T$,} \\
	& \partial_{\n} \phi = \partial_{\n} \mu = \partial_{\n} \sigma = 0 
	\qquad && \hbox{in $\Sigma_T$,} \\
	& \phi(0) = \phi_0, \quad \sigma(0) = \sigma_0 
	\qquad && \hbox{in $\Omega$.}
\end{alignat*}
We now establish several sets of hypotheses regarding the non-linearities present in the aforementioned system. Specifically, we incorporate all the physical scenarios discussed in the introduction while maintaining a level of generality.
Indeed, we make the following structural assumptions on the parameters of the system \eqref{eq:phi}--\eqref{ic}:
\begin{enumerate}[font = \bfseries, label = A\arabic*., ref = \bf{A\arabic*}]
	\item\label{ass:setting} $\chi \ge 0$, $T > 0$ and $\Omega \subset \R^N$, $N=2,3$, is an open bounded domain with $\mathcal{C}^2$ boundary. 
	\item\label{ass:fbelow} $F \in \mathcal{C}^2(\R)$ and there exist $c_1 \gs \chi^2 \ge 0$ and $c_2 \ge 0$ such that  
	\[ F(y) \ge c_1 y^2 - c_2 \quad \hbox{for any $y\in\R$.} \]
	\item\label{ass:fder} There exist $c_3 > 0$ and $c_4 \ge 0$ such that
	\[ \abs{F'(y)} \le c_3 F(y) + c_4 \quad \hbox{for any $y\in\R$.} \] 
	\item\label{ass:fconv} $F$ can be written as $F = F_0 + F_1$ for some $F_0, F_1 \in \mathcal{C}^2(\R)$. Moreover, there exist $c_0, c_0' \gs 0$, $s \in [2,6)$ and $l \ge 0$ such that 
	\[ c_0' (1 + \abs{y}^{s-2}) \le F_0''(y) \le c_0 (1 + \abs{y}^{s-2}) \quad \hbox{for any $y\in\R$,} \]
	and 
	\[ \abs{F_1''(y)} \le l \quad \hbox{for any $y\in\R$.} \]
	\item\label{ass:p} $P \in \mathcal{C}^0(\R)$ and there exist $c_5 >0$ and $q \in [1,4]$ such that 
	\[ 0 \le P(s) \le c_5 (1+\abs{y}^q) \quad \forall y \in \R. \]  
	\item\label{ass:hc} $\hh \in \mathcal{C}^0(\R) \cap L^\infty(\R)$ and $c \in L^\infty(Q_T)$. We call $\hh_\infty = \norm{\hh}_{L^\infty(\R)}$ and $c_\infty = \norm{c}_{\Lqt\infty}$. 
	\item\label{ass:initial} $\phi_0 \in V$ and $\sigma_0 \in H$.
\end{enumerate}
In the following, we will extensively use the symbol $C \gs 0$ to denote positive constants depending only on $\Omega$, $T$, the parameters of the system and all the constants introduced in the assumptions above and subsequent ones. Such constants may also change from line to line. In some cases, we will highlight specific dependencies of the constants by adding subscripts.
Moreover, we will always tacitly assume that $N = 3$, unless further stated. 
Of course, if $N = 2$ some hypotheses could be relaxed, but all the results obtained with $N = 3$ still hold.

\begin{remark}
	Observe that, by \ref{ass:fconv}, $F$ has polynomial growth up to degree $6$. Moreover, if $\phi_0 \in V$, thanks to Sobolev embeddings in three dimensions, we know that $\phi_0 \in \Lx 6$. Therefore, we can immediately conclude that $F(\phi_0) \in L^1(\Omega)$. We also point out that the bound on the polynomial growth of $F$ is needed to get the optimal regularity $\phi \in \LT 2 {\Hx3}$ (see the following Theorem \ref{thm:weaksols}), however weaker solutions can also be obtained by only assuming a general growth hypothesis as \ref{ass:fder} (see \cite[Corollary 1]{FGR2015_TumGrowth}).
\end{remark}

\noindent
We have the following result about the existence of weak solutions for \eqref{eq:phi}--\eqref{ic}.

\begin{theorem}
	\label{thm:weaksols}
	Under assumptions \emph{\ref{ass:setting}--\ref{ass:initial}}, there exists a weak solution $(\phi, \mu, \sigma)$ to \eqref{eq:phi}--\eqref{ic}, such that 
	\begin{align*}
		& \phi \in H^1(0,T;V^*) \cap L^\infty(0,T;V) \cap L^2(0,T;\Hx3), \\
		& \mu \in L^2(0,T;V), \\
		& \sigma \in H^1(0,T;V^*) \cap \mathcal{C}^0([0,T];H) \cap L^2(0,T;V), 
	\end{align*}
	which satisfies
	\[ \phi(0) = \phi_0 \quad \text{in } V \quad \text{and} \quad \sigma(0) = \sigma_0 \quad \text{in } H \]
	and the following variational formulation for a.e. $t \in (0,T)$ and for any $\zeta \in V$:
	\begin{align}
		& \duality{\phi_t, \zeta}_V + (\nabla \mu, \nabla \zeta)_H 
		= (P(\phi)(\sigma + \chi(1-\phi) - \mu), \zeta)_H - (\hh(\phi) c, \zeta)_H, \label{varform:phi} \\
		& (\mu,\zeta)_H 
		= (\nabla \phi, \nabla \zeta)_H + (F'(\phi), \zeta)_H - \chi (\sigma, \zeta)_H, \label{varform:mu} \\
		& \duality{\sigma_t, \zeta}_V + (\nabla \sigma - \chi \nabla \phi, \nabla \zeta)_H 
		= - (P(\phi)(\sigma + \chi(1-\phi) - \mu), \zeta)_H + (1 - \sigma, \zeta)_H. \label{varform:sigma}
	\end{align}
	In particular, there exists a constant $C>0$, depending only on the parameters of the model and on the data $\phi_0$ and $\sigma_0$, such that: 
	\begin{equation}
		\label{eq:weaknorms_est}
		\begin{split}
			& \norm{\phi}_{H^1(0,T;V^*) \cap L^\infty(0,T;V) \cap L^2(0,T;\Hx3)} 
			+ \norm{\mu}_{ L^2(0,T;V)} \\
			& \quad + \norm{\sigma}_{H^1(0,T;V^*) \cap L^\infty(0,T;H) \cap L^2(0,T;V)} \le C.
		\end{split}
	\end{equation}
\end{theorem}

\begin{proof}
    The argument is heavily inspired by those of \cite[Theorem 1]{FGR2015_TumGrowth} and \cite[Theorem 3.1]{GY2020}. 
    One should just highlight the differences due to the presence of the additional source terms.
    For this reason, we postpone the proof of Theorem \ref{thm:weaksols} to the Appendix.
    We just mention that, in our regularity setting, the initial value $\phi(0)$ makes sense in $V$ due to the embedding $\HT 1 {V^*} \cap \LT \infty V \hookrightarrow \mathcal{C}^0_w([0,T];V)$, where with $\mathcal{C}^0_w([0,T];V)$ we denote the space of continuous functions with respect to the weak topology in $V$. 
\end{proof}

Next, we pursue a continuous dependence result on the initial data for solutions to \eqref{eq:phi}--\eqref{ic}. 
This, in turn, will also provide uniqueness for the forward problem.
To get a strong enough continuous dependence estimate, we need to assume a bit more regularity on the functions $P$ and $\hh$. 
We need this to be able to prove the Fr\'echet differentiability of the forward map, starting only from weak solutions to \eqref{eq:phi}--\eqref{ic}.
In Section \ref{sec:tikhonov}, we will then use such properties to treat a Tikhonov-regularised version of the proposed inverse problem with minimal assumptions on the initial data.
Indeed, we now further assume that:
\begin{enumerate}[font = \bfseries, label = A\arabic*., ref = \bf{A\arabic*}]
	\setcounter{enumi}{7} 
	\item\label{ass:phinf} $P \in W^{1,\infty}(\R)$ and $\hh \in W^{1,\infty}(\R)$. Moreover, we call $P_\infty = \norm{P}_{L^\infty(\R)}$, $P'_\infty = \norm{P'}_{L^\infty(\R)}$ and $\hh'_\infty = \norm{\hh'}_{L^\infty(\R)}$.
\end{enumerate}
Then, we can prove the following result.

\begin{theorem}
	\label{thm:contdep}
	Assume hypotheses \emph{\ref{ass:setting}--\ref{ass:hc}} and \emph{\ref{ass:phinf}}. Let ${\phi_0}_1$, ${\sigma_0}_1$ and ${\phi_0}_2$, ${\sigma_0}_2$ be two sets of data satisfying \emph{\ref{ass:initial}} and let $(\phi_1, \mu_1, \sigma_1)$ and $(\phi_2, \mu_2, \sigma_2)$ two corresponding weak solutions as in Theorem \emph{\ref{thm:weaksols}}. Then, there exists a constant $K>0$, depending only on the data of the system and on the norms of $\{ ({\phi_0}_i, {\sigma_0}_i) \}_{i=1,2}$, but not on their difference, such that
	\begin{equation}
		\label{eq:contdep_estimate}
		\begin{split}
			& \norm{\phi_1 - \phi_2}_{\HT 1 {W^*} \cap \LT \infty H \cap \LT 2 W} 
			+ \norm{\mu_1 - \mu_2}_{\LT 2 H} \\
			& \quad + \norm{\sigma_1 - \sigma_2}_{\HT 1 {V^*} \cap \LT \infty H \cap \LT 2 V}
			\le K \left( \norm{{\phi_0}_1 - {\phi_0}_2}_H + \norm{{\sigma_0}_1 - {\sigma_0}_2}_H \right).
		\end{split}
	\end{equation}
\end{theorem}

\begin{proof}
    Let $\phi = \phi_1 - \phi_2$, $\mu = \mu_1 - \mu_2$, $\sigma = \sigma_1 - \sigma_2$, $\phi_0 = {\phi_0}_1 - {\phi_0}_2$ and $\sigma_0 = {\sigma_0}_1 - {\sigma_0}_2$, then, up to adding and subtracting some terms, they solve:     
	\begin{alignat}{2}
		& \partial_t \phi - \Delta \mu 
        = P(\phi_1) (\sigma - \chi \phi - \mu) \nonumber \\ 
        & \quad + (P(\phi_1) - P(\phi_2)) (\sigma_2 + \chi (1 - \phi_2) - \mu_2) - (\hh(\phi_1) - \hh(\phi_2)) c 
        \qquad && \text{in } Q_T,  \label{eq:phi2} \\
		& \mu 
        = - \Delta \phi + F'(\phi_1) - F'(\phi_2) - \chi \sigma 
        \qquad && \text{in } Q_T,  \label{eq:mu2} \\
		& \partial_t \sigma - \Delta \sigma + \chi \Delta \phi 
        = - P(\phi_1) (\sigma - \chi \phi - \mu) \nonumber \\
        & \quad - (P(\phi_1) - P(\phi_2)) (\sigma_2 + \chi (1 - \phi_2) - \mu_2) - \sigma 
        \qquad && \text{in } Q_T, \label{eq:sigma2} \\
		& \partial_{\n} \phi = \partial_{\n} \mu = \partial_{\n} \sigma = 0 \qquad && \text{on } \Sigma_T, \label{bc2} \\
		& \phi(0) = \phi_0, \quad \sigma(0) = \sigma_0 \qquad && \text{in } \Omega. \label{ic2}
	\end{alignat}
    Now, to deduce our continuous dependence estimate, we test \eqref{eq:phi2} by $\phi$, \eqref{eq:mu2} by $\Delta \phi$, \eqref{eq:sigma2} by $\sigma$ and we sum them up to obtain:
    \begin{equation}
        \label{eq:contdep:mainest}
        \begin{split}
            &\mezzo \ddt \norm{\phi}^2_H + \mezzo \ddt \norm{\sigma}^2_H + \norm{\Delta \phi}^2_H + \norm{\nabla \sigma}^2_H \\
            & \quad = (F'(\phi_1) - F'(\phi_2), \Delta \phi)_H - 2\chi (\sigma, \Delta \phi)_H + 
            (P(\phi_1)(\sigma - \chi \phi - \mu), \phi - \sigma)_H \\
            & \qquad + ((P(\phi_1) - P(\phi_2))(\sigma_2 + \chi (1 - \phi_2) - \mu_2), \phi - \sigma)_H \\
            & \qquad - ((\hh(\phi_1) - \hh(\phi_2)) \, c, \phi)_H - (\sigma, \sigma)_H.
        \end{split}
    \end{equation}
    Hence, we proceed with estimating all the terms on the right-hand side one by one. 
    First, by using Cauchy--Schwarz and Young's inequalities, we see that
    \begin{align*}
        (F'(\phi_1) - F'(\phi_2), \Delta \phi)_H 
        & \le \frac{1}{8} \norm{\Delta \phi}^2_H + C \norm{F'(\phi_1) - F'(\phi_2)}^2_H \\
        & \le \frac{1}{8} \norm{\Delta \phi}^2_H + C \norm{ ( \sup_{z \in [\phi_1, \phi_2]} F''(z) ) \,\phi }^2_H \\
        & \le \frac{1}{8} \norm{\Delta \phi}^2_H + C \norm{ \sup_{z \in [\phi_1, \phi_2]} F''(z) }^2_{\Lx\infty} \norm{\phi }^2_H \\
        & \le \frac{1}{8} \norm{\Delta \phi}^2_H + C \underbrace{\left( 1 + \norm{\phi_1}^{s-2}_{\Lx\infty} + \norm{\phi_2}^{s-2}_{\Lx\infty} \right)^2}_{:= g(t)} \norm{\phi }^2_H,
    \end{align*}
    where we used Lagrange's Theorem, since $F \in \mathcal{C}^2$ by \ref{ass:fbelow}, and the growth bound given by hypothesis \ref{ass:fconv} for $s \in [2,6)$.
    Now, in order to use Gronwall's lemma later, we need the function $g(t)$ above to be uniformly bounded in $\Lt1$. 
    Indeed, observe that, by Gagliardo--Nirenberg's interpolation inequality \eqref{gn:ineq} with $p=\infty$, $j=0$, $N=3$, $r=2$, $m=3$, $\alpha = 1/4$ and $q=6$, it follows that 
    \begin{equation}
        \label{eq:emb_l8linf}
        \LT \infty V \cap \LT 2 {\Hx3} \hookrightarrow \LT 8 {\Lx\infty}.
    \end{equation}
    Then, since by \eqref{eq:weaknorms_est} and \eqref{eq:emb_l8linf} $\phi_1$ and $\phi_2$ are bounded in $\LT 8 {\Lx\infty}$, we have that $g(t) \in \Lt1$ if $2(s-2) \le 8$, which is guaranteed by hypothesis \ref{ass:fconv}.
    Next, we can easily estimate the second and the last term of \eqref{eq:contdep:mainest} by means of Cauchy--Schwarz and Young's inequalities:
    \[ -2\chi (\sigma, \Delta \phi)_H - (\sigma, \sigma)_H \le \frac{1}{8} \norm{\Delta \phi}^2_H + C \norm{\sigma}^2_H. \]
    Observe also that, by proceeding similarly as above, we can say that 
    \begin{equation}
        \label{eq:contdep:mu_norm}
        \norm{\mu}_H \le \norm{\Delta \phi}_H + \norm{F'(\phi_1) - F'(\phi_2)}_H + \chi \norm{\sigma}_H \le \norm{\Delta \phi}_H + \sqrt{g(t)} \norm{\phi}_H + \chi \norm{\sigma}_H. 
    \end{equation}
    Hence, by using \eqref{eq:contdep:mu_norm} and Cauchy--Schwarz and Young's inequalities, we can estimate the third term on the right-hand side of \eqref{eq:contdep:mainest} as follows:
    \begin{align*}
        & (P(\phi_1)(\sigma - \chi \phi - \mu), \phi - \sigma)_H 
        \le P_\infty \left( \norm{\sigma}_H + \chi \norm{\phi}_H + \norm{\mu}_H \right) \left( \norm{\phi}_H + \norm{\sigma}_H \right) \\
        & \quad \le C \norm{\sigma}^2_H + C \norm{\phi}^2_H + C \left( \norm{\Delta \phi}_H + \sqrt{g(t)} \norm{\phi}_H + \chi \norm{\sigma}_H \right) \left( \norm{\phi}_H + \norm{\sigma}_H \right) \\
        & \quad \le \frac{1}{8} \norm{\Delta \phi}^2_H + C \left( 1 + \sqrt{g(t)} + g(t) \right) \norm{\phi}^2_H + C \norm{\sigma}^2_H,
    \end{align*}
    where $\sqrt{g(t)} + g(t)$ is still uniformly bounded in $\Lt1$.
    Next, by using hypothesis \ref{ass:phinf}, Sobolev embeddings and H\"older and Young's inequalities, we can estimate also the fourth term:
    \begin{align*}
        & ((P(\phi_1) - P(\phi_2)(\sigma_2 + \chi (1 - \phi_2) - \mu_2), \phi - \sigma)_H \\
        & \quad \le \norm{P(\phi_1) - P(\phi_2)}_H \norm{\sigma_2 + \chi (1 - \phi_2) - \mu_2}_{\Lx4} \norm{\phi - \sigma}_{\Lx4} \\
        & \quad \le C P'_\infty \norm{\phi}_H \norm{\sigma_2 + \chi (1 - \phi_2) - \mu_2}_V \left( \norm{\phi}_V + \norm{\sigma}_V \right) \\
        & \quad \le C \norm{\phi}_H \norm{\sigma_2 + \chi (1 - \phi_2) - \mu_2}_V \left( \norm{\phi}_W + \norm{\sigma}_V \right) \\
        & \quad \le C \norm{\phi}_H \norm{\sigma_2 + \chi (1 - \phi_2) - \mu_2}_V \left( \norm{\phi}_H + \norm{\Delta \phi}_H + \norm{\sigma}_V \right) \\
        & \quad \le \frac{1}{8} \norm{\Delta \phi}^2_H + \mezzo \norm{\sigma}^2_V + C \left( 1 + \norm{\sigma_2 + \chi (1 - \phi_2) - \mu_2}^2_V \right) \norm{\phi}^2_H, 
    \end{align*}
    where $\norm{\sigma_2 + \chi (1 - \phi_2) - \mu_2}^2_V \in \Lt1$ by \eqref{eq:weaknorms_est}.
    Finally, by using again hypothesis \ref{ass:phinf}, we also have that
    \[
        ((\hh(\phi_1) - \hh(\phi_2)) c, \phi)_H \le \hh'_\infty c_\infty \norm{\phi}^2_H \le C \norm{\phi}^2_H.
    \]
    Therefore, by putting all together and integrating on $(0,t)$ for any $t \in (0,T)$, we find that 
    \begin{align*}
        & \mezzo \norm{\phi}^2_H + \mezzo \norm{\sigma}^2_H + \mezzo \norm{\Delta \phi}^2_H + \mezzo \norm{\nabla \sigma}^2_H \\
        & \quad \le  \mezzo \norm{\phi_0}^2_H + \mezzo \norm{\sigma_0}^2_H + C \int_0^T \norm{\sigma}^2_H \, \de t \\
        & \qquad + C \int_0^T \left( 1 + \sqrt{g(t)} + g(t) + \norm{\sigma_2 + \chi (1 - \phi_2) - \mu_2}^2_V \right) \norm{\phi}^2_H \, \de t,
    \end{align*}
    which, by Gronwall's lemma, eventually implies that
    \begin{align*}
        \norm{\phi}^2_{\LT \infty H \cap \LT 2 W} + \norm{\sigma}^2_{\LT \infty H \cap \LT 2 V} \le C \left( \norm{{\phi_0}_1 - {\phi_0}_2}^2_H + \norm{{\sigma_0}_1 - {\sigma_0}_2}^2_H \right).
    \end{align*}
    Moreover, starting from \eqref{eq:contdep:mu_norm} and integrating in time, we can also easily see that
    \begin{align*}
        & \norm{\mu}^2_{\LT 2 H} \le \norm{\Delta \phi}^2_{\LT 2 H} + \norm{\phi}^2_{\LT \infty H} \int_0^T g(t) \, \de t + \chi^2 \norm{\sigma}^2_{\LT 2 H} \\
        & \quad \le C \left( \norm{{\phi_0}_1 - {\phi_0}_2}^2_H + \norm{{\sigma_0}_1 - {\sigma_0}_2}^2_H \right),
    \end{align*}
    since $g \in \Lt1$.
    Finally, by comparison in \eqref{eq:phi2} and \eqref{eq:sigma2}, it is also straightforward to see that 
    \begin{align*}
        \norm{\phi}^2_{\HT 1 {W^*}} + \norm{\sigma}^2_{\HT 1 {V^*}} 
        \le C \left( \norm{{\phi_0}_1 - {\phi_0}_2}^2_H + \norm{{\sigma_0}_1 - {\sigma_0}_2}^2_H \right).
    \end{align*}
    This concludes the proof of Theorem \ref{thm:contdep}.
\end{proof}

\begin{remark}
    Since we have the embeddings $\HT 1 {W^*} \cap \LT 2 W \hookrightarrow \C 0 H$ and $\HT 1 {V^*} \cap \LT 2 V \hookrightarrow \C 0 H$, the continuous dependence estimate \eqref{eq:contdep_estimate} also implies that 
    \[ 
        \norm{\phi_1(T) - \phi_2(T)}_H + \norm{\sigma_1(T) - \sigma_2(T)}_H 
        \le K \left( \norm{{\phi_0}_1 - {\phi_0}_2}_V + \norm{{\sigma_0}_1 - {\sigma_0}_2}_H \right).
    \]
    This means that the forward operator $\Rcal : V \times H \to H \times H$ such that $\Rcal(\phi_0, \sigma_0) = (\phi(T), \sigma(T))$ is well-defined and moreover it is Lipschitz continuous. 
    Observe that, in estimate \eqref{eq:contdep_estimate}, the norm of the difference of the initial data ${\phi_0}_1$ and ${\phi_0}_2$ is taken in $H$, however, to get such estimate, we needed to assume that ${\phi_0}_1, {\phi_0}_2 \in V$. 
    For this reason, we can only say that $\Rcal$ is Lipschitz continuous from $V \times H$ and not from $H \times H$.
\end{remark}

\section{Uniqueness for the inverse problem}
\label{sec:uniq}

In this section, we prove uniqueness results for the inverse problem of identifying the initial data by setting the value of the solution at the final time $T$. 
Hence, we aim to prove the injectivity of the forward operator $\Rcal : V \times H \to H \times H$ such that $\Rcal(\phi_0, \sigma_0) = (\phi(T), \sigma(T))$, which is equivalent to a backward uniqueness property for the system \eqref{eq:phi}--\eqref{ic}.
Unfortunately, to prove such a result, we need to assume more regularity on the initial data, namely that $\phi_0 \in W$ and $\sigma_0 \in V$.
We stress that this is a common procedure in the analysis of inverse problems and it is linked to our chosen proof strategy, involving logarithmic convexity methods \cite{AN1967}.
Indeed, we first obtain more regular solutions to \eqref{eq:phi}--\eqref{ic} in Theorem \ref{thm:strongsols} and then prove injectivity of the operator $\Rcal : W \times V \to H \times H$ in Theorem \ref{thm:backuniq}.
Finally, as a byproduct of the backward uniqueness result, we get a stronger continuous dependence estimate, which would then pave the way to a Lipschitz stability estimate for the inverse problem in finite-dimensional subspaces through standard results \cite{B2013}. 
We discuss this possibility in Remark \ref{rmk:lipstab} at the end of this section.
Such results further cement the fact that the inverse problem is well-posed, at least in more regular spaces, thus it is reasonable to try and approximate its solutions through a Tikhonov regularisation procedure, as done in the next Section \ref{sec:tikhonov}.

\subsection{Strong solutions}

In this subsection, we prove higher regularity results for the solutions to \eqref{eq:phi}--\eqref{ic}.
To do this, we further assume the following:
\begin{enumerate}[font = \bfseries, label = B\arabic*., ref = \bf{B\arabic*}]
	\item\label{ass:fc4} $F \in \mathcal{C}^4(\R)$.
    \item\label{ass:pc1} $P, \hh \in \mathcal{C}^1(\R)$ and the exponent $q$ of hypothesis \ref{ass:p} is such that $q \in [1,2]$.
	\item\label{ass:initial2} $\phi_0 \in W$ and $\sigma_0 \in V$.
\end{enumerate}
Then, we have the following result about strong solutions to our tumour growth model.

\begin{theorem}
	\label{thm:strongsols}
	Under assumptions \emph{\ref{ass:setting}--\ref{ass:hc}} and \emph{\ref{ass:fc4}--\ref{ass:initial2}}, a solution $(\phi, \mu, \sigma)$ to \eqref{eq:phi}--\eqref{ic} enjoys the following higher regularities:
	\begin{align*}
		& \phi \in H^1(0,T;H) \cap \mathcal{C}^0([0,T];W) \cap L^2(0,T;\Hx 4), \\
		& \mu \in L^2(0,T;W), \\
		& \sigma \in H^1(0,T;H) \cap \mathcal{C}^0([0,T];V) \cap L^2(0,T;W), 
	\end{align*}
	In particular, there exists a constant $C>0$, depending only on the parameters of the model and on the data $\phi_0$ and $\sigma_0$, such that: 
	\begin{equation}
		\label{eq:strongnorms_est}
		\begin{split}
			& \norm{\phi}_{H^1(0,T;H) \cap L^\infty(0,T;W) \cap L^2(0,T;\Hx4)} 
			+ \norm{\mu}_{ L^2(0,T;W)} \\
			& \quad + \norm{\sigma}_{H^1(0,T;H) \cap L^\infty(0,T;V) \cap L^2(0,T;W)} \le C.
		\end{split}
	\end{equation}
\end{theorem}

\begin{proof}
    We proceed again by formal estimates; to be rigorous, one should go back to a Galerkin discretisation scheme, as hinted in the proof of Theorem \ref{thm:weaksols}. 
    First, we consider equation \eqref{eq:sigma} and observe that it can be seen as
    \begin{alignat*}{2}
        & \partial_t \sigma - \Delta \sigma 
        = - \chi \Delta \phi - P(\phi) (\sigma + \chi (1 - \phi) - \mu) + (1 - \sigma)
        \quad && \hbox{in $Q_T$,} \\
        & \partial_{\n} \sigma = 0 
        \quad && \hbox{in $\Sigma_T$,} \\
        & \sigma(0) = \sigma_0 
        \quad && \hbox{in $\Omega$.}
    \end{alignat*}
    In particular, due to the regularity of the weak solution given by \eqref{eq:weaknorms_est}, we can easily infer that the right-hand side is uniformly bounded in $\LT 2 H$. 
    Indeed, the only non-trivial term to check is the reaction one, but thanks to the fact that $P(\phi)$ can grow up to an exponent $q \le 2$ by \ref{ass:pc1}, it follows that
    \begin{align*}
        \norm{P(\phi) (\sigma + \chi (1 - \phi) - \mu)}_H 
        & \le \norm{P(\phi)}_{\Lx3} \norm{\sigma + \chi (1 - \phi) - \mu}_{\Lx6} \\
        & \le C \underbrace{\left( \int_\Omega c_5 (1 + \abs{\phi}^q)^3 \, \de x \right)^{1/3}}_{\in \, \Lt\infty} \norm{\sigma + \chi (1 - \phi) - \mu}_V \in \Lt2, 
    \end{align*}
    since $\phi \in \LT \infty {\Lx6}$ and $\sigma + \chi (1 - \phi) - \mu  \in \LT 2 V$ by Sobolev embeddings and \eqref{eq:weaknorms_est}.
    Then, since $\sigma_0 \in V$ and the right-hand side is bounded in $\LT 2 H$, by standard parabolic regularity theory, we can immediately infer that 
    \begin{equation}
        \label{eq:strong:sigmareg}
        \norm{\sigma}_{\HT 1 H \cap \LT \infty V \cap \LT 2 W} \le C, 
    \end{equation}
    for some constant $C > 0$ depending only on the initial data and the parameters of the system.

    Next, for the second estimate, we test \eqref{eq:phi} by $\partial_t \phi$, \eqref{eq:mu} by $\Delta \partial_t \phi$ and we sum them up to obtain, after cancellations, that 
    \begin{equation}
        \label{eq:strong:secondestimate}
        \begin{split}
            & \norm{\phi_t}^2_H + \mezzo \ddt \norm{\Delta \phi}^2_H \\
            & \quad = (P(\phi)(\sigma + \chi (1 - \phi) - \mu), \phi_t)_H - (c \hh(\phi), \phi_t)_H - \chi (\sigma, \Delta \phi_t)_H + (F'(\phi), \Delta \phi_t)_H,
        \end{split}
    \end{equation}
    where we used the notation $\phi_t = \partial_t \phi$ for simplicity. 
    We can now estimate the first two terms on the right-hand side of \eqref{eq:strong:secondestimate} by means of H\"older and Young's inequalities, indeed:
    \begin{align*}
        & (P(\phi)(\sigma + \chi (1 - \phi) - \mu), \phi_t)_H - (c \hh(\phi), \phi_t)_H \\
        & \quad \le \norm{P(\phi)}_{\Lx3} \norm{\sigma + \chi (1 - \phi) - \mu}_{\Lx6} \norm{\phi_t}_H + c_\infty \hh_\infty \int_\Omega \abs{\phi_t} \, \de x \\
        & \quad \le \frac{1}{4} \norm{\phi_t}^2_H + C \norm{P(\phi)}^2_{\Lx 3} \norm{\sigma + \chi (1 - \phi) - \mu}^2_V + C,
    \end{align*}
    where $\norm{P(\phi)}^2_{\Lx3}$ is again uniformly bounded in $\Lt\infty$ by \eqref{eq:weaknorms_est} and \ref{ass:pc1}, hence we can bound it by a constant (depending on $T$). 
    For the last two terms, instead, we use Leibinz's rule to exchange time-derivatives and then H\"older and Young's inequalities to see that
    \begin{align*}
        & - \chi (\sigma, \Delta \phi_t)_H + (F'(\phi), \Delta \phi_t)_H \\
        & \quad = \ddt \left( - \chi (\sigma, \Delta \phi)_H + (F'(\phi), \Delta \phi)_H \right) + \chi (\sigma_t, \Delta \phi)_H - (F''(\phi) \phi_t, \Delta \phi)_H \\
        & \quad \le \ddt \left( - \chi (\sigma, \Delta \phi)_H + (F'(\phi), \Delta \phi)_H \right) + \frac{1}{4} \norm{\phi_t}^2_H + C \norm{\sigma_t}^2_H + C \left( 1 + \norm{F''(\phi)}^2_{\Lx\infty} \right) \norm{\Delta \phi}^2_H, 
    \end{align*}
    where we note that $\norm{F''(\phi)}^2_{\Lx\infty} \in \Lt 1$ by \eqref{eq:emb_l8linf} and \ref{ass:fconv}.
    Hence, starting from \eqref{eq:strong:secondestimate}, we deduced that
    \begin{align*}
        & \mezzo \norm{\phi_t} + \ddt \overbrace{\left( \mezzo \norm{\Delta \phi}^2_H + \chi (\sigma, \Delta \phi)_H - (F'(\phi), \Delta \phi)_H \right)}^{:= g(t)} \\
        & \quad \le C \norm{\sigma_t}^2_H + C \left( 1 + \norm{F''(\phi)}^2_{\Lx\infty} \right) \norm{\Delta \phi}^2_H + C \norm{\sigma + \chi (1 - \phi) - \mu}^2_V + C.
    \end{align*}
    Moreover, by integration by parts and Cauchy--Schwarz and Young's inequalities, we can see that
    \begin{align*}
        g(t) & = \mezzo \norm{\Delta \phi}^2_H + \chi (\sigma, \Delta \phi)_H - (F'(\phi), \Delta \phi)_H \\
        & \ge \mezzo \norm{\Delta \phi}^2_H - \frac{1}{4} \norm{\Delta \phi}^2_H - C \norm{\sigma}^2_H + (F''(\phi) \nabla \phi, \nabla \phi)_H \\
        & \ge \frac{1}{4} \norm{\Delta \phi}^2_H - C \norm{\sigma}^2_H - C \norm{\nabla \phi}^2_H,
    \end{align*}
    where we used the fact that $F''(y) \ge - C$ for any $y \in \R$ for some constant $C \gs 0$, which easily follows by hypothesis \ref{ass:fconv}. 
    Additionally, since $\phi_0 \in W \hookrightarrow \Lx\infty$ and $F \in \mathcal{C}^2$, we also have that
    \[  
        g(0) \le C \norm{\phi_0}^2_W + C \norm{\sigma_0}^2_H + \norm{F''(\phi_0)}_{\Lx\infty} \norm{\phi_0}^2_V \le C \left( \norm{\phi_0}^2_W + \norm{\sigma_0}^2_H \right).  
    \]
    By putting all together and integrating on $(0,t)$, for any $t \in (0,T)$, we finally arrive at
    \begin{align*}
        & \mezzo \int_0^t \norm{\phi_t}^2_H \, \de s + \frac{1}{4} \norm{\Delta \phi}^2_H \\
        & \quad \le C \left( \norm{\phi_0}^2_W + \norm{\sigma_0}^2_H \right) + C \int_0^T \left( 1 + \norm{F''(\phi)}^2_{\Lx\infty} \right) \norm{\Delta \phi}^2_H \, \de t + C \int_0^T \norm{\sigma_t}^2_H \, \de t \\
        & \qquad + \int_0^T \norm{\sigma + \chi (1 - \phi) - \mu}^2_V \, \de t + C \norm{\sigma}^2_{\LT \infty H} + C \norm{\nabla \phi}^2_{\LT \infty H}.
    \end{align*}
    Therefore, by Gronwall's inequality, as well as, \eqref{eq:weaknorms_est} and \eqref{eq:strong:sigmareg}, we conclude that 
    \begin{equation}
        \label{eq:strong:phih1h}
        \norm{\phi}_{\HT 1 H \cap \LT \infty W} \le C,
    \end{equation}
    for some constant $C$, depending only on the initial data and on the parameters.
    Moreover, by comparison in \eqref{eq:phi}, it is now easy to see that also 
    \begin{equation}
        \label{eq:strong:mul2w}
        \norm{\mu}_{\LT 2 W} \le C.
    \end{equation}
    In particular, since by Sobolev embeddings $W \hookrightarrow \Lx\infty$, from \eqref{eq:strong:phih1h} we can now infer that
    \begin{equation}
        \label{eq:strong:philinfty}
        \norm{\phi}_{\Lqt\infty} \le C,
    \end{equation}
    hence, given that $F \in \mathcal{C}^4(\R)$, we also have that
    \[
        \norm{F^{(i)}(\phi)}_{\Lqt\infty} \le C \quad \hbox{for any $i=1,\dots,4$.}
    \]
    Finally, we formally take the laplacian of \eqref{eq:mu} and then estimate the norm of the bilaplacian of $\phi$ in $H$ as follows:
    \begin{align*}
        \norm{\Delta^2 \phi}_H & \le \norm{\Delta \mu}_H + \norm{\Delta F'(\phi)}_H + \chi \norm{\Delta \sigma} \\
        & \le \norm{\Delta \mu}_H + \norm{F''(\phi) \Delta \phi}_H + \norm{F'''(\phi) \nabla \phi \cdot \nabla \phi}_H + \chi \norm{\Delta \sigma}_H \\
        & \le \norm{\Delta \mu}_H + \norm{F''(\phi)}_{\Lqt\infty} \norm{\Delta \phi}_H + \norm{F'''(\phi)}_{\Lqt\infty} \norm{\nabla \phi}^2_{\Lx4} + \chi \norm{\Delta \sigma}_H \\
        & \le \norm{\Delta \mu}_H + C \norm{\Delta \phi}_H + C \norm{\phi}^2_{\Hx2} + C \norm{\Delta \sigma}_H.
    \end{align*}
    Then, we observe that the right-hand side is uniformly bounded in $\Lt\infty$, due to \eqref{eq:strong:sigmareg}, \eqref{eq:strong:mul2w} and \eqref{eq:strong:phih1h}, therefore, by elliptic regularity theory, we deduce that
    \begin{equation}
        \label{eq:strong:phil2h4}
        \norm{\phi}_{\LT 2 {\Hx4}} \le C.
    \end{equation}
    This concludes the proof of Theorem \ref{thm:strongsols}. 
\end{proof}

\begin{remark}
    \label{rmk:fderivatives}
    By Sobolev embeddings, one can easily see that \eqref{eq:strongnorms_est} implies that 
    \[ \norm{\phi}_{\Cqt0} \le C, \]
    for some $C \gs 0$, depending only on the parameter of the system. 
    Consequently, given that $F \in \mathcal{C}^4(\R)$ and $P, \hh \in \mathcal{C}^1(\R)$ by \ref{ass:fc4}-\ref{ass:pc1}, we also infer that
    \begin{gather*}
        \norm{F^{(i)}(\phi)}_{\Cqt0} \le C \quad \hbox{for any $i=1,\dots,4$,} \\
        \norm{P^{(i)}(\phi)}_{\Cqt0} \le C \quad \hbox{for any $i=0,1$,} \\
        \norm{\hh^{(i)}(\phi)}_{\Cqt0} \le C \quad \hbox{for any $i=0,1$.}
    \end{gather*}
\end{remark}

\begin{remark}
    We just mention that the strong solution guaranteed by Theorem \ref{thm:strongsols} is unique. Indeed, even if hypothesis \ref{ass:phinf} was not assumed here, the same argument of Theorem \ref{thm:contdep} can be repeated since $\phi$ is now uniformly bounded and the estimates in Remark \ref{rmk:fderivatives} hold.
\end{remark}

\subsection{Backward uniqueness}

We can now state and prove our result about backward uniqueness for \eqref{eq:phi}--\eqref{ic}.
However, we need the following additional hypothesis:
\begin{enumerate}[font = \bfseries, label = B\arabic*., ref = \bf{B\arabic*}]
	\setcounter{enumi}{3} 
	\item\label{ass:chi} The chemotaxis coefficient $\chi$ is such that $\chi^2 \ls 2$.
\end{enumerate}
Hypothesis \ref{ass:chi} is needed to make the leading differential operator of \eqref{eq:phi}--\eqref{ic} uniformly elliptic, even in the presence of chemotaxis. 
Such a condition is a cornerstone of our proof strategy, relying on the logarithmic convexity method \cite{AN1967}.
We later comment on the biological feasibility of this hypothesis in Remark \ref{rmk:chemotaxis}
Then, under this additional assumption, we have the following result.

\begin{theorem}
	\label{thm:backuniq}
	Assume hypotheses \emph{\ref{ass:setting}--\ref{ass:hc}}, \emph{\ref{ass:fc4}--\ref{ass:pc1}} and \emph{\ref{ass:chi}}. 
    Let $(\phi_1, \sigma_1)$ and $(\phi_2, \sigma_2)$ be two solutions of \eqref{eq:phi}--\eqref{ic} corresponding to two pairs of initial data $(\phi_0^i, \sigma_0^i)$ for $i=1,2$, satisfying hypothesis \ref{ass:initial2}.
	
	If $(\phi_1, \sigma_1)(T) = (\phi_2, \sigma_2)(T)$, then $(\phi_1, \sigma_1)(t) = (\phi_2, \sigma_2)(t)$ for any $t \in [0,T]$. In particular, $(\phi_0^1, \sigma_0^1) = (\phi_0^2, \sigma_0^2)$ in $W \times V$.
\end{theorem}

\begin{proof}
    We adapt the arguments of \cite[Lemmas 6.1 and 6.2]{temam}, which are based on the logarithmic convexity method by Agmon and Nirenberg \cite{AN1967}. 
    As in the proof of Theorem \ref{thm:contdep} we let $\phi = \phi_1 - \phi_2$, $\mu = \mu_1 - \mu_2$, $\sigma = \sigma_1 - \sigma_2$, $\phi_0 = {\phi_0}_1 - {\phi_0}_2$ and $\sigma_0 = {\sigma_0}_1 - {\sigma_0}_2$.
    Then, we rewrite the system \eqref{eq:phi2}--\eqref{ic2} solved by the difference of two solutions in the following way:
    \begin{alignat}{2}
        & \partial_t \phi + \Delta^2 \phi + \chi \Delta \sigma = f_\phi 
        \quad && \hbox{in $Q_T$,} \label{eq:phi3} \\
        & \partial_t \sigma - \Delta \sigma + \chi \Delta \phi = f_\sigma 
        \quad && \hbox{in $Q_T$,} \label{eq:sigma3}\\
        & \partial_{\n} \phi = \partial_{\n} \Delta \phi = \partial_{\n} \sigma = 0 
        \quad && \hbox{on $\Sigma_T$,} \label{bc3}\\
        & \phi(0) = \phi_0, \quad \sigma(0) = \sigma_0 
        \quad && \hbox{in $\Omega$,} \label{ic3}
    \end{alignat}
    where the right-hand sides are
    \begin{align*}
        & f_\phi = \Delta (F'(\phi_1)) - \Delta (F'(\phi_2)) + P(\phi_1)(\sigma - \chi \phi + \Delta \phi - (F'(\phi_1) - F'(\phi_2)) + \chi \sigma) \\
        & \qquad + (P(\phi_1) - P(\phi_2))(\sigma_2 + \chi (1 - \phi_2) + \Delta \phi_2 - F'(\phi_2) + \chi \sigma_2) - c (\hh(\phi_1) - \hh(\phi_2)), \\
        & f_\sigma = - P(\phi_1)(\sigma - \chi \phi + \Delta \phi - (F'(\phi_1) - F'(\phi_2)) + \chi \sigma) \\
        & \qquad - (P(\phi_1) - P(\phi_2))(\sigma_2 + \chi (1 - \phi_2) + \Delta \phi_2 - F'(\phi_2) + \chi \sigma_2) - \sigma.
    \end{align*}
    In particular, we note that boundary conditions \eqref{bc2} are equivalent to those in \eqref{bc3}.
    We now frame system \eqref{eq:phi3}--\eqref{ic3} in the abstract setting of \cite[Section 6.1]{temam}. 
    To do this, we introduce the Hilbert spaces
    \[
        \HH := H \times H, \quad \VV := W \times V, \quad \WW := \{ u \in \Hx4 \mid \partial_{\n} u = \partial_{\n} \Delta u = 0 \} \times W,
    \]
    as well as the linear self-adjoint unbounded operator $\AA : D(\AA) = \WW \subset \HH \to \HH$ such that 
    \begin{equation*}
        \AA \begin{pmatrix} u \\ v \end{pmatrix} 
        : = \begin{pmatrix} \Delta^2 + \mathrm{Id} & \chi \Delta \\ \chi \Delta & - \Delta + \mathrm{Id} \end{pmatrix} \begin{pmatrix} u \\ v \end{pmatrix} 
        = \begin{pmatrix} \Delta^2 u + u + \chi \Delta v \\ - \Delta v + v + \chi \Delta u \end{pmatrix},
    \end{equation*}
    where $\mathrm{Id}: H \to H$ is the identity operator.
    We now check that $\AA$ is also positive definite if \ref{ass:chi} holds. Indeed, by using integration by parts, boundary conditions, as well as Sobolev embeddings and Cauchy--Schwarz and Young's inequalities, we see that
    \begin{align}
    \label{eqn:Apd}
        \notag \duality*{\AA \binom{u}{v}, \binom{u}{v}} & = \norm{\Delta u}^2_H + \norm{u}^2_H + \norm{\nabla v}^2_H + \norm{v}^2_H - 2 \chi (\nabla u, \nabla v)_H \\
        & \notag \ge \norm{u}^2_W + \norm{v}^2_V - (1 - \delta) \norm{\nabla v}^2_H - \frac{\chi^2}{1 - \delta} \norm{\nabla u}^2_H \\
        & \notag \ge \norm{u}^2_W + \delta \norm{v}^2_V - \frac{\chi^2}{2(1 - \delta)} \norm{u}^2_W \\
        & \ge \underbrace{\left( 1 - \frac{\chi^2}{2(1 - \delta)} \right)}_{:= \gamma} \norm{u}^2_W + \delta \norm{v}^2_V,
    \end{align}
    for any $0 \ls \delta \ls 1$. 
    Then, by \ref{ass:chi}, we can choose $0 \ls \delta \ls 1$ such that also $\gamma \gs 0$, which implies the positivity of $\AA$. 
    In particular, to justify the positivity of $\AA$, we used the inequality
    \[ 
        \norm{\nabla u}^2_H \le \mezzo \norm{\Delta u}^2_H + \mezzo \norm{u}^2_H = \mezzo \norm{u}^2_W,
    \]
    which holds for any $u \in W$ and motivates the choice of $\chi^2 < 2$ in hypothesis \ref{ass:chi}.
    Consequently, the operator $\AA^{1/2}$ is well-defined and one can easily check that $D(\AA^{1/2}) = \VV$ and that the natural norm in $D(\AA^{1/2})$ is equivalent to the one of $\VV$.
    For the rest of this proof, we redefine the scalar product and the corresponding norm in $\VV$ as the equivalent ones in terms of the operator $\AA$, namely
    \[
        (\vec{u},\vec{v})_{\VV} := (\AA^{1/2} \vec{u}, \AA^{1/2} \vec{v})_{\HH}, \quad \norm{\vec{u}}^2_{\VV} = (\AA^{1/2} \vec{u}, \AA^{1/2} \vec{u})_{\HH}, 
    \]
    for any $\vec{u}, \vec{v} \in \VV$.
    By what we have just shown above, it is clear that there exist two constants $C_1, C_2 > 0$ such that 
    \[
        C_1 \left( \norm{u_1}^2_W + \norm{u_2}^2_V \right) \le \norm{\vec{u}}^2_{\VV} \le C_2 \left( \norm{u}^2_W + \norm{v}^2_V \right)
    \]
    for any $\vec{u} = (u, v) \in \VV$.
    Now, \eqref{eq:phi3}--\eqref{ic3} can be rewritten in abstract form as 
    \begin{equation}
    \label{eq:backuniq:absform}
        \begin{split}
            & \partial_t \vpsi + \AA \vpsi = \vec{f}_{\vpsi}, \\
            & \vpsi(0) = \vpsi_0,
        \end{split}
    \end{equation}
    where 
    \[
        \vpsi = \binom{\phi}{\sigma}, \quad \vpsi_0 = \binom{\phi_0}{\sigma_0}, \quad \vec{f}_{\vpsi} = \binom{f_\phi + \phi}{f_\sigma + \sigma}. 
    \]
    Note that the additional terms in $\vec{f}_{\vpsi}$ are due to the inclusion of the identity operator in the definition of $\AA$. 
    We now aim to show that if $\vpsi(T)=0$, then necessarily $\vpsi(t)=0$ for any $t \in [0,T]$. 
    To do this, we apply \cite[Lemmas 6.1 and 6.2]{temam}. 
    Hence, we have to verify the following two conditions. 
    Firstly, equation \eqref{eq:backuniq:absform} has to admit strong solutions with regularity
    \begin{equation}
        \label{eq:backuniq:strongreg}
        \vpsi \in \HT 1 {\HH} \cap \LT \infty {\VV} \cap \LT 2 {\WW}.
    \end{equation}
    We observe that such a regularity is already available, since $\vpsi$ is defined through the difference of two strong solutions given by Theorem \ref{thm:strongsols}.
    Secondly, the right-hand side has to satisfy the following estimate
    \begin{equation}
        \label{eq:backuniq:rhsest}
        \norm{\vec{f}_{\vpsi}(t)}_{\HH} \le \alpha(t) \norm{\vpsi(t)}_{\VV} \quad \hbox{a.e.~$t \in (0,T)$}
    \end{equation}
    for some function $\alpha \in \Lt2$. 
    Therefore, we now need to prove \eqref{eq:backuniq:rhsest}, namely we estimate:
    \begin{align*}
        & \norm{\vec{f}_{\vpsi}}_{\HH} \le \norm{f_\phi}_H + \norm{f_\sigma}_H \\
        & \quad \le \norm{\Delta(F'(\phi_1)) - \Delta(F'(\phi_2))}_H + 2 \norm{P(\phi_1)(\sigma - \chi \phi + \Delta \phi - (F'(\phi_1) - F'(\phi_2)) + \chi \sigma)}_H \\
        & \qquad + 2 \norm{(P(\phi_1) - P(\phi_2))(\sigma_2 + \chi (1 - \phi_2) + \Delta \phi_2 - F'(\phi_2) + \chi \sigma_2)}_H \\
        & \qquad + \norm{c (\hh(\phi_1) - \hh(\phi_2)}_H + \norm{\phi}_H \\
        & \quad = I_1 + I_2 + I_3 + I_4 + \norm{\phi}_H\,.
    \end{align*}
    Now, by computing the laplacians and using the local Lipschitz continuity of $F''$ and $F'''$, guaranteed by \ref{ass:fc4}, the fact that $\phi_1$ and $\phi_2$ are globally bounded by Remark \ref{rmk:fderivatives}, as well as H\"older's inequality, we see that
    \begin{align*}
        I_1 & \leq \norm{F''(\phi_1) \Delta \phi}_H + \norm{(F''(\phi_1) - F''(\phi_2)) \Delta \phi_2}_H \\
        & \quad + \norm{F'''(\phi_1) (\nabla \phi_1 + \nabla \phi_2) \cdot \nabla \phi}_H +
        \norm{(F'''(\phi_1) - F'''(\phi_2)) \nabla \phi_2 \cdot \nabla \phi_2}_H \\
        & \le \norm{F''(\phi_1)}_{\Lx\infty} \norm{\Delta \phi}_H + C \norm{\Delta \phi_2}_H \norm{\phi}_{\Lx\infty} \\
        & \quad + \norm{F'''(\phi_1)}_{\Lx\infty} \norm{\nabla \phi_1 + \nabla \phi_2}_{\Lx4} \norm{\nabla \phi}_{\Lx4} + C \norm{\nabla \phi_2}^2_{\Lx4} \norm{\phi}_{\Lx\infty} \\
        & \le C \norm{\Delta \phi}_H + C \norm{\Delta \phi_2}_H \norm{\phi}_{\Hx2} + C \norm{\phi_1 + \phi_2}_{\Hx2} \norm{\phi}_{\Hx2} + C \norm{\phi_2}^2_{\Hx2} \norm{\phi}_{\Hx2} \\
        & \le C \left( 1 + \norm{\phi_2}_W + \norm{\phi_1}_W + \norm{\phi_2}^2_W \right) \norm{\phi}_W,
    \end{align*}
    where we extensively used the Sobolev embeddings $\Hx2 \hookrightarrow \Wx{1,4}$ and $\Hx2 \hookrightarrow \Lx\infty$. 
    Moreover, we also estimate similarly the remaining terms: 
    \begin{align*}
        I_2 & \le 2 \norm{P(\phi_1)}_{\Lx\infty} \norm{\sigma - \chi \phi + \Delta \phi - (F'(\phi_1) - F'(\phi_2)) + \chi \sigma}_H \\
        & \le C \left( \norm{\sigma}_H + \chi \norm{\phi}_H + \norm{\Delta \phi}_H + C \norm{\phi}_H + \chi \norm{\sigma}_H \right) \\
        & \le C \left( \norm{\sigma}_H + \norm{\phi}_W \right), \\
        I_3 & \le 2 \norm{(P(\phi_1) - P(\phi_2)}_{\Lx\infty} \norm{\sigma_2 + \chi (1 - \phi_2) + \Delta \phi_2 - F'(\phi_2) + \chi \sigma_2}_H \\
        & \le C \norm{\phi}_{\Lx\infty} \left( 1 + \norm{\sigma_2}_H + \norm{\phi_2}_H + \norm{\phi_2}_W \right) \\
        & \le C \left( 1 + \norm{\sigma_2}_H + \norm{\phi_2}_H + \norm{\phi_2}_W \right) \norm{\phi}_{W}, \\
        I_4 & \le \norm{c}_{\Lqt\infty} \norm{\hh(\phi_1) - \hh(\phi_2)}_H \le C \norm{\phi}_H.
    \end{align*}
    Putting it all together, we get that 
    \[ 
        \norm{\vec{f}_{\vpsi}}_{\HH} \le C \left( 1 + \norm{\sigma_2}_H + \norm{\phi_1}_W + \norm{\phi_2}^2_W \right) \left( \norm{\sigma}_H + \norm{\phi}_W \right) \le \alpha(t) \norm{\vpsi}_{\VV},
    \]
    where 
    \[ 
        \alpha(t) = C \left( 1 + \norm{\sigma_2}_H + \norm{\phi_1}_W + \norm{\phi_2}^2_W \right) \in \Lt\infty,
    \]
    since $\phi_i \in \LT \infty W$ and $\sigma_i \in \LT \infty H$, $i=1,2$, by Theorem \ref{thm:strongsols}.

    Hence, one can now easily use the logarithmic convexity method highlighted in \cite[Lemmas 6.1 and 6.2]{temam}, with Dirichlet quotient
    \[ \Lambda(t) = \frac{\norm{\vpsi(t)}^2_{\VV}}{\norm{\vpsi(t)}^2_{\HH}}, \]
    to conclude the proof of Theorem \ref{thm:backuniq}.
    For the sake of completeness, we report below the main steps of the cited logarithmic convexity method. 
    We recall that we aim to show that if $\vpsi(T)=0$, then necessarily $\vpsi(t)=0$ for any $t \in [0,T]$.
    Indeed, assume by contradiction that $\norm{\vpsi(t_0)}_{\HH} \neq 0$ for some $t_0 \in [0, T)$.
    Then, by continuity, $\norm{\vpsi(t)}_{\HH} \neq 0$ on some interval $(t_0, t_0 + \delta)$ and we denote by $t_1 \le T$ the largest time for which $\norm{\vpsi(t)}_{\HH} \neq 0$ on $[t_0,  t_1)$. 
    Hence, necessarily $\norm{\vpsi(t_1)}_{\HH} = 0$.
    Now, $\Lambda(t)$ is well-defined on $[t_0, t_1)$, then we can compute
    \begin{align*}
        \mezzo \ddt \Lambda(t) 
        & = \frac{(\partial_t \vpsi, \vpsi)_{\VV}}{\norm{\vpsi}^2_{\HH}} - \frac{\norm{\vpsi}^2_{\VV}}{\norm{\vpsi}^4_{\HH}} (\partial_t \vpsi, \vpsi)_{\HH} 
        = \frac{1}{\norm{\vpsi}^2_{\HH}} (\partial_t \vpsi, \AA \vpsi - \Lambda \vpsi)_{\HH} \\
        & = \frac{1}{\norm{\vpsi}^2_{\HH}} (\vec{f}_{\vpsi} - \AA \vpsi, \AA \vpsi - \Lambda \vpsi)_{\HH} 
        = - \frac{\norm{\AA \vpsi - \Lambda \vpsi}^2_{\HH}}{\norm{\vpsi}^2_H} + \frac{(\vec{f}_{\vpsi}, \AA \vpsi - \Lambda \vpsi)_{\HH}}{\norm{\vpsi^2}_{\HH}} \\
        & \le - \mezzo \frac{\norm{\AA \vpsi - \Lambda \vpsi}^2_{\HH}}{\norm{\vpsi}^2_H} + \mezzo \frac{\norm{\vec{f}_{\vpsi}}^2_{\HH}}{\norm{\vpsi}^2_{\HH}} 
        \le - \mezzo \frac{\norm{\AA \vpsi - \Lambda \vpsi}^2_{\HH}}{\norm{\vpsi}^2_H} + \mezzo \alpha^2 \Lambda,
    \end{align*}
    where we respectively used the equivalent definition of the norm $\VV$ through $\AA$, the definition of $\Lambda$, the equation \eqref{eq:backuniq:absform}, the fact that $(\Lambda \vpsi, \AA \vpsi - \Lambda \vpsi)_{\HH} = 0$, Cauchy--Schwarz and Young's inequalities and the estimate \eqref{eq:backuniq:rhsest}.
    Then, by integrating on $[t_0, t)$, for any $t \in (t_0, t_1)$, and applying Gronwall's inequality we deduce that 
    \begin{equation}
        \label{eq:backuniq:lambdaest}
        \Lambda(t) \le \Lambda(t_0) \exp\left\{ \int_{t_0}^{t_1} \alpha^2(t) \, \de t \right\} \le C \quad \hbox{for any $t \in [t_0, t_1)$,}
    \end{equation}
    since $\alpha \in \Lt2$. 
    Now, we consider the function $t \mapsto \log 1/\norm{\vpsi(t)}_{\HH}$, which is still well-defined on $[t_0, t_1)$, and by differentiation we see that
    \begin{align*}
        \ddt \log \frac{1}{\norm{\vpsi}_{\HH}} 
        & = - \mezzo \ddt \log \norm{\vpsi}^2_{\HH} = - \frac{(\partial_t \vpsi, \vpsi)_{\HH}}{\norm{\vpsi}^2_{\HH}} \\
        & = - \frac{(\vec{f}_{\vpsi} - \AA \vpsi, \vpsi)_{\HH}}{\norm{\vpsi}^2_{\HH}} 
        = \Lambda - \frac{(\vec{f}_{\vpsi}, \vpsi)_{\HH}}{\norm{\vpsi}^2_{\HH}} \\
        & \le \Lambda + \alpha \Lambda^{1/2} \le \frac32 \Lambda + \mezzo \alpha^2,
    \end{align*}
    where we used similar techniques to the ones used above.
    Finally, we integrate again in $[t_0, t)$, for any $t \in (t_0, t_1)$, and we use \eqref{eq:backuniq:lambdaest} to find that 
    \begin{align*}
        \log \frac{1}{\norm{\vpsi(t)}_{\HH}} \le \log \frac{1}{\norm{\vpsi(t_0)}_{\HH}} + \int_{t_0}^{t_1} \left( \frac32 \Lambda(s) + \mezzo \alpha^2(s) \right)\, \de s \le C.
    \end{align*}
    This shows that $1/\norm{\vpsi(t)}_{\HH}$ is bounded from above as $t \to t_1^-$, thus contradicting the fact that $\norm{\vpsi(t_1)}_{\HH} = 0$. 
    The proof of Theorem \ref{thm:backuniq} is then concluded.
\end{proof}

\begin{remark}
    Due to standard embeddings, if $(\phi_0, \sigma_0) \in W \times V$, then also $(\phi(T), \sigma(T)) \in W \times V$. 
    This means that Theorem \ref{thm:backuniq} gives a uniqueness result for the inverse problem \ref{inv:prob} if the measurements $(\phimeas, \sigmameas)$ are taken in the more regular space $W \times V$.
    However, in practical situations, one could not have such a regular measurement, which is an additional reason why in Section \ref{sec:tikhonov} we decide to propose a Tikhonov regularisation in the more general case in which the measurement is just in $H \times H$, even if the uniqueness is not guaranteed.
\end{remark}

\begin{remark}
\label{rmk:chemotaxis}
    In order to make considerations about the physical viability of assumption \ref{ass:chi}, we may review the arguments in the proof of Theorem \ref{thm:backuniq} starting from system \eqnhdad, i.e. reintroducing all the physical parameters in the model. Note that in this case the operator $\mathbb{A}$ is not symmetric, hence the proof of the backward uniqueness through logarithmic convexity becomes more involved. Anyhow, it is easy to obtain the following expression for the parameter $\gamma$ in \eqref{eqn:Apd}:
    \[
    \gamma=\tilde{\eps}^2\frac{M_{\phi}}{M_{\sigma}}-\frac{\tilde{\chi}^2\left(\frac{M_{\phi}}{M_{\sigma}}+1\right)^2}{8(1-\delta)}.\]
    Hence, assumption \ref{ass:chi} becomes:
    \[
    \tilde{\chi}^2<\frac{8\tilde{\eps}^2\frac{M_{\phi}}{M_{\sigma}}}{\left(\frac{M_{\phi}}{M_{\sigma}}+1\right)^2},
    \]
    which, given the expressions of $\tilde{\chi}$ and $\tilde{\eps}$ in terms of the dimensional parameters $\chi$ and $\eps$, becomes:
    \begin{equation}
    \label{eqn:chidimbu}
    \chi^2<\frac{8\eps^2M_{\phi}P_0}{(M_{\phi}+M_{\sigma})^2}=:\iota.
    \end{equation}
    Referring to the reference biological ranges for the
values of the model parameters reported in \cite[Table 1]{Agosti1}, we can deduce that $\chi \in [9.64,456521.18] \, Pa$, while $\iota \in [2.28,16942.54] \, Pa^2$. Hence, we conclude that for small values of the chemotactic parameter $\chi$ the assumption \ref{ass:chi} is biologically feasible.
\end{remark}

Due to the higher regularity of solutions to \eqref{eq:phi}--\eqref{ic}, as a byproduct of the previous proof of the backward uniqueness, we can also prove a stronger continuous dependence result. 

\begin{corollary}
	\label{thm:contdep_strong}
	Assume hypotheses \emph{\ref{ass:setting}--\ref{ass:hc}} and \emph{\ref{ass:fc4}--\ref{ass:pc1}}. Let ${\phi_0}_1$, ${\sigma_0}_1$ and ${\phi_0}_2$, ${\sigma_0}_2$ be two sets of data satisfying \emph{\ref{ass:initial2}} and let $(\phi_1, \mu_1, \sigma_1)$ and $(\phi_2, \mu_2, \sigma_2)$ two corresponding strong solutions as in Theorem \emph{\ref{thm:strongsols}}. Then, there exists a constant $K>0$, depending only on the data of the system and on the norms of $\{ ({\phi_0}_i, {\sigma_0}_i) \}_{i=1,2}$, but not on their difference, such that
	\begin{equation}
		\label{eq:contdep_estimate_strong}
		\begin{split}
			& \norm{\phi_1 - \phi_2}_{\HT 1 H \cap \LT \infty W \cap \LT 2 {\Hx4}} 
			+ \norm{\mu_1 - \mu_2}_{\LT 2 W} \\
			& \quad + \norm{\sigma_1 - \sigma_2}_{\HT 1 H \cap \LT \infty V \cap \LT 2 W}
            \le K \left( \norm{{\phi_0}_1 - {\phi_0}_2}_W + \norm{{\sigma_0}_1 - {\sigma_0}_2}_V \right).
		\end{split}
	\end{equation}
\end{corollary}

\begin{proof}
    First of all, note that the weak continuous dependence estimate \eqref{eq:contdep_estimate} can be now obtained without assuming hypothesis \ref{ass:phinf}, since $\phi \in \Cqt0$ and Remark \ref{rmk:fderivatives} holds.
    We then recall the compact notation introduced in the proof of Theorem \ref{thm:backuniq}, where we equivalently wrote the system solved by the differences of two strong solutions $\phi = \phi_1 - \phi_2$, $\mu = \mu_1 - \mu_2$ and $\sigma = \sigma_1 - \sigma_2$ as 
    \begin{equation*}
        \begin{split}
            & \partial_t \vpsi + \AA \vpsi = \vec{f}_{\vpsi}, \\
            & \vpsi(0) = \vpsi_0.
        \end{split}
    \end{equation*}
    Next, we observe that the core of the argument was already done in the proof of Theorem \ref{thm:backuniq} when showing the validity of the estimate \eqref{eq:backuniq:rhsest}. 
    Indeed, if \eqref{eq:backuniq:rhsest} holds, by testing \eqref{eq:backuniq:absform} by $\partial_t \vpsi + \AA \vpsi$ and using Cauchy--Schwarz and Young's inequalities, we get that 
    \begin{align*}
        & \norm{\partial_t \vpsi}^2_{\HH} + \ddt \norm{\vpsi}^2_{\VV} + \norm{\AA \vpsi}^2_{\HH} 
        = (\vec{f}_{\vpsi}, \partial_t \vpsi + \AA \vpsi)_{\HH} \\
        & \quad \le \norm{\vec{f}_{\vpsi}}_{\HH} \left( \norm{\partial_t \vpsi}_{\HH} + \norm{\AA \vpsi}_{\HH} \right) \le \alpha(t) \norm{\vpsi(t)}_{\VV} \left( \norm{\partial_t \vpsi}_{\HH} + \norm{\AA \vpsi}_{\HH} \right) \\
        & \quad \le \mezzo \norm{\partial_t \vpsi}^2_{\HH} + \mezzo \norm{\AA \vpsi}^2_{\HH} + C \alpha(t)^2 \norm{\vpsi(t)}^2_{\VV}. 
    \end{align*}
    Then, the continuous dependence estimate \eqref{eq:contdep_estimate_strong} is simply obtained in compact notation by integrating on $(0,t)$, for any $t \in (0,T)$, and applying Gronwall's lemma, together with the fact that $\vpsi_0$ is assumed to be in $\VV$.
    Indeed, what is done above is equivalent to taking the system \eqref{eq:phi2}--\eqref{ic2} solved by the difference of two solutions, testing \eqref{eq:phi2} by $\partial_t \phi + \Delta^2 \phi$ and \eqref{eq:sigma2} by $\partial_t \sigma - \Delta \sigma$ and summing them up, while computing explicitly the expression of $\Delta \mu$. 
    Additionally, note that, in this case, we can drop assumption \ref{ass:chi} on the smallness of the chemotactic coefficient, since the two estimates need not be taken simultaneously. 
    In particular, one can first test only \eqref{eq:sigma2} by $\partial_t \sigma - \Delta \sigma$ and recover the estimate for $\sigma$, by using the facts that an $\LT 2 W$-estimate for $\phi$ is already available from \eqref{eq:contdep_estimate} and that the coefficient $\alpha(t)$ in \eqref{eq:backuniq:rhsest} is actually in $\Lt \infty$ in our case.
    Then, one can later test \eqref{eq:phi2} by $\partial_t \phi + \Delta^2 \phi$ and use the newfound $\LT 2 W$-estimate for $\sigma$, thus removing the need of tuning the chemotactic coefficient. 
    Hence, we leave the details of the proof to the interested reader and consider the proof of Corollary \ref{thm:contdep_strong} concluded.
\end{proof}

\begin{remark}
    \label{rmk:lipstab}
    As a side comment, we mention that the strong regularity of the solutions proven in Theorem \ref{thm:strongsols} and the strong continuous dependence estimate \eqref{eq:contdep_estimate_strong} would be enough to prove that the operator $\Rcal : W \times V \to H \times H$ is Fr\'echet-differentiable between these stronger spaces and that its Fr\'echet derivative is continuous as a function from $W \times V$ to the space $\mathcal{L}(W \times V, H \times H)$, as similarly done in \cite{BCFLR2024}.
    Moreover, one should also be able to prove the injectivity of the Fr\'echet derivative $\D \Rcal(\phiob, \sigmaob)$ for any fixed point $(\phiob, \sigmaob) \in W \times V$.
    Indeed, this can be done by proving a backward uniqueness result in the spirit of Theorem \ref{thm:backuniq} for the linearised system which will be introduced in the following \eqref{eq:xi}--\eqref{icl}.
    Such results, then, by a suitable version of the inverse map Theorem proved in \cite[Theorem 2.1]{B2013}, would imply a Lipschitz stability estimate for the inverse problem if one assumes to be reconstructing the initial data in a finite-dimensional subspace of $W \times V$, defined e.g. via a $\mathcal{C}^1$-conforming finite element space. 
    However, the possibility of refining this Lipschitz stability estimate by quantifying its constant, as in \cite{BCFLR2024}, is left open for further investigation, due to the more complex fourth-order structure of the system.

\end{remark}


\section{Tikhonov Regularisation}
\label{sec:tikhonov}

In this section, we aim to approximate the solution to the proposed inverse problem \eqref{inv:prob} employing a Tikhonov regularisation. 
This will be done in the most useful setting in practice, which is when the final measurements are taken in $H \times H$ and the initial data are assumed to be in $V \times H$. 
Thus, to carry out our analysis, we only use the regularity of weak solutions by Theorem \ref{thm:weaksols} and the corresponding continuous dependence estimate by Theorem \ref{thm:contdep}.
In this way, we just assume the minimal assumptions on the initial data that guarantee that the final values $(\phi(T), \sigma(T))$ are well-defined in $H \times H$, namely that $(\phi_0, \sigma_0) \in V \times H$.
We recall that the results presented in the previous Section \ref{sec:uniq} show that the solution to the inverse problem is unique in the more regular class $\Hx2 \times \Hx1$.
Hence, one could also think of introducing a higher-order Tikhonov regularisation to approximate the initial data in this more regular space, where uniqueness is guaranteed. 
However, this would mean that the computational cost of solving the problem numerically would become significantly higher, due to the need for higher-order finite element spaces. 
For this reason and also because we are still technically able to characterise the solutions to the Tikhonov-regularised problem even by starting only from weak solutions, we stick to the above-mentioned less regular setting.

We now introduce the Tikhonov-regularisation of our inverse reconstruction problem and then study it by interpreting it as a constrained minimisation problem on the initial data. 
Indeed, we study the following:

\bigskip
\noindent(MP) \textit{Minimise the cost functional}
\begin{align*}
		\mathcal{J}(\phi_0, \sigma_0)
		& = \, \frac{\lambda_1}{2} \int_{\Omega} |\phi(T) - \phimeas|^2 \,\de x  
		+ \alpha_1 \int_{\Omega} \left( F(\phi_0) + \mezzo |\nabla \phi_0|^2 \right) \,\de x \\
        & \quad + \frac{\lambda_2}{2} \int_{\Omega} |\sigma(T) - \sigmameas|^2 \,\de x  
		+ \frac{\alpha_2}{2} \int_{\Omega} |\sigma_0|^2 \,\de x,
\end{align*}
\textit{subject to the constraints}
\begin{align*} 
	& \phi_0 \in \Uad = \left\{ \phi_0 \in V \mid 0 \le \phi_0 \le 1 \text{ a.e. in } \Omega \right\}, \\
    & \sigma_0 \in \Vad = \left\{ \sigma_0 \in H \mid 0 \le \sigma_0 \le 1 \text{ a.e. in } \Omega \right\},
\end{align*}
 and where 
$( \phi(T), \sigma(T) ) = \Rcal ( \phi_0, \sigma_0 )$.
\bigskip

\noindent 
Notice that, in the cost functional, we use an $\Lx2$-penalisation for $\sigma_0$, whereas a Ginzburg--Landau energy penalisation for $\phi_0$. 
Indeed, since we are dealing with a phase-field model, we believe that the most natural choice to stabilize the inverse reconstruction is to use the phase-field version of the perimeter of the interface \cite{MM1977}. 
As a matter of fact, in most situations, the initial condition to be reconstructed is constituted by domains of well-separated phases, with homogeneous values given by minima of $F$.
The same procedure was already used in \cite{GHK2019} when studying a similar problem on a phase-field model for two-phase flows.
It is also a common penalisation in many inverse problems solved by the phase-field method \cite{cavelast, cavsemilinear, condinclusion}.

Regarding the parameters at play, we make the following hypotheses:
\begin{enumerate}[font = \bfseries, label = C\arabic*., ref=\bf{C\arabic*}]
	\item\label{ass:C1} $\lambda_1, \lambda_2, \alpha_1, \alpha_2 \ge 0$, but not all equal to $0$.
	\item\label{ass:C2} $\phimeas, \sigmameas \in L^2(\Omega)$.
	\item\label{ass:fgrowth4} $F \in \mathcal{C}^3(\R)$ and the exponent $s$ in hypothesis \ref{ass:fconv} is such that $s \in [2,4]$.
	\item\label{ass:ph_second} $P, \hh \in W^{2,\infty}(\R)$.
\end{enumerate}

\begin{remark}
	In hypothesis \ref{ass:fgrowth4} we only allow growth up to the fourth order for the convex part of the potential $F$, while in \ref{ass:ph_second} we assume higher regularity and boundedness for $P$ and $\hh$. 
    Note that the key examples \eqref{f:example} and \eqref{p:example} are still included in this setting, up to some smoothing of $P$.
	Such stronger hypotheses are needed to prove Fr\'echet differentiability of the operator $\Rcal$, without assuming additional regularity on the initial data and consequently using the strong regularity of the solutions.
    This can be seen as a technical novelty of this work, as in similar papers on optimal control problems for Cahn--Hilliard tumour growth models (e.g., \cite{CGRS2017, CRW2021, CSS2021_secondorder}) the Fr\'echet differentiability was always proved starting from strong solutions.
    As a matter of fact, we use only the regularity of weak solutions, given by Theorem \ref{thm:weaksols}, to characterise the solutions to the Tikhonov-regularised problem (MP). 
	In this way, we can resort to a lower order Tikhonov regularisation, namely only in $\Hx1$ for $\phi_0$ and $\Lx2$ for $\sigma_0$.
    This gives great advantages when dealing with the numerical approximation of the problem. 
    Indeed, we can limit the computational cost of the numerical methods by avoiding using higher-order finite element spaces and using standard conforming finite element methods with the lowest order instead.
\end{remark}

First, we prove the existence of a solution to the regularised minimisation problem.

\begin{theorem}
	\label{thm:control_existence}
	Assume hypotheses \emph{\ref{ass:setting}--\ref{ass:hc}} and \emph{\ref{ass:C1}}--\emph{\ref{ass:C2}}. Then, the constrained minimisation problem (MP) admits at least one solution $(\phiob, \sigmaob) \in \Uad \times \Vad$, such that if $(\phib(T), \sigmab(T)) = \Rcal(\phiob, \sigmaob)$ is the solution to \eqref{eq:phi}--\eqref{ic} associated to $\phiob$ and $\sigmaob$, one has that
	\begin{equation}
		\mathcal{J}(\phiob, \sigmaob) = \min_{(\phi_0, \sigma_0) \,\in \,\Uad \times \Vad} \, \mathcal{J}(\phi_0, \sigma_0).
	\end{equation}
\end{theorem}

\begin{proof} 
	First of all, we note that $\mathcal{J}$ is bounded from below, since most of its terms are non-negative and hypothesis \ref{ass:fbelow} holds on $F$. 
    Then, we let $\{ (\phi_0^n, \sigma_0^n) \}_{n \in \N} \subset \Uad \times \Vad$ be a minimising sequence such that 
    \begin{equation*}
        \lim_{n \to + \infty} \mathcal{J} (\phi_0^n, \sigma_0^n) = \inf_{ (\phi_0, \sigma_0) \, \in \, \Uad \times \Vad } \mathcal{J} (\phi_0, \sigma_0) \ge - C,
    \end{equation*}
    Since $\{ (\phi_0^n, \sigma_0^n) \}_{n \in \N} \subset \Uad \times \Vad$, we have that $\{ (\phi_0^n, \sigma_0^n) \}$ are uniformly bounded in $\Lx\infty \times \Lx\infty$, therefore we deduce that there exist $\phiob, \sigmaob \in \Lx\infty$ such that, up to a subsequence, 
    \begin{equation*}
        \phi_0^n \weakstar \phiob \quad \hbox{and} \quad \sigma_0^n \weakstar \sigmaob \quad \hbox{weakly star in $\Lx\infty$.}
    \end{equation*}
    Moreover, without loss of generality, we can assume that $\mathcal{J}$ is uniformly bounded along the minimising sequence, i.e.~$\mathcal{J}(\phi_0^n, \sigma_0^n) \le C$ for some constant $C \gs 0$. 
    In particular, this means that the sequence $\{\nabla \phi_0^n\}$ is uniformly bounded in $\Lx2$, hence we also infer that, up to a further subsequence, 
    \[
        \phi_0^n \weak \phiob \quad \hbox{weakly in $\Hx1$.}
    \]
    Finally, since $\Uad \times \Vad$ is convex and closed in $\Hx1 \times \Lx\infty$, it is also weakly-star sequentially closed and thus $(\phiob, \sigmaob) \in \Uad \times \Vad$. 

    Next, we consider the corresponding weak solutions $(\phi^n, \mu^n, \sigma^n)$ to \eqref{eq:phi}--\eqref{ic} and observe that, since $\{(\phi_0^n, \sigma_0^n)\}$ are uniformly bounded in $\Hx1 \times \Lx2$, they are uniformly bounded in the spaces of weak solutions by Theorem \ref{thm:weaksols}. 
    Therefore, by Banach-Alaoglu's Theorem, we deduce that, up to a subsequence, 
    \begin{align*}
        & \phi^n \weakstar \phib \quad \hbox{weakly star in $\HT 1 {V^*} \cap \LT \infty V \cap \LT 2 {\Hx3}$,} \\
        & \mu^n \weak \mub \quad \hbox{weakly in $\LT 2 V$,} \\
        & \sigma^n \weakstar \sigmab \quad \hbox{weakly star in $\HT 1 {V^*} \cap \LT \infty H \cap \LT 2 V$.}
    \end{align*}
    Now, by the compact embeddings of Aubin--Lions--Simon (see \cite[Section 8, Corollary 4]{simon}), it follows that $\phi^n \to \phib$ strongly in $\C 0 {\Lx{p}}$, for any $2 \le p \ls 6$. 
    In particular, up to a further subsequence, $\phi^n \to \phib$ a.e.~in $Q_T$, hence, due to the continuity of $F'$, $P$ and $\hh$, also $F'(\phi^n) \to F'(\phib)$, $P(\phi^n) \to P(\phib)$ and $\hh(\phi^n) \to \hh(\phib)$ a.e.~in $Q_T$.
    Then, by standard application of Lebesgue, Egorov and Vitali's convergence theorems, with such pieces of information and hypotheses \ref{ass:fconv}, \ref{ass:p} and \ref{ass:hc} on the growth of $F$, $P$ and $\hh$ respectively, one can pass to the limit in the weak formulation \eqref{varform:phi}--\eqref{varform:sigma} and deduce that $(\phib, \mub, \sigmab)$ is a weak solution to \eqref{eq:phi}--\eqref{ic} with $(\phiob, \sigmaob)$.
    Then, we can infer that 
    \[
        \inf_{(\phi_0, \sigma_0) \, \in \, \Uad \times \Vad} \mathcal{J}(\phi_0, \sigma_0) \le \mathcal{J}(\phiob, \sigmaob).
    \]
    Finally, we call $\Jtil(\phi(T), \sigma(T), \phi_0, \sigma_0)$ the functional defined on $H \times H \times V \times H$ such that $\mathcal{J}(\phi_0, \sigma_0) = \Jtil(\Rcal(\phi_0, \sigma_0), \phi_0, \sigma_0)$.   
    Then, we observe that $\Jtil$ is weakly lower-semicontinuous as a functional defined on the space $H \times H \times V \times H$, being a sum of weakly lower-semicontinuous functionals. 
    Indeed, it is easy to see that most of the terms are weakly lower-semicontinuous, as they are essentially $\Lx2$-norms. 
    We only have to check the term $\int_\Omega F(\phi_0) \, \de x$, but, since $V$ is compactly embedded $\Lx{p}$ for any $p \ls 6$ and $F$ has growth up to the power $s \ls 6$, we easily deduce that it is weakly continuous in $V$.
    Consequently, we can now infer that
    \begin{align*}
        \mathcal{J}(\phiob, \sigmaob) 
        & \le \liminf_{n \to +\infty} \mathcal{J}(\phi_0^n, \sigma_0^n) \\
        & = \inf_{(\phi_0, \sigma_0) \, \in \, \Uad \times \Vad} \mathcal{J}(\phi_0, \sigma_0) \le \mathcal{J}(\phiob, \sigmaob),
    \end{align*}
    which means that $(\phiob, \sigmaob) \in \Uad \times \Vad$ is optimal. 
    This concludes the proof of Theorem \ref{thm:control_existence}.
\end{proof}

\begin{remark}
    As a consequence of the central part of the proof of Theorem \ref{thm:control_existence}, we can immediately say that the operator $\Rcal: V \times H \to H \times H$ is weakly sequentially closed. 
    Namely, this means that if $(\phi_0^n, \sigma_0^n) \weak (\phiob, \sigmaob)$ in $V \times H$ and $\Rcal(\phi_0^n, \sigma_0^n) \weak (\xi, \eta)$ in $H \times H$, then $(\xi, \eta) = \Rcal(\phiob, \sigmaob) = (\phib(T), \sigmab(T))$. 
    This is essentially due to the fact that the limit variables $(\phib, \sigmab)$ still satisfy the forward system.  
\end{remark}

\subsection{Differentiability of the solution mapping} 

We now prove Fr\'echet differentiability of the operator $\Rcal : V \times H \to H \times H$. 
As an \emph{ansatz} for the Fr\'echet derivative, we first introduce the corresponding linearised system to \eqref{eq:phi}--\eqref{ic}, which takes the form
\begin{alignat}{2}
	& \partial_t \xi - \Delta \eta 
	= P'(\phib) (\sigmab + \chi (1 -\phib) - \mub) \xi + P(\phib) (\rho - \chi \xi - \eta) - \hh'(\phib) c \, \xi 
	\qquad && \hbox{in  $Q_T$},  \label{eq:xi} \\
	& \eta 
	= - \Delta \xi + F''(\phib) \xi - \chi \rho 
	\qquad && \hbox{in  $Q_T$},  \label{eq:eta} \\
	& \partial_t \rho - \Delta \rho + \chi \Delta \xi
	= - P'(\phib) (\sigmab + \chi (1 -\phib) - \mub) \xi - P(\phib) (\rho - \chi \xi - \eta) - \rho 
	\qquad && \hbox{in  $Q_T$}, \label{eq:rho} \\
	& \partial_{\n} \xi = \partial_{\n} \eta = \partial_{\n} \rho = 0 
	\qquad && \hbox{on  $\Sigma_T$}, \label{bcl} \\
	& \xi(0) = h, \quad \rho(0) = k 
	\qquad && \hbox{in  $\Omega$}, \label{icl}
\end{alignat}
where $(h,k) \in V \times H$ are increments and $(\phib, \mub, \sigmab)$ is the weak solution to \eqref{eq:phi}--\eqref{ic}.
Then, we state a well-posedness result for the linearised system \eqref{eq:xi}--\eqref{icl}.

\begin{theorem}
	\label{thm:linearised}
	Assume hypotheses \emph{\ref{ass:setting}--\ref{ass:hc}}, \emph{\ref{ass:phinf}} and \emph{\ref{ass:fc4}--\ref{ass:pc1}}. Let $(\phib, \mub, \sigmab)$ be the weak solution to \eqref{eq:phi}--\eqref{ic}, corresponding to some initial data $(\phi_0, \sigma_0) \in V \times H$. Then, for any $(h,k) \in H \times H$, the linearised system \eqref{eq:xi}--\eqref{icl} admits a unique weak solution, which is uniformly bounded in the following spaces
	\begin{align*}
		& \xi \in H^1(0,T;W^*) \cap \mathcal{C}^0([0,T]; H) \cap L^2(0,T;W), \\
		& \eta \in L^2(0,T;H), \\
		& \rho \in H^1(0,T;V^*) \cap \mathcal{C}^0([0,T]; H) \cap L^2(0,T;V),
	\end{align*}
	and fulfils \eqref{eq:xi}--\eqref{icl} in variational form, i.e.~it satisfies
	\begin{align}
		& \duality{\xi_t, w}_W - (\eta, \Delta w)_H = (P'(\phib)(\sigmab + \chi (1 -\phib) - \mub) \xi + P(\phib)(\rho - \chi \xi - \eta) - \hh'(\phib) c \, \xi, w)_H, \label{varform:xi}\\
		& (\eta,w)_H = - (\Delta \xi,w)_H + (F''(\phib) \xi, w)_H - \chi (\rho, w)_H, \label{varform:eta} \\
		& \duality{\rho_t,v}_V + (\nabla \rho - \chi \nabla \xi, \nabla v)_H \nonumber \\
		& \quad = - (P'(\phib)(\sigmab + \chi (1 -\phib) - \mub) \xi + P(\phib)(\rho - \chi \xi - \eta) - \kappa \rho,v)_H, \label{varform:rho}
	\end{align}
	for a.e. $t \in (0,T)$ and for any $(w, v) \in W \times V$, and $\xi(0) = h$, $\rho(0) = k$.

\end{theorem}

\begin{proof}
    Since the system is linear, we proceed formally. 
    The argument can then be made rigorous by employing a Galerkin discretisation scheme.
    To begin, we test \eqref{eq:xi} by $\xi$, \eqref{eq:eta} by $\Delta \xi$, \eqref{eq:rho} by $\rho$ and sum them up to obtain:
    \begin{equation}
        \label{eq:linear:est1}
        \begin{split}
            & \mezzo \ddt \norm{\xi}^2_H + \mezzo \ddt \norm{\rho}^2_H + \norm{\Delta \xi}^2_H + \norm{\nabla \rho}^2_H \\
            & \quad = (F''(\phib)\xi, \Delta \xi)_H - 2 \chi (\rho, \Delta \xi)_H + (P(\phib) (\rho - \chi \xi - \eta), \xi - \rho)_H \\
            & \qquad + (P'(\phib)(\sigmab + \chi (1 - \phib) - \mub) \xi, \xi - \rho)_H - (\hh'(\phib) c \, \xi, \xi)_H - \norm{\rho}^2_H.
        \end{split}
    \end{equation}
    We now estimate the terms on the right-hand side one by one. 
    For the first two ones, we use Cauchy--Schwarz and Young's inequalities to infer that
    \begin{align*}
        & (F''(\phib) \xi, \Delta \xi)_H - 2 \chi (\rho, \Delta \xi)_H \\
        & \le \frac{1}{4} \norm{\Delta \xi}^2_H + C \norm{F''(\phib)}^2_{\Lx\infty} \norm{\xi}^2_H + C \norm{\rho}^2_H, 
    \end{align*}
    where $\norm{F''(\phib)}^2_{\Lx\infty} \in \Lt1$, since $\phib$ is uniformly bounded in $\LT 8 {\Lx\infty}$ by \eqref{eq:emb_l8linf} and \ref{ass:fconv}. 
    Next, by using equation \eqref{eq:eta}, we immediately see that
    \begin{equation}
        \label{eq:linear:etaH}
        \norm{\eta}^2_H \le \norm{\Delta \xi}^2_H + \norm{F''(\phib)}^2_{\Lx\infty} \norm{\xi}^2_H + \chi \norm{\rho}^2_H.
    \end{equation}
    Hence, by recalling also hypothesis \ref{ass:phinf}, we can similarly estimate
    \begin{align*}
        & (P(\phib)(\rho - \chi \xi - \eta), \xi - \rho)_H - (\hh'(\phi) c \, \xi, \xi)_H \\
        & \quad \le P_\infty (\norm{\rho}_H + \chi \norm{\xi}_H + \norm{\eta}_H)(\norm{\xi}_H + \norm{\rho}_H) + \hh'_\infty c_\infty \norm{\xi}^2_H \\
        & \quad \le \frac{1}{4} \norm{\Delta \phi}^2_H + C \left( 1 + \norm{F''(\phib)}^2_{\Lx\infty} \right) \norm{\xi}^2_H + C \norm{\rho}^2_H\,.
    \end{align*}
    Finally, by using again \ref{ass:phinf}, H\"older and Young's inequalities and Sobolev embeddings, we also deduce that
    \begin{align*}
        & (P'(\phib)(\sigmab + \chi (1 - \phib) - \mub) \xi, \xi - \rho)_H \\
        & \quad \le P'_\infty \norm{\sigmab + \chi (1 - \phib) - \mub}_{\Lx6} \norm{\xi}_{\Lx3} (\norm{\xi}_H + \norm{\rho}_H) \\
        & \quad \le P'_\infty \norm{\sigmab + \chi (1 - \phib) - \mub}_V (\norm{\Delta \xi}_H + \norm{\xi}_H) (\norm{\xi}_H + \norm{\rho}_H) \\
        & \quad \le \frac{1}{4} \norm{\Delta \xi}^2_H + C\norm{\xi}^2_H + C \norm{\sigmab + \chi (1 - \phib) - \mub}^2_V (\norm{\xi}^2_H + \norm{\rho}^2_H), 
    \end{align*}
    where $\norm{\sigmab + \chi (1 - \phib) - \mub}^2_V$ is uniformly bounded in $\Lt1$ by Theorem \ref{thm:weaksols}.
    Therefore, by integrating \eqref{eq:linear:est1} on $(0,t)$, for any $t \in (0,t)$, and putting everything together, we get the inequality:
    \begin{align*}
        & \mezzo \norm{\xi(t)}^2_H + \mezzo \norm{\rho(t)}^2_H + \frac{1}{4} \int_0^t \norm{\Delta \xi}^2_H \, \de s + \int_0^t \norm{\nabla \rho}^2_H \, \de s \\
        & \quad \le \mezzo \norm{h}^2_H + \mezzo \norm{k}^2_H + C \int_0^T \left( 1 + \norm{\sigmab + \chi (1 - \phib) - \mub}^2_V \right) \norm{\rho}^2_H \, \de s \\
        & \qquad + C \int_0^T \left( 1 + \norm{F''(\phib)}^2_{\Lx\infty} + \norm{\sigmab + \chi (1 - \phib) - \mub}^2_V \right) \norm{\xi}^2_H \, \de s.
    \end{align*}
    Therefore, by Gronwall's inequality, we deduce that 
    \begin{equation}
        \label{eq:linear:weaknorms}
        \norm{\xi}^2_{\LT \infty H \cap \LT 2 W} + \norm{\rho}^2_{\LT \infty H \cap \LT 2 V} \le C \left( \norm{h}^2_H + \norm{k}^2_H \right).
    \end{equation}
    Then, by \eqref{eq:linear:etaH} and \eqref{eq:linear:weaknorms}, we immediately infer that also
    \begin{equation}
        \label{eq:linear:etal2h}
        \norm{\eta}^2_{\LT 2 H} \le C \left( \norm{h}^2_H + \norm{k}^2_H \right).
    \end{equation}
    Moreover, by comparison in \eqref{eq:xi} and \eqref{eq:eta}, by exploiting also \ref{ass:phinf}, we additionally get that
    \begin{equation*}
        \norm{\xi}^2_{\HT 1 {W^*}} + \norm{\rho}^2_{\HT 1 {V^*}} \le C \left( \norm{h}^2_H + \norm{k}^2_H \right).
    \end{equation*}
    Thus, the first part of Theorem \ref{thm:linearised} is proved. 
    Uniqueness follows by the linearity of the system. 
\end{proof}

Next, we are now able to prove that the map $\Rcal$ is Fr\'echet differentiable and characterise its derivative as the solution to the linearised system.

\begin{theorem}
	\label{thm:frechet}
	Assume hypotheses \emph{\ref{ass:setting}--\ref{ass:hc}} and \emph{\ref{ass:fgrowth4}--\ref{ass:ph_second}}. Then, $\Rcal: V \times H \to H \times H$ is Fr\'echet differentiable, i.e. for any $(\phiob, \sigmaob) \in V \times H$ there exists a unique Fr\'echet derivative $\D\Rcal(\phiob, \sigmaob) \in \mathcal{L}(V \times H, H \times H)$ such that, as $\norm{(h,k)}_{V \times H} \to 0$,
	\begin{equation}
		\label{frechet:diff}
		\frac{ \norm{ \Rcal(\phiob + h, \sigmaob + k) - \mathcal{R}(\phiob, \sigmaob) - \D\mathcal{R}(\phiob, \sigmaob)[h,k]  }_{H \times H}}{ \norm{(h,k)}_{V \times H} } \to 0.
	\end{equation}
	Moreover, for any $(h,k) \in V \times H$, the Fr\'echet derivative at $(\phiob, \sigmaob)$ in $(h,k)$ is defined as 
	\[ \D\mathcal{R}(\phiob, \sigmaob)[h,k] = (\xi(T), \rho(T))  \]
	where $\xi(T)$ and $\rho(T)$ are the solutions to the linearised system \eqref{eq:xi}--\eqref{icl} with initial data $(h,k)$, evaluated at the final time. 

\end{theorem}

\begin{proof}
    We observe that it is sufficient to prove the result for any small enough perturbation $(h,k)$, i.e.~we fix $\Lambda > 0$ and consider only perturbations such that
	\begin{equation}
		\label{hk:bound}
		\norm{(h,k)}_{V \times H} \le \Lambda.
	\end{equation} 
	Now, we fix $\phiob$, $\sigmaob$, $h$ and $k$ as above and consider
	\begin{align*}
		& (\phi, \mu, \sigma) := \Rcal(\phiob + h, \sigmaob +k), \\
		& (\phib, \mub, \sigmab) := \Rcal(\phiob, \sigmaob), \\
		& (\xi, \eta, \rho) \text{ as the solution to \eqref{eq:xi}--\eqref{icl} with respect to } (h,k).
	\end{align*}
	In order to show Fr\'echet differentiability, then, it is enough to show that there exists a constant $C>0$, depending only on the parameters of the system and possibly on $\Lambda$, and an exponent $q > 2$ such that
	\[ \norm{ (\phi(T), \sigma(T)) - (\phib(T), \sigmab(T)) - (\xi(T), \rho(T)) }^2_{H \times H} \le C \norm{(h,k)}^q_{V \times H}. \]
	To do this, we introduce the additional variables
	\begin{align*}
		& \psi := \phi - \phib - \xi \in H^1(0,T;W^*) \cap \C 0 H \cap  L^2(0,T; W), \\
		& \zeta  := \mu - \mub - \eta \in L^2(0,T;V), \\
		& \theta := \sigma - \sigmab - \rho \in H^1(0,T;V^*) \cap \C 0 H \cap  L^2(0,T; V), 
	\end{align*}
	which, by Theorems \ref{thm:weaksols} and \ref{thm:linearised} enjoy the regularities shown above. Then, this is equivalent to showing that
	\begin{equation}
		\label{frechet:aim}
		\norm{ (\psi(T), \theta(T)) }_{H \times H}^2 \le C \norm{(h,k)}^q_{V \times H}.
	\end{equation}
	By inserting the equations solved by the variables in the definitions of $\psi$, $\zeta$ and $\theta$ and exploiting the linearity of the involved differential operators, we infer that $\psi$, $\zeta$ and $\theta$ formally satisfy the equations:
	\begin{alignat}{2}
		& \partial_t \psi - \Delta \zeta = Q^h - U^h \,\,\, && \text{in } Q_T,  \label{eq:psi}\\
		& \zeta =  - \Delta \psi + F^h - \chi \theta \quad && \text{in } Q_T,  \label{eq:zeta} \\
		& \partial_t \theta - \Delta \theta + \chi \Delta \psi = - Q^h - \theta \quad && \text{in } Q_T, \label{eq:theta}
	\end{alignat}
	together with boundary and initial conditions:
	\begin{alignat}{2}
		& \partial_{\n} \psi = \partial_{\n} \zeta = \partial_{\n} \theta = 0 \qquad && \text{on } \Sigma_T, \label{bcd} \\
		& \psi(0) = 0, \quad \theta(0) = 0 \qquad && \text{in } \Omega, \label{icd}
	\end{alignat}
	where:
	\begin{align*}
		F^h & = F'(\phi) - F'(\phib) - F''(\phib)\xi, \\
		Q^h & = P(\phi)(\sigma + \chi (1 - \phi) - \mu) - P(\phib)(\sigmab + \chi (1 - \phib) - \mub) \\
        & \quad - P(\phib)(\rho - \chi \xi - \eta) - P'(\phib) (\sigmab + \chi (1 - \phib) - \mub) \xi, \\
		U^h & = c \hh(\phi) - c \hh(\phib) - c \hh'(\phib) \xi.
	\end{align*}
	Note that, to be precise, system \eqref{eq:psi}--\eqref{icd} has to be understood in weak sense, i.e.~through a variational formulation, since only weak regularity is available.
	Before going on, we can rewrite in a better way the terms $F^h$, $Q^h$ and $U^h$, by using the following version of Taylor's theorem with integral remainder for any real function $f \in \mathcal{C}^2$ at a point $x_0 \in \R$:
	\[ f(x) = f(x_0) + f'(x_0) (x-x_0) + \left( \int_0^1 (1-z) f''(x_0 + z(x-x_0)) \, \de z \right) (x-x_0)^2. \]
	Indeed, with straightforward calculations, one can see that
	\begin{align*}
		& F^h = F''(\phib) \psi + R_1^h (\phi - \phib)^2, \\
		& U^h = c \hh'(\phib) \psi + c R_2^h (\phi - \phib)^2, 
	\end{align*}
	and also, up to adding and subtracting some additional terms, that
	\begin{align*}
		Q^h & = P(\phib) (\theta - \chi \psi - \zeta) + P'(\phib) (\sigmab + \chi (1 - \phib) - \mub) \, \psi  \\
		& \quad + (P(\phi) - P(\phib))[ (\sigma - \sigmab) - \chi (\phi - \phib) - (\mu - \mub) ] + R_3^h (\sigmab + \chi (1 - \phib) - \mub) (\phi - \phib)^2 \\
	\end{align*}
	where
	\begin{gather*}
		R_1^h = \int_0^1 (1-z) F'''(\phib + z(\phi - \phib)) \, \de z, \quad R_2^h = \int_0^1 (1-z) \hh''(\phib + z(\phi - \phib)) \, \de z,  \\
        R_3^h = \int_0^1 (1-z) P''(\phib + z(\phi - \phib)) \, \de z
	\end{gather*}
    Next, we observe that, since hypothesis \ref{ass:fconv} holds with $s \in [2,4]$, $F'''$ grows at most linearly, thus we can estimate:
    \begin{align*}
        \norm{R_1^h}_{\Lx\infty} & \le \int_0^1 \abs{1 - z} \norm{F'''(\phib + z (\phi - \phib)}_{\Lx\infty} \, \de z \\
        & \le \left( \int_0^1 \abs{1 - z} \, \de z \right) \left( 2 \norm{F'''(\phib)}_{\Lx\infty} + \norm{F'''(\phi)}_{\Lx\infty} \right) \\
        & \le C \left( 1 + \norm{\phib}_{\Lx\infty} + \norm{\phi}_{\Lx\infty} \right).
    \end{align*}
    Therefore, since $\phib, \phi$ are uniformly bounded in $\LT \infty V \cap \LT 2 {\Hx3} \hookrightarrow \LT 8 {\Lx\infty}$ (cf.~\eqref{eq:emb_l8linf}) and \eqref{hk:bound} holds, we can surely conclude that 
    \begin{equation}
        \label{eq:frechet:r1h}
        \int_0^T \norm{R_1^h}^8_{\Lx\infty} \, \de t \le C \left( 1 + \int_0^T \norm{\phib}^8_{\Lx\infty} \, \de t + \int_0^T \norm{\phi}^8_{\Lx\infty} \, \de t \right) \le C_\Lambda,
    \end{equation}
    meaning that $R^h_1$ is uniformly bounded in $\LT 8 {\Lx\infty}$.
    Moreover, since $P, \hh \in W^{2,\infty}(\R)$ by \ref{ass:ph_second}, we can also easily infer that
    \begin{equation}
        \label{eq:frechet:r23h}
        \norm{R^h_2}_{\Lqt\infty} \le C \quad \text{and} \quad \norm{R^h_3}_{\Lqt\infty} \le C.
    \end{equation}
    
	To show \eqref{frechet:aim}, we now proceed by performing a priori estimates on the system \eqref{eq:psi}--\eqref{icd}. 
    Indeed, the main estimate is done by testing \eqref{eq:psi} by $\psi$, \eqref{eq:eta} by $\Delta \psi$ and \eqref{eq:theta} by $\theta$ and summing them up to obtain:
	\begin{equation}
    \label{eq:frechet:mainest}
    \begin{split}
		& \mezzo \ddt \norm{\psi}^2_H + \mezzo \ddt \norm{\theta}^2_H + \norm{\Delta \psi}^2_H + \norm{\nabla \theta}^2_H \\
        & \quad = (F^h, \Delta \psi)_H - 2\chi (\theta, \Delta \psi)_H + (Q^h, \psi - \theta)_H - (U^h,\psi)_H - \norm{\theta}^2_H.
    \end{split}
	\end{equation}
    We now proceed to estimate each term on the right-hand side of \eqref{eq:frechet:mainest}. 
    First, by using \ref{ass:fgrowth4}, \eqref{eq:frechet:r1h} and H\"older and Young's inequalities, we infer that 
    \begin{align*}
        (F^h, \Delta \psi)_H & = (F''(\phib) \psi, \Delta \psi)_H + (R^h_1 (\phi - \phib)^2, \Delta \psi)_H \\
        & \le \norm{F''(\phib)}_{\Lx\infty} \norm{\psi}_H \norm{\Delta \psi}_H + \norm{R^h_1}_\infty \norm{\phi - \phib}^2_{\Lx4} \norm{\Delta \psi}_H \\
        & \le \frac14 \norm{\Delta \psi}^2_H + C \norm{F''(\phib)}^2_{\Lx\infty} \norm{\psi}^2_H + C \norm{R^h_1}^2_\infty \norm{\phi - \phib}^4_{\Lx4},
    \end{align*}
    where $\norm{F''(\phi)}^2_{\Lx\infty} \in \Lx1$ by \eqref{eq:emb_l8linf} and \ref{ass:fconv}. 
    Moreover, after integrating on $(0,T)$, the last term in the above inequality can be handled by using Gagliardo--Nirenberg's inequality \eqref{gn:ineq} in space with $N=3$, $p=4$, $j=0$, $r=2$, $m=2$, $\alpha = 3/8$ and $q=2$ and then H\"older's inequality in time with $3/4 + 1/4 =1$ in the following way:
    \begin{align*}
        & \int_0^T \norm{R^h_1}^2_\infty \norm{\phi - \phib}^4_{\Lx4} \, \de t 
        \le C \int_0^T \norm{R^h_1}^2_\infty \norm{\phi - \phib}^{3/2}_W \norm{\phi - \phib}^{5/2}_H \\
        & \quad \le C \left( \int_0^T \norm{\phi - \phib}^2_W \, \de t \right)^{3/4} \left( \int_0^T \norm{R^h_1}^8_{\Lx\infty} \norm{\phi - \phib}^{10}_H \, \de t \right)^{1/4} \\
        & \quad \le C \underbrace{\norm{R^h_1}^{2}_{\LT 8 {\Lx\infty}}}_{\le \, C_\Lambda} \norm{\phi - \phib}^{3/2}_{\LT 2 W} \norm{\phi - \phib}^{5/2}_{\LT \infty H}.
    \end{align*}
    Hence, by \eqref{eq:frechet:r1h} and the continuous dependence estimate \eqref{eq:contdep_estimate}, we deduce that 
    \begin{equation}
        \label{eq:frechet:estrhs1}
        \int_0^T \norm{R^h_1}^2_\infty \norm{\phi - \phib}^4_{\Lx4} \, \de t \le C \left( \norm{h}^4_H + \norm{k}^4_H \right).
    \end{equation}
    Next, we continue estimating the terms on the right-hand side of \eqref{eq:frechet:mainest}. 
    The second one can be easily bounded by means of Cauchy--Schwarz and Young's inequalities, indeed
    \[ 
        2 \chi (\theta, \Delta \psi)_H \le \frac14 \norm{\Delta \psi}^2_H + C \norm{\theta}^2_H .
    \]
    Then, we pass on to the third term and we use again H\"older and Young's inequalities, together with \eqref{eq:frechet:r23h} and the Sobolev embeddings $W \hookrightarrow V$, $W \hookrightarrow \Lx\infty$ and $V \hookrightarrow \Lx6$, to deduce that 
    {\allowdisplaybreaks
    \begin{align*}
        (Q^h, \psi - \theta)_H & = (P(\phib)(\theta - \chi \psi - \zeta), \psi - \theta)_H 
        + (P'(\phib)(\sigmab + \chi (1 - \phib) - \mub) \psi, \psi - \theta)_H \\
        & \quad + ((P(\phi) - P(\phib))[(\sigma - \sigmab) - \chi (\phi - \phib) - (\mu - \mub)], \psi - \theta)_H \\
        & \quad + (R^h_3 (\sigmab + \chi (1 - \phib) - \mub) (\phi - \phib)^2, \psi - \theta)_H \\
        & \le P_\infty \norm{\theta - \chi \psi - \zeta}_H \norm{\psi - \theta}_H \\
        & \quad + P'_\infty \norm{\sigmab + \chi (1 - \phib) - \mub}_H \norm{\psi}_\infty \norm{\psi - \theta}_H \\
        & \quad + \norm{P(\phi) - P(\phib)}_H \norm{(\sigma - \sigmab) - \chi (\phi - \phib) - (\mu - \mub)}_{\Lx4} \norm{\psi - \theta}_{\Lx4} \\
        & \quad + \norm{R^h_3}_{\Lx\infty} \norm{\sigmab + \chi (1 - \phib) - \mub}_{\Lx6} \norm{\phi - \phib}^2_{\Lx6} \norm{\psi - \theta}_H \\
        & \le C \left( \norm{\theta}_H + \norm{\psi}_H + \norm{\zeta}_H \right) \left( \norm{\psi}_H + \norm{\theta}_H \right) \\
        & \quad + C \norm{\sigmab + \chi (1 - \phib) - \mub}_H \norm{\psi}_W \left( \norm{\psi}_H + \norm{\theta}_H \right) \\
        & \quad + C \norm{\phi - \phib}_H \left( \norm{\sigma - \sigmab}_V + \norm{\phi - \phib}_W + \norm{\mu - \mub}_V \right) \left( \norm{\psi}_W + \norm{\theta}_V \right) \\
        & \quad + C \norm{\sigmab + \chi (1 - \phib) - \mub}_V \norm{\phi - \phib}^2_V \left( \norm{\psi}_H + \norm{\theta}_H \right) \\
        & \le \frac18 \norm{\zeta}^2_H + \frac18 \norm{\psi}^2_W + \frac12 \norm{\theta}^2_V + C \left( 1 + \norm{\sigmab + \chi (1 - \phib) - \mub}^2_V \right) \left( \norm{\psi}^2_H + \norm{\theta}^2_H \right) \\
        & \quad + C \norm{\phi - \phib}^2_H \left( \norm{\sigma - \sigmab}^2_V + \norm{\phi - \phib}^2_W + \norm{\mu - \mub}^2_V \right) + C \norm{\phi - \phib}^4_V,
    \end{align*} }
    where $\norm{\sigmab + \chi (1 - \phib) - \mub}^2_V \in \Lt1$ by Theorem \ref{thm:weaksols}. 
    Additionally, after integration on $(0,T)$, the last two terms can be bound by using the continuous dependence estimate \eqref{eq:contdep_estimate} and Gagliardo--Nirenberg's inequality \eqref{gn:ineq} as follows:
    \begin{equation}
    \label{eq:frechet:estrhs2}
    \begin{split}
        & \int_0^T \norm{\phi - \phib}^2_H \left( \norm{\sigma - \sigmab}^2_V + \norm{\phi - \phib}^2_W + \norm{\mu - \mub}^2_V \right) \, \de t \\
        & \qquad \le \norm{\phi - \phib}^2_{\LT \infty H} \left( \norm{\sigma - \sigmab}^2_{\LT 2 V} + \norm{\phi - \phib}^2_{\LT 2 W} + \norm{\mu - \mub}^2_{\LT 2 V} \right) \\
        & \qquad \le C \left( \norm{h}^4_H + \norm{k}^4_H \right), \\
        & \int_0^T \norm{\phi - \phib}^4_V \, \de t \le \int_0^T \norm{\phi - \phib}^2_H \norm{\phi - \phib}^2_W \, \de t \\
        & \qquad \le  \norm{\phi - \phib}^2_{\LT \infty H} \norm{\phi - \phib}^2_{\LT 2 W} \le C \left( \norm{h}^4_H + \norm{k}^4_H \right).
    \end{split}
    \end{equation}
    Moreover, by comparison in \eqref{eq:zeta}, we can easily see that 
    \begin{equation}
        \label{eq:frechet:zeta}
        \norm{\zeta}^2_H \le \norm{\Delta \psi}^2_H + \norm{F''(\phib)}^2_{\Lx \infty} \norm{\psi}^2_H + \norm{R^h_1}^2_{\Lx\infty} \norm{\phi - \phib}^4_{\Lx4} + \chi^2 \norm{\theta}^2_H,
    \end{equation}
    where $\norm{F''(\phib)}^2_{\Lx\infty} \in \Lx1$ by \eqref{eq:emb_l8linf} and \ref{ass:fconv} and the term $\norm{R^h_1}^2_{\Lx\infty} \norm{\phi - \phib}^4_{\Lx4}$ can be treated as in \eqref{eq:frechet:estrhs1}.
    Finally, we can estimate the fourth term in \eqref{eq:frechet:mainest} with similar procedures, getting that
    \begin{align*}
        (U^h,  \psi)_H & = (c \hh'(\phib) \psi, \psi)_H + (c R^h_2 (\phi - \phib)^2, \psi)_H \\
        & \le c_\infty \hh'_\infty \norm{\psi}^2_H + c_\infty \norm{\R^h_2}_{\Lx\infty} \norm{\phi - \phib}^2_{\Lx4} \norm{\psi}_H \\
        & \le C \norm{\psi}^2_H + C \norm{\phi - \phib}^4_V.
    \end{align*}
    Hence, putting all together and integrating on $(0,t)$, for any $t \in (0,T)$, from \eqref{eq:frechet:mainest} we deduce that 
    \begin{align*}
        & \mezzo \norm{\psi(t)}^2_H + \mezzo \norm{\theta(t)}^2_H + \frac14 \int_0^t \norm{\Delta \psi}^2_H \, \de s + \mezzo \int_0^t \norm{\nabla \theta}^2_H \, \de s \\
        & \quad \le C \int_0^T \left( 1 + \norm{F''(\phib)}^2_{\Lx\infty} + \norm{\sigmab + \chi (1 - \phib) - \mub}^2_V \right) \norm{\psi}^2_H \, \de s \\
        & \qquad + C \int_0^T \left( 1 + \norm{\sigmab + \chi (1 - \phib) - \mub}^2_V \right) \norm{\theta}^2_H \, \de s \\
        & \qquad + C \int_0^T \norm{R^h_1}^2_{\Lx\infty} \norm{\phi - \phib}^4_{\Lx4} \, \de s + C \int_0^T \norm{\phi - \phib}^4_V \, \de s \\
        & \qquad + C \int_0^T \norm{\phi - \phib}^2_H \left( \norm{\sigma - \sigmab}^2_V + \norm{\phi - \phib}^2_W + \norm{\mu - \mub}^2_V \right) \, \de s.
    \end{align*}
    Therefore, by means of Gronwall's inequality, together with \eqref{eq:frechet:estrhs1} and \eqref{eq:frechet:estrhs2}, we conclude that
    \begin{equation}
        \label{eq:frechet:energyest}
        \norm{\psi}^2_{\LT \infty H \cap \LT 2 W} + \norm{\theta}^2_{\LT \infty H \cap \LT 2 V} \le C \left( \norm{h}^4_H + \norm{k}^4_H \right).
    \end{equation}
    Then, by \eqref{eq:frechet:zeta}, together with \eqref{eq:frechet:estrhs1} and \eqref{eq:frechet:energyest}, we also immediately infer that 
    \begin{equation}
        \label{eq:frechet:zetal2h}
        \norm{\zeta}^2_{\LT 2 H} \le C \left( \norm{h}^4_H + \norm{k}^4_H \right).
    \end{equation}
    Moreover, by comparison in \eqref{eq:psi} and \eqref{eq:theta}, starting from \eqref{eq:frechet:energyest} and \eqref{eq:frechet:zetal2h}, we additionally deduce that
    \begin{equation}
        \label{eq:frechet:timeest}
        \norm{\psi}^2_{\HT 1 {W^*}} + \norm{\theta}^2_{\HT 1 {V^*}} \le C \left( \norm{h}^4_H + \norm{k}^4_H \right).
    \end{equation}
    Then, due to the standard embeddings $\HT 1 {W^*} \cap \LT 2 W \hookrightarrow \C 0 H$ as well as $\HT 1 {V^*} \cap \LT 2 V \hookrightarrow \C 0 H$, estimates \eqref{eq:frechet:energyest} and \eqref{eq:frechet:timeest} imply that 
    \[
        \norm{(\psi(T), \theta(T))}^2_{H \times H} \le C \norm{(h,k)}^4_{H \times H},
    \]
    which subsequently implies \eqref{frechet:aim} with $q = 4 \gs 2$.
    This concludes the proof of Theorem \ref{thm:frechet}.
\end{proof}

\subsection{Adjoint system and necessary optimality conditions}

We now introduce the adjoint system in order to deduce the necessary optimality conditions. 
We fix an optimal state $(\phib, \mub, \sigmab)$ corresponding to an optimal control $(\phiob, \sigmaob) \in \Uad \times \Vad$, then the associated adjoint variables $(p,q,r)$ formally satisfy the following system:
\begin{alignat}{2}
	& - \partial_t p - \Delta q + F''(\phib) q + \chi \Delta r + \chi P(\phib) (p - r) \nonumber \\ 
	& \qquad - P'(\phib) (\sigmab + \chi (1 - \phib) - \mub)(p-r) + c \hh'(\phib) p = 0 
	\qquad && \hbox{in  $Q_T$},  \label{eq:p}\\
	& - q - \Delta p + P(\phib)(p-r) = 0 
	\qquad && \hbox{in  $Q_T$},  \label{eq:q} \\
	& - \partial_t r - \Delta r - \chi q - P(\phib) (p-r) + r = 0 
	\qquad && \hbox{in  $Q_T$}, \label{eq:r} \\
	& \partial_{\n} p = \partial_{\n} q = \partial_{\n} r = 0 
	\qquad && \hbox{on  $\Sigma_T$}, \label{bca} \\
	& p(T) = \lambda_1 (\phib(T) - \phimeas), \quad r(T) = \lambda_2 (\sigmab(T) - \sigmameas) 
	\qquad && \hbox{in  $\Omega$}. \label{fca}
\end{alignat}

\noindent
First, we state a well-posedness result for the adjoint system.

\begin{theorem}
	\label{thm:adjoint}
	Assume hypotheses \emph{\ref{ass:setting}--\ref{ass:hc}}, \emph{\ref{ass:phinf}} and \emph{\ref{ass:C1}--\ref{ass:fgrowth4}}. Let $(\phib, \mub, \sigmab)$ be the weak solution to \eqref{eq:phi}--\eqref{ic}, corresponding to $(\phiob, \sigmaob) \in \Uad \times \Vad$. Then, the adjoint system \eqref{eq:p}--\eqref{fca} admits a unique weak solution such that
	\begin{align*}
		& p \in H^1(0,T;W^*) \cap \mathcal{C}^0([0,T]; H) \cap L^2(0,T;W), \\
		& q \in L^2(0,T;H), \\
		& r \in H^1(0,T;V^*) \cap \mathcal{C}^0([0,T]; H) \cap L^2(0,T;V),
	\end{align*}
	which fulfils \eqref{eq:p}--\eqref{fca} in variational formulation, i.e.~it satisfies
	\begin{align}
		& \duality{-\partial_t p, w}_W - (q, \Delta w)_H + (F''(\phib) q, w)_H + \chi (\nabla r, \nabla w)_H + \chi (P(\phib) (p - r), w)_H \nonumber \\
		& \qquad - (P'(\phib) (\sigmab + \chi (1 - \phib) - \mub)(p-r), w)_H + (c \hh'(\phib) p, w)_H = 0, \label{varform:p} \\
		& (q,w)_H = - (\Delta p, w)_H + (P(\phib)(p-r),w)_H, \label{varform:q} \\
		& \duality{- \partial_t r, v}_V + (\nabla r, \nabla v)_H - \chi (q, v)_H - (P(\phib)(p-r),v)_H + (r, v)_H = 0, \label{varform:r}
	\end{align}
	for a.e.~$t \in (0,T)$ and for any $(w,v) \in W \times V$, and the final conditions $p(T) = \lambda_1 (\phib(T) - \phimeas)$, $r(T) = \lambda_2 (\sigmab(T) - \sigmameas)$.
\end{theorem}

\begin{proof}
    Since we are dealing with a backward \emph{linear} system, we proceed with formal estimates that can be made rigorous through a Faedo--Galerkin scheme. 
    For the main estimate, we choose $w = p$ in \eqref{varform:p}, $w = \Delta p$ in \eqref{varform:q}, $v = r$ in \eqref{varform:r} and sum them up to obtain
    \begin{equation}
        \label{eq:adjoint:mainest}
        \begin{split}
        & - \mezzo \ddt \norm{p}^2_H - \mezzo \ddt \norm{r}^2_H + \norm{\Delta p}^2_H  + \norm{\nabla r}^2_H \\
        & \quad = - (F''(\phib) q, p)_H - \chi (\nabla r, \nabla p)_H + \chi (q, r)_H + (P(\phib)(p - r), - \chi p + \Delta p + r)_H \\
        & \qquad + (P'(\phib)(\sigmab + \chi  (1 - \phib) - \mub)(p - r), p)_H  - (c \hh'(\phib) p, p)_H - \norm{r}^2_H.
        \end{split}
    \end{equation}
    We begin by observing that, by comparison in \eqref{eq:q}, it follows that 
    \begin{equation}
    \label{eq:adjoint:q}
        \norm{q}^2_H \le \norm{\Delta p}^2_H +  \norm{P(\phib)(p - r)}^2_H \le \norm{\Delta p}^2_H + C \norm{p}^2_H + C \norm{r}^2_H,
    \end{equation}
    where we used the boundedness of $P$. 
    Next, we can start estimating the terms on the right-hand side of \eqref{eq:adjoint:mainest}. 
    Indeed, by using H\"older and Young's inequalities, together with \ref{ass:ph_second} and the Sobolev embedding $W \hookrightarrow \Lx\infty$, we can see that the following estimates hold:
    {\allowdisplaybreaks
    \begin{align*}
        & (F''(\phib) q, p)_H \le \norm{F''(\phib)}_{\Lx\infty} \norm{q}_H \norm{p}_H \le \frac18 \norm{q}^2_H + C \norm{F''(\phib)}^2_{\Lx\infty} \norm{p}^2_H \\
        & \qquad \le \frac18 \norm{\Delta p}^2_H + C \left(1 + \norm{F''(\phib)}^2_{\Lx\infty} \right) \norm{p}^2_H + C \norm{r}^2_H, \\
        & \chi (\nabla r, \nabla p)_H = - \chi (r, \Delta p)_H \le \frac18 \norm{\Delta p}^2_H + C \norm{r}^2_H, \\
        & \chi (q,r)_H \le \frac18 \norm{q}^2_H + C \norm{r}^2_H \le \frac18 \norm{\Delta p}^2_H + C \norm{p}^2_H + C \norm{r}^2_H, \\
        & (P(\phib)(p - r), -\chi p + \Delta p + r)_H \le \frac18 \norm{\Delta p}^2_H + C \norm{p}^2_H + C \norm{r}^2_H, \\
        & (P'(\phib)(\sigmab + \chi (1 - \phib) - \mub)(p - r), p)_H \le P'_\infty \norm{\sigmab + \chi (1 - \phib) - \mub}_H ( \norm{p}_H + \norm{r}_H ) \norm{p}_{\Lx\infty} \\
        & \qquad \le \frac14 \norm{p}^2_W + C \norm{\sigmab + \chi (1 - \phib) - \mub}^2_H ( \norm{p}^2_H + \norm{r}^2_H ), \\
        & (c \hh'(\phib) p, p)_H \le c_\infty \hh'_\infty \norm{p}^2_H, 
    \end{align*}
    }
    where $\norm{F''(\phib)}^2_{\Lx\infty} \in \Lx1$ by \eqref{eq:emb_l8linf} and \ref{ass:fconv} and $\norm{\sigmab + \chi (1 - \phib) - \mub}^2_H \in \Lx1$ by Theorem \ref{thm:weaksols}. 
    Then, by gathering all estimates and integrating on $(t, T)$, for any $t \in (0,T)$, we infer that 
    \begin{align*}
        & \mezzo \norm{p(t)}^2_H + \mezzo \norm{r(t)}^2_H + \frac14 \int_t^T \norm{\Delta p}^2_H \, \de s + \int_t^T \norm{\nabla r}^2_H \, \de s \\
        & \quad \le \norm{\lambda_1 (\phib(T) - \phimeas)}^2_H + \norm{\lambda_2 (\sigmab(T) - \sigmameas)}^2_H \\
        & \qquad + C \int_0^T (1 + \norm{\sigmab + \chi (1 - \phib) - \mub}^2_H) (\norm{p}^2_H + \norm{r}^2_H) \, \de s. 
    \end{align*}
    Hence, by applying Gronwall's inequality, we conclude that 
    \begin{equation}
        \label{eq:adjoint:energyest}
        \begin{split}
        & \norm{p}^2_{\LT \infty H \cap \LT 2 W} + \norm{r}^2_{\LT \infty H \cap \LT 2 V} \\
        & \quad \le C \left( \norm{\lambda_1 (\phib(T) - \phimeas)}^2_H + \norm{\lambda_2 (\sigmab(T) - \sigmameas)}^2_H \right).
        \end{split}
    \end{equation}
    Moreover, by \eqref{eq:adjoint:q} and \eqref{eq:adjoint:energyest}, we also easily deduce that
    \begin{equation}
        \label{eq:adjoint:ql2h}
        \norm{q}^2_{\LT 2 H} \le C \left( \norm{\lambda_1 (\phib(T) - \phimeas)}^2_H + \norm{\lambda_2 (\sigmab(T) - \sigmameas)}^2_H \right).
    \end{equation}
    Finally, by exploiting the boundedness of $P$ and $\hh$ by \ref{ass:ph_second} and arguing by comparison in \eqref{eq:p} and \eqref{eq:r}, starting from the estimates proved in \eqref{eq:adjoint:energyest} and \eqref{eq:adjoint:ql2h}, we also infer that 
    \begin{equation}
        \label{eq:adjoint:timeest}
        \norm{p}^2_{\HT 1 {W^*}} + \norm{r}^2_{\HT 1 {V^*}} \le C \left( \norm{\lambda_1 (\phib(T) - \phimeas)}^2_H + \norm{\lambda_2 (\sigmab(T) - \sigmameas)}^2_H \right).
    \end{equation}
    Indeed, the only non-trivial terms to treat are the following ones:
    {\allowdisplaybreaks
    \begin{align*}
        & \norm{F''(\phib) q}^2_{\LT 2 {W^*}} = \int_0^T \norm{F''(\phib) q}^2_{W^*} \, \de t = \int_0^T \left( \sup_{\norm{w}_W = 1} (F''(\phib) q, w)_H \right)^2 \, \de t \\
        & \qquad \le \int_0^T \left( \sup_{\norm{w}_W = 1} \norm{F''(\phib)}_H \norm{q}_H \norm{w}_{\Lx\infty} \right)^2 \, \de t \\
        & \qquad \le C \int_0^T \norm{F''(\phib)}^2_H \norm{q}^2_H \, \de t \le C \norm{F''(\phib)}^2_{\LT \infty H} \norm{q}^2_{\LT 2 H} \le C \norm{q}^2_{\LT 2 H}, \\
        & \norm{P'(\phib)(\sigmab + \chi (1 - \phib) - \mub)(p - r)}^2_{\LT 2 {W^*}} = \int_0^T \norm{P'(\phib)(\sigmab + \chi (1 - \phib) - \mub)(p - r)}^2_{W^*} \, \de t \\
        & \qquad = \int_0^T \left( \sup_{\norm{w}_W = 1} (P'(\phib)(\sigmab + \chi (1 - \phib) - \mub)(p - r), w)_H \right)^2 \, \de t \\
        & \qquad \le \int_0^T \left( P'_\infty \norm{\sigmab + \chi (1 - \phib) - \mub}_H \norm{p - r}_H \norm{w}_{\Lx\infty} \right)^2 \, \de t \\
        & \qquad \le C \norm{\sigmab + \chi (1 - \phib) - \mub}^2_{\LT 2 H} \norm{p - r}^2_{\LT \infty H} \le C \left( \norm{p}^2_{\LT \infty H} + \norm{r}^2_{\LT \infty H} \right), 
    \end{align*}
    }
    where $F''(\phib) \in \LT \infty H$ by \ref{ass:fgrowth4} and the fact that $\phib$ is uniformly bounded in $\LT \infty {\Lx6}$ by Theorem \ref{thm:weaksols} and the Sobolev embedding $V \hookrightarrow \Lx6$, while $\sigma + \chi (1 - \phi) - \mu \in \LT 2 H$ by Theorem \ref{thm:weaksols}. 

    Then, starting from the uniform estimates \eqref{eq:adjoint:energyest}, \eqref{eq:adjoint:ql2h} and \eqref{eq:adjoint:timeest}, one can easily pass to the limit in the Galerkin discretisation and, due to the linearity of the system, prove the uniqueness of the solution. 
    This concludes the proof of Theorem \ref{thm:adjoint}.
\end{proof}

\begin{remark}
    Note that, due to the regularity of the solution to the adjoint system \eqref{eq:p}--\eqref{fca}, the initial value $(p(0), r(0))$ is well-defined in $H \times H$. 
\end{remark}

Finally, we can derive the first-order necessary optimality conditions for the constrained minimisation problem.

\begin{theorem}
	\label{thm:optcond}
	Assume hypotheses \emph{\ref{ass:setting}--\ref{ass:hc}} and \emph{\ref{ass:C1}--\ref{ass:ph_second}}. 
	Let $(\phiob, \sigmaob) \in \Uad \times \Vad$ be an optimal pair for \emph{(MP)} and let $(\phib, \mub, \sigmab)$ be the corresponding optimal state, i.e.~the weak solution to \eqref{eq:phi}--\eqref{ic} with such $(\phiob, \sigmaob)$. 
	Let also $(p,q,r)$ be the adjoint variables to $(\phib, \sigmab, \mub)$, i.e.~the solutions to the adjoint system \eqref{eq:p}--\eqref{fca}. 
	Then, they satisfy the following variational inequality:
	\begin{equation}
		\label{var:ineq}
        \begin{split}
		& \int_\Omega p(0) (\phi_0 - \phiob) \, \de x + \int_\Omega \alpha_1 F'(\phiob) (\phi_0 - \phiob) + \alpha_1 \nabla \phiob \cdot \nabla (\phi_0 - \phiob) \, \de x  \\ 
        & \quad + \int_\Omega r(0) (\sigma_0 - \sigmaob) \, \de x + \int_\Omega \alpha_2 \sigmaob (\sigma_0 - \sigmaob) \, \de x \ge 0, \quad \hbox{for any $(\phi_0, \sigma_0) \in \Uad \times \Vad$.} 
        \end{split}
	\end{equation}
\end{theorem}

\begin{proof}
    We recall the notation used in Theorem \ref{thm:control_existence}, namely that we call $\Jtil(\phi(T), \sigma(T), \phi_0, \sigma_0)$ the funtional such that $\Jcal(\phi_0, \sigma_0) = \Jtil(\Rcal(\phi_0, \sigma_0), \phi_0, \sigma_0)$.
	First, observe that the cost functional $\Jtil(\phi(T), \sigma(T), \phi_0, \sigma_0)$ is Fr\'echet differentiable as a functional from the space $H \times H \times V \times H$ to $\R$. 
    Indeed, most of the terms are quadratic, so it is easy to see that they are differentiable in $\Lx2$. 
    The only non-trivial term is the one related to the Ginzburg--Landau energy 
    \[
        E(\phi_0) = \int_\Omega F(\phi_0) + \mezzo \abs{\nabla \phi_0}^2 \, \de x.
    \]
    However, it is a standard matter to prove that this is also Fr\'echet differentiable in $V$, with derivative given by
    \[
        \D E (\phi_0) [v] = \int_\Omega F'(\phi_0) v + \nabla \phi_0 \cdot \nabla v \, \de x \quad \hbox{for any $ v \in V$,}
    \]
    due to the growth of $F$ up to the power $6$ by \ref{ass:fconv} and the Sobolev embedding $V \hookrightarrow \Lx6$.
    Next, in Theorem \ref{thm:frechet} we showed that the solution operator $\Rcal(\phi_0, \sigma_0) = (\phi(T), \sigma(T))$ is Fr\'echet differentiable from $V \times H$ to $H \times H$.
    Thus, the functional $\Jcal: V \times H \to \R$ such that
    \[
        \Jcal(\phi_0, \sigma_0) := \Jtil(\Rcal(\phi_0, \sigma_0), \phi_0, \sigma_0)
    \]
    is also Fr\'echet differentiable by the chain rule.

    Then, since $\Uad \times \Vad$ is closed and convex and $\Jcal$ is Fr\'echet differentiable, an optimal pair $(\phiob, \sigmaob)$ has to satisfy the necessary optimality condition
    \[
        \D\Jcal(\phiob, \sigmaob) [(\phi_0 - \phiob, \sigma_0 - \sigmaob)] \ge 0 \quad \hbox{for any $(\phi_0, \sigma_0) \in \Uad \times \Vad$.}
    \]
    By computing explicitly the derivative of $\Jcal$, we equivalently get that for any $(\phi_0, \sigma_0) \in \Uad \times \Vad$
    \begin{align*}
        & \int_\Omega \lambda_1 (\phib(T) - \phimeas) \xi(T) \, \de x + \int_\Omega \lambda_2 (\sigmab(T) - \sigmameas) \rho(T) \, \de x \\
        & \quad + \alpha_1 \int_\Omega F'(\phiob) (\phi_0 - \phiob) + \nabla \phiob \cdot \nabla(\phi_0 - \phiob) \, \de x + \alpha_2 \int_\Omega \sigmaob (\sigma_0 - \sigmaob) \, \de x \ge 0, 
    \end{align*}
    where $(\xi(T), \rho(T)) = \D \Rcal (\phiob, \sigmaob) [(\phi_0 - \phiob, \sigma_0 - \sigmaob)]$ are the components of the solution to the linearised system \eqref{eq:xi}--\eqref{icl} corresponding to $h = \phi_0 - \phiob$ and $k = \sigma_0 - \sigmaob$.

    Now, to show \eqref{var:ineq}, we just need to prove that the following identity holds
    \begin{align*}
        & \int_\Omega \lambda_1 (\phib(T) - \phimeas) \xi(T) \, \de x + \int_\Omega \lambda_2 (\sigmab(T) - \sigmameas) \rho(T) \, \de x \\
        & \quad = \int_\Omega p(0) (\phi_0 - \phiob) \, \de x + \int_\Omega r(0) (\sigma_0 - \sigmaob) \, \de x.
    \end{align*}
    To do this, we consider the weak formulations of the adjoint system \eqref{varform:p}--\eqref{varform:r} and of the linearised system \eqref{varform:xi}--\eqref{varform:rho}. 
    Indeed, we test \eqref{varform:p} by $\xi$, \eqref{varform:q} by $\eta$, \eqref{varform:r} by $\rho$, sum them up and integrate on $(0,T)$ to obtain:
    \begin{align*}
        & \int_0^T \Big[ \duality{- p_t, \xi}_W + \int_\Omega \Big( -q \Delta \xi + F''(\phib) q, \xi + \chi \nabla r \cdot \nabla \xi + \chi P(\phib) (p - r) \xi \\
        & \qquad - P'(\phib)(\sigmab + \chi (1 - \phib) - \mub) (p - r) \xi + c \hh'(\phib) p \xi \, \de x \Big) \Big] \, \de t \\
        & \quad + \int_0^T \int_\Omega - q \eta - \Delta p \eta + P(\phib)(p - r) \eta \, \de x \, \de t \\
        & \quad + \int_0^T \Big[ \duality{- r_t, \rho}_V + \int_\Omega \Big( \nabla r \cdot \nabla \rho - \chi q \rho - P(\phib) (p - r) \rho + r \rho \Big) \, \de x \Big] \, \de t = 0.
    \end{align*}
    We now integrate by parts in time and space and, after regrouping some terms, we deduce the following:
    \begin{align*}
        & \int_0^T \Big[ \duality{\xi_t, p}_W + \int_\Omega \Big( - \eta \Delta p - P(\phib) (\rho - \chi \xi - \eta) p \\
        & \qquad + P'(\phib)(\sigmab + \chi (1 - \phib) - \mub) \xi p + c \hh'(\phib) \xi p \, \de x \Big) \Big] \, \de t \\
        & \quad + \int_0^T \int_\Omega \left( - \eta - \Delta \xi + F''(\phib) \xi - \chi \rho \right) q \, \de x \, \de t \\
        & \quad + \int_0^T \Big[ \duality{\rho_t, r}_V + \int_\Omega \Big( \nabla \rho \cdot \nabla r + \chi \nabla \xi \cdot \nabla r + P(\phib) (\rho - \chi \xi - \eta) r \\
        & \qquad + P'(\phib)(\sigmab + \chi (1 - \phib) - \mub) \xi r - \rho r \Big) \, \de x \Big] \, \de t \\
        & \quad + \int_\Omega p(T) \xi(T) - p(0) \xi(0) \, \de x + \int_\Omega r(T) \rho(T) - r(0) \rho(0) \, \de x = 0.
    \end{align*}
    In the first three integrals, we now recognise the weak formulations of the linearised system \eqref{varform:xi}--\eqref{varform:rho}, tested respectively by $p$, $q$ and $r$. 
    Then, such terms vanish and we are only left with
    \[ 
        \int_\Omega p(T) \xi(T) + r(T) \rho(T) \, \de x = \int_\Omega p(0) \xi(0) + r(0) \rho(0) \, \de x.
    \]
    Hence, by recalling the final conditions \eqref{fca} on the adjoint system and the initial conditions \eqref{icl}, with $h = \phi_0 - \phiob$ and $k = \sigma_0 - \sigmaob$, on the linearised system, we finally infer that 
    \begin{align*}
        & \int_\Omega \lambda_1 (\phib(T) - \phimeas) \xi(T) \, \de x + \int_\Omega \lambda_2 (\sigmab(T) - \sigmameas) \rho(T) \, \de x \\
        & \quad = \int_\Omega p(0) (\phi_0 - \phiob) \, \de x + \int_\Omega r(0) (\sigma_0 - \sigmaob) \, \de x.
    \end{align*}
    This concludes the proof of Theorem \ref{thm:optcond}.
\end{proof}

\begin{remark}
    Through a discrete version of the variational inequality \eqref{var:ineq}, we can then approximate the solution to the Tikhonov-regularised inverse problem through gradient descent methods. 
    This is the starting point of Sections \ref{sec:numerics} and \ref{sec:simulations} on numerical methods and experiments.
\end{remark}

\section{Numerical methods}
\label{sec:numerics}

In this section, we introduce the finite element and time discretisations of \eqref{eq:phi}-\eqref{eq:sigma} and of the constrained minimisation problem \eqref{oc:problem} 
We recall that in Section \ref{sec:tikhonov} we found the optimality conditions for a Tikhonov-regularised problem aiming to reconstruct the initial data in $V \times H$.
This may be in contrast with the fact that in Section \ref{sec:uniq}, we showed that the inverse problem of reconstructing the initial data for the operator $\Rcal$ has a unique solution in the more regular space $W \times V$. 
However, such a choice is largely motivated by the numerical implementation carried out here.
Indeed, using only $C^0$ conforming finite elements with first-order basis functions reduces the computational cost of the algorithms, while still giving good reconstruction results. Reconstructing initial data $\phi_0$ in $W$ would require the employment of $C^1$ conforming finite elements to approximate the space $W$, which are non-standard and computationally demanding.

We now start introducing the discretised problem.
We adopt here the ``first discretise then optimise'' strategy to numerically solve the optimisation problem, which directly preserves the structure which is inherent in the infinite-dimensional optimisation problem, in contrast with straightforward ``first optimise then discretise'' strategies (see e.g. \cite[Chapter 3]{Hinzebook}).
Let $\mathcal{T}_{h}$ be a quasi-uniform conforming decomposition of $\Omega$ into $d-$simplices $K$, where d=2,3,
and let us define the following finite element space:
\begin{align*}
& S^{h} := \{\chi \in C(\bar{\Omega}):\chi |_{K}\in P^{1}(K) \; \forall K\in \mathcal{T}_{h}\}\subset W^{1,\infty}(\Omega),
\end{align*}
where $\mathbb{P}_{1}(K)$ indicates the space of polynomials of total order one on $K$.
We set $\Delta t = T/N$ for a $N \in \mathbb{N}$ and $t_{n}=n\Delta t$, $n=1,...,N$.
Starting from given data $\phi_{0}\in \Uad$ and $\sigma_0\in \Vad$, we set $\phi_{h}^{0}={P}_{V}^{h}(\phi_{0})\cap \Uad$ and $\sigma_h^0={P}_{H}^{h}(\sigma_{0})\cap \Vad$, where ${P}_{V}^{h}$ and $P_{H}^{h}$ are projection operators from $V$ to $S^h$ and from $H$ to $S^h$ respectively. We consider the following fully discretised problem:

\medskip
\noindent
\textbf{Problem $\mathbf{P}^{h}$.}
For $n=1,\dots,N$, given $(\phi_{h}^{n-1},\sigma_{h}^{n-1})\in S^{h}\times S^h$, find $(\phi_{h}^{n},\mu_{h}^{n},\sigma_{h}^{n})\in S^{h}\times S^{h}\times S^{h}$ such that for all $(\psi,\xi,\eta)\in S^h\times S^{h}\times S^{h}$,
\begin{alignat}{2}
\label{eqn:ph1}
\displaystyle & \notag \left(\frac{\phi_{h}^n-\phi_{h}^{n-1}}{\Delta t},\psi\right)+D_{\phi}(\nabla \mu_{h}^n,\nabla \psi)=\delta P_0 \left(P(\phi_{h}^{n-1}) \left(\frac{1}{\delta}\sigma_{h}^n + \chi (1-\phi_{h}^n) - \mu_{h}^n\right),\psi\right) \\
\displaystyle & \qquad - \left(c(x,t) \hh (\phi_{h}^{n-1}),\psi\right),
\\ \notag \\
\label{eqn:ph2} \displaystyle & (\mu_{h}^n,\xi)=\eps^2(\nabla \phi_{h}^n,\nabla \xi)+\Gamma \left(F_+'(\phi_{h}^n)+F_-'(\phi_{h}^{n-1}),\xi\right)-\chi (\sigma_h^n,\xi),\\ \notag \\
\displaystyle & \notag \left(\frac{\sigma_h^n-\sigma_h^{n-1}}{\Delta t},\eta\right)+\frac{D_{\sigma}}{\delta}(\nabla \sigma_h^n,\nabla \eta)-D_{\sigma}\chi(\nabla \phi_{h}^n,\nabla \eta) =\\
\label{eqn:ph3} & \displaystyle \qquad -\delta P_0 \left(P(\phi_{h}^{n-1}) \left(\frac{1}{\delta}\sigma_{h}^n + \chi (1-\phi_{h}^n) - \mu_{h}^n\right),\eta\right) + \left(\kappa (1-\sigma_h^n),\eta\right),
\end{alignat}
\Accorpa\eqnph {eqn:ph1} {eqn:ph3}
where we have reintroduced all the physical parameters as in {\eqnhd} and where
\[
F_+'(r):=4r^3+2r-6(\min(0,r))^2, \;\; F_-'(r):=-6(\max(0,r))^2, \quad \forall r\in \mathbb{R}.
\]


\begin{remark}
    We observe that in {\eqnph} we have reintroduced the physical parameters of the problem, since we will take physically viable values for them in the numerical simulations.
\end{remark}

We next introduce the discrete analogue to the minimisation problem (MP), 
\[\text{(MP)}_h:\;\;\min_{(\phi_h^0,\sigma_h^0)\in \Uad^h\times \Uad^h}	\mathcal{J}(\phi_h^N, \phi_h^0,\sigma_h^N, \sigma_h^0),\] 
where
\begin{align}
\label{eqn:jdisc}
	\notag \mathcal{J}(\phi_h^N, \phi_h^0,\sigma_h^N, \sigma_h^0)
		 = & \frac{\lambda_1}{2} \int_{\Omega} |\phi_h^N - \phimeas|^2 \,\de x + {\alpha_1} \int_{\Omega} \left(\Gamma\,F(\phi_h^0) + \frac{\eps^2}{2}|\nabla \phi_h^0|^2 \right)\,\de x\\
  & + \frac{\lambda_2}{2} \int_{\Omega} |\sigma_h^N - \sigmameas|^2 \,\de x + {\frac{\alpha_2}{2}} \int_{\Omega} | \sigma_h^0|^2 \,\de x,
\end{align}
subject to the constraints
\begin{equation*} 
	\phi_h^0,\sigma_h^0 \in \Uad^h := \left\{ \chi_h \in S^h \mid 0 \le \chi_h \le 1\right\},
\end{equation*}
and to the discretised state system \eqnph. To state the first order optimality conditions associated with (MP)$_h$ we introduce the adjoint variables $(p_h^n,q_h^n,r_h^n)\in S^h\times S^h\times S^h$, for $n=1,\dots,N$, and define the Lagrangian 
\[\mathcal{L}_h:\Uad^h\times \Uad^h\times S^h\times S^h\times S^h\times S^h\times S^h\times S^h\to \mathbb{R}\]
as:
\begin{align*}
    &\mathcal{L}_h\left(\phi_h^0,\sigma_h^0,(\phi_h^n)_{n=1}^N,(\mu_h^n)_{n=1}^N,(\sigma_h^n)_{n=1}^N,(p_h^n)_{n=1}^N,(q_h^n)_{n=1}^N,(r_h^n)_{n=1}^N\right)\\
    & \quad :=\frac{\lambda_1}{2} \norm{\phi_h^N - \phimeas}_H^2+{\alpha_1} \int_{\Omega} \left(\Gamma\,F(\phi_h^0) + \frac{\eps^2}{2}|\nabla \phi_h^0|^2 \right)\,\de x\\
    & \qquad
    +\frac{\lambda_2}{2} \norm{\sigma_h^N - \sigmameas}_H^2+{\frac{\alpha_2}{2}}\norm{\sigma_h^0}_H^2
    \\
    & \qquad -\Delta t \sum_{n=1}^N\biggl[\left(\frac{\phi_{h}^n-\phi_{h}^{n-1}}{\Delta t},p_h^n\right)+D_{\phi}(\nabla \mu_{h}^n,\nabla p_h^n)\\
    & \qquad \qquad -\delta P_0\left(P(\phi_{h}^{n-1}) \left(\frac{1}{\delta}\sigma_{h}^n + \chi (1-\phi_{h}^n) - \mu_{h}^n\right),p_h^n\right)+\left(c\hh (\phi_{h}^{n-1}),p_h^n\right)\biggr]\\
    & \qquad -\Delta t \sum_{n=1}^N\left[\eps^2(\nabla \phi_{h}^n,\nabla q_h^n)+\left(\Gamma\left(F_+'(\phi_{h}^n)+F_-'(\phi_{h}^{n-1})\right)-\chi\sigma_h^n-\mu_{h}^n,q_h^n\right)\right]\\
    & \qquad -\Delta t \sum_{n=1}^N\biggl[\left(\frac{\sigma_h^n-\sigma_h^{n-1}}{\Delta t},r_h^n\right)+\frac{D_{\sigma}}{\delta}(\nabla \sigma_h^n,\nabla r_h^n)-D_{\sigma}\chi(\nabla \phi_{h}^n,\nabla r_h^n)\\
    &\qquad \qquad+\delta P_0 \left(P(\phi_{h}^{n-1}) \left(\frac{1}{\delta}\sigma_{h}^n + \chi (1-\phi_{h}^n) - \mu_{h}^n\right),r_h^n\right) - \kappa \left((1-\sigma_h^n),r_h^n\right)\biggr].
\end{align*}
The first-order optimality conditions take the form
\[
(\D\mathcal{L}_h(\vec{x}_h),\tilde{\vec{x}}_h-\vec{x}_h)\geq 0,
\]
for any admissible direction $\tilde{\vec{x}}_h\in S^h$, where $\vec{x}_h$ stands for the vector of all the arguments of $\mathcal{L}_h$. We observe that for all components of $\vec{x}_h$ except for $\phi_h^0$ and $\sigma_h^0$ the optimality conditions reduce to the equalities $(\D\mathcal{L}_h(\vec{x}_h),\tilde{\vec{x}}_h)=0$, since in those cases $\tilde{\vec{x}}_h\in S^h$ without further constraints. In the case of the arguments $\phi_h^0$ and $\sigma_h^0$, we have $\tilde{\vec{x}}_h\in \Uad^h$, which is a convex subset of $S^h$, and the associated optimality condition takes the form of a variational inequality.
\\
\noindent
Proceeding explicitly with the Lagrangian calculus, we obtain the following discrete optimality systems.\\ \\
\noindent
\textbf{Problem $\mathbf{P}^{h}(\phi_h^0,\sigma_h^0)$: Direct problem.}\\
\noindent For $n=1,\dots,N$ and given $\phi_h^0,\sigma_h^0\in \Uad^h$, the derivatives $(\D_{p_h^m}\mathcal{L}_h(\vec{x}_h),\psi)=0$, $(\D_{q_h^m}\mathcal{L}_h(\vec{x}_h),\xi)=0$ and $(\D_{r_h^m}\mathcal{L}_h(\vec{x}_h),\eta)=0$, along directions $(\psi,\xi,\eta)\in S^h\times S^h\times S^h$, give the direct problem {\eqnph} with initial conditions $\phi_h^0$ and $\sigma_h^0$.\\ \\
\noindent
\textbf{Problem $\mathbf{Q}_1^{h}(\phi_h^N,\phi_h^{N-1},\sigma_h^N)$: Initialisation.} \\
\noindent
The derivatives $(\D_{\phi_h^N}\mathcal{L}_h(\vec{x}_h),u)=0$, $(\D_{\mu_h^N}\mathcal{L}_h(\vec{x}_h),v)=0$ and $(\D_{\sigma_h^N}\mathcal{L}_h(\vec{x}_h),w)=0$, along directions $(u,v,w)\in S^h\times S^h\times S^h$, give the initialisation problem: find $(p_h^N,q_h^N,r_h^N)\in S^h\times S^h\times S^h$ such that for all $(u,v,w)\in S^h\times S^h\times S^h$,
\begin{alignat}{2}
\label{eqn:qh1}
\displaystyle & \notag \left(p_h^N,u\right)+\Delta t\,\eps^2(\nabla q_h^N,\nabla u)+\Delta t\,\Gamma\left(F_+{''}(\phi_h^N)q_h^N,u\right)-\Delta t\,D_{\sigma}\chi(\nabla r_h^N,\nabla u)\\
\displaystyle & \qquad +\Delta t\,P_0 \delta \chi\left(P(\phi_{h}^{N-1})(p_h^N-r_h^N),u\right)=\lambda_1\left(\phi_h^N-\phimeas,u\right), \\
\notag \\
\label{eqn:qh2} \displaystyle & (q_{h}^N,v)=D_{\phi}(\nabla p_{h}^N,\nabla v)+P_0 \delta\left(P(\phi_{h}^{N-1})(p_h^N-r_h^N),v\right),\\ \notag \\
\displaystyle & \notag \left(r_h^N,w\right)+\Delta t\frac{D_{\sigma}}{\delta}(\nabla r_h^N,\nabla w)-\Delta t\,\chi\left(q_h^N,w\right)\\
\label{eqn:qh3} & \displaystyle \qquad -\Delta t P_0 \left(P(\phi_{h}^{n-1})\left(p_h^N-r_h^N\right),w\right) + \Delta t\, \kappa \left(r_h^N,w\right)=\lambda_2\left(\sigma_h^N-\sigmameas,w\right).
\end{alignat}
\Accorpa\eqnqhin {eqn:qh1} {eqn:qh3}
\\ \\
\noindent
\textbf{Problem $\mathbf{Q}_2^{h}\left(p_h^N,q_h^N,r_h^N,\left(\phi_h^n\right)_{n=1}^N,\left(\mu_h^n\right)_{n=1}^N,\left(\sigma_h^n\right)_{n=1}^N\right)$: Adjoint problem.}\\
\noindent
The derivatives $(\D_{\phi_h^n}\mathcal{L}_h(\vec{x}_h),u)=0$, $(\D_{\mu_h^n}\mathcal{L}_h(\vec{x}_h),v)=0$ and $(\D_{\sigma_h^n}\mathcal{L}_h(\vec{x}_h),w)=0$, for $n=1,\dots,N-1$, along directions $(u,v,w)\in S^h\times S^h\times S^h$, give the adjoint problem: for $n=1,\dots,N-1$, given $(p_h^N,q_h^N,r_h^N)\in S^h\times S^h\times S^h$, find $(p_h^n,q_h^n,r_h^n)\in S^h\times S^h\times S^h$ such that for all $(u,v,w)\in S^h\times S^h\times S^h$,
\begin{alignat}{2}
\label{eqn:qha1}
\displaystyle & \notag \left(\frac{p_h^n-p_h^{n+1}}{\Delta t},u\right)+\eps^2(\nabla q_h^n,\nabla u)+\Gamma\left(F_+{''}(\phi_h^n)q_h^n,u\right)+\Gamma\left(F_-{''}(\phi_h^n)q_h^{n+1},u\right)-D_{\sigma}\chi(\nabla r_h^n,\nabla u)\\
\displaystyle & \notag \qquad -P_0 \delta \left(P'(\phi_{h}^{n})\left(\frac{1}{\delta}\sigma_{h}^{n+1} + \chi (1-\phi_{h}^{n+1}) - \mu_{h}^{n+1}\right)\left(p_h^{n+1}-r_h^{n+1}\right),u\right)\\
\displaystyle & \qquad +P_0 \delta \chi\left(P(\phi_{h}^{n-1})(p_h^n-r_h^n),u\right)+\left(c\hh' (\phi_{h}^{n})p_h^{n+1},u\right)=0,\\
\notag \\
\label{eqn:qha2} \displaystyle & (q_{h}^n,v)=D_{\phi}(\nabla p_{h}^n,\nabla v)+P_0 \delta\left(P(\phi_{h}^{n-1})(p_h^n-r_h^n),v\right),\\ \notag \\
\displaystyle & \notag \left(\frac{r_h^n-r_h^{n+1}}{\Delta t},w\right)+\frac{D_{\sigma}}{\delta}(\nabla r_h^N,\nabla w)-\chi\left(q_h^N,w\right)\\
\label{eqn:qha3} & \displaystyle \qquad -P_0 \left(P(\phi_{h}^{n-1})\left(p_h^n-r_h^n\right),w\right) + \kappa \left(r_h^n,w\right)=0.
\end{alignat}
\Accorpa\eqnqha {eqn:qha1} {eqn:qha3}\\ \\
\noindent
\textbf{Problem $\mathbf{Q}_3^{h}\left(p_h^1,q_h^1,r_h^1\right)$: Optimality condition.}\\
\noindent
The variational inequality $(\D_{\phi_h^0}\mathcal{L}_h(\vec{x}_h),\tilde{\phi}-\bar{\phi}_h^0)\geq 0$, for any admissible direction $\tilde{\phi}\in \Uad^h$, where $\bar{\phi}_h^0$ is an optimal solution to (MP)$_h$, gives the optimality condition: find $\bar{\phi}_h^0\in \Uad^h$ such that, for all $\tilde{\phi}\in \Uad^h$,
\begin{align}
   \label{eqn:qhop} 
  \displaystyle & \notag \left(p_h^1+\Delta t\,\left[P_0 \delta P'(\bar{\phi}_{h}^{0})\left(\frac{1}{\delta}\sigma_{h}^{1} + \chi (1-\phi_{h}^{1}) - \mu_{h}^{1}\right)\left(p_h^{1}-r_h^{1}\right)-c\hh' (\bar{\phi}_{h}^{0})p_h^{1}-\Gamma F_-{''}(\bar{\phi}_h^0)q_h^{1}\right],\tilde{\phi}-\bar{\phi}_h^0\right)\\
  \displaystyle & \qquad + \alpha_1\left(\Gamma F'(\bar{\phi}_h^0),\tilde{\phi}-\bar{\phi}_h^0\right)+ \alpha_1 \eps^2\left(\nabla \bar{\phi}_h^0,\nabla(\tilde{\phi}-\bar{\phi}_h^0)\right)\geq 0.
\end{align}
We observe that we have an extra term appearing in the discrete version of the optimality condition with respect to \eqref{var:ineq}. 
This is due to the time discretisation procedure.
We introduce the operator $\Pi_{\phi}:\Uad^h\to S^h$, which associates to an element $\bar{\phi}_h^0\in \Uad^h$ an element $\Pi_{\phi}(\bar{\phi}_h^0)\in S^h$ such that, for all $\psi \in S^h$,
\begin{align*}
   \displaystyle & (\Pi_{\phi}(\bar{\phi}_h^0),\psi)\\
  \displaystyle & \quad =\left(p_h^1+\Delta t\,\left[P_0 \delta P'(\bar{\phi}_{h}^{0})\left(\frac{1}{\delta}\sigma_{h}^{1} + \chi (1-\phi_{h}^{1}) - \mu_{h}^{1}\right)\left(p_h^{1}-r_h^{1}\right)-c\hh' (\bar{\phi}_{h}^{0})p_h^{1}-\Gamma F_-{''}(\bar{\phi}_h^0)q_h^{1}\right],\psi\right)\\
  \displaystyle & \qquad + \alpha_1 \left(\Gamma F'(\bar{\phi}_h^0),\psi\right) + \alpha_1 \eps^2\left(\nabla \bar{\phi}_h^0,\nabla \psi\right).
\end{align*}
Then, \eqref{eqn:qhop} can be rewritten as
\[
(\Pi_{\phi}(\bar{\phi}_h^0),\tilde{\phi}-\bar{\phi}_h^0)\geq 0.
\]
Moreover, the variational inequality $(\D_{\sigma_h^0}\mathcal{L}_h(\vec{x}_h),\tilde{\sigma}-\bar{\sigma}_h^0)\geq 0$, for any admissible direction $\tilde{\sigma}\in \Uad^h$, where $\bar{\sigma}_h^0$ is an optimal solution to (MP)$_h$, gives the optimality condition: find $\bar{\sigma}_h^0\in \Uad^h$ such that, for all $\tilde{\sigma}\in \Uad^h$,
\begin{align}
   \label{eqn:qhopsigma} 
  \displaystyle & \left(r_h^1+\alpha_2\bar{\sigma}_h^0,\tilde{\sigma}-\bar{\sigma}_h^0\right)\geq 0.
\end{align}
Similarly as before, we introduce the operator $\Pi_{\sigma}:\Uad^h\to S^h$, which associates to an element $\bar{\sigma}_h^0\in \Uad^h$ an element $\Pi_{\sigma}(\bar{\sigma}_h^0)\in S^h$ such that, for all $\psi \in S^h$,
\begin{align*}
   \displaystyle & (\Pi_{\sigma}(\bar{\sigma}_h^0),\psi) =\left(r_h^1+\alpha_2\bar{\sigma}_h^0,\psi\right).
\end{align*}
Then, \eqref{eqn:qhopsigma} can be rewritten as
\[
(\Pi_{\sigma}(\bar{\sigma}_h^0),\tilde{\sigma}-\bar{\sigma}_h^0)\geq 0.
\]
We solve the latter variational inequalities through an iterative algorithm based on the projected gradient method, with a learning rate chosen according to a line search along the
descent directions to guarantee global convergence (see e.g. \cite{BR2015}).
\\
\noindent
In order to report the employed optimisation algorithm, let us introduce the projection function $\mathcal{P}_{k,i}^{\rho}:\Uad^h\to \Uad^h$, $i\in \{\phi,\sigma\}$, which associates to an element $\phi_{h}\in \Uad^h$ an element $\mathcal{P}_{k,i}^{\rho}(\phi_h)\in \Uad^h$ such that
\[
\mathcal{P}_{k,i}^{\rho}(\phi_h):=\pi^h_{[0,1]}\left(\phi_h-\rho \Pi_i(\phi_h)\right),
\]
depending on a parameter $\rho\in \mathbb{R}^+$, where 
\[\pi^h_{[0,1]}(\phi_h):=\sum_{j=1}^{N_h}\max\{0,\min\{1,\phi_h(\mathbf{x}_j)\}\}\chi_j.\]
Here, $\chi_j\in S^h$, for $j=1,\dots,N_h$ are the Lagrangian basis functions of $S^h$ associated to the $j-$th node of coordinates $\mathbf{x}_j$, with $N_h$ the number of nodes of $\mathcal{T}_h$. We observe that $\chi_j(\cdot)\in [0,1]$. \\
\noindent
At each iterative step of the optimisation algorithm, we refine the mesh based on a D\"{o}rfler marking scheme, with an error indicator given cell-wise as the sum of the jumps of the normal derivatives of the updated value of the initial condition $\phi_h^0$ at the previous step. See e.g. \cite{GHK2016} for the use of this kind of indicator for mesh adaptation for the Cahn--Hilliard equation. Then, we interpolate the current initial condition on the new mesh. We thus define the following: \\ \\
\noindent
\textbf{Problem $\mathcal{M}^h(\mathcal{T}_h^k,\phi_h)$: mesh refinement.}\\
Given a mesh $\mathcal{T}_h^k$ and $\phi_h\in S^h$, define a new mesh $\mathcal{T}_h^{k+1}$ obtained by refinement according to the previously described indicator associated to $\phi_h$ and interpolate $\phi_h$ on $\mathcal{T}_h^{k+1}$.\\ \\
\noindent
We finally formulate the following \textbf{Optimisation algorithm}:
\begin{algorithm}[h!]
\caption{Optimisation algorithm}\label{alg:opt}
\begin{algorithmic}
\Require $0<\rho<1$, $0<\iota<1$, $\mathcal{T}_h$, $(\phi_h^0,\sigma_h^0)\in \Uad^h\times \Uad^h$, $(\phimeas,\sigmameas) \in H\times H$;
\State \textbf{Initialization} Set $\phi_h^{0,0}\leftarrow \phi_h^0$, $\sigma_h^{0,0}\leftarrow \sigma_h^0$, $\mathcal{T}_h^0\leftarrow \mathcal{T}_h$;
\For{$k\geq 1$}  
    \State \textbf{Mesh Refinement} Set $\mathcal{T}_h^k=\mathcal{M}^h(\mathcal{T}_h^{k-1},\phi_h^{0,k-1})$; 
    \State \textbf{Step 1--\textbf{Direct problem}}: \[\left(\left(\phi_h^{n,k}\right)_{n=1}^N,\left(\mu_h^{n,k}\right)_{n=1}^N,\left(\sigma_h^{n,k}\right)_{n=1}^N\right)=P^h(\phi_h^{0,k-1},\sigma_h^{0,k-1});\]
    \State \textbf{Step 2--\textbf{Initialization}}: \[\left(p_h^{N,k},q_h^{N,k},r_h^{N,k}\right)=Q_1^h(\phi_h^{N,k},\phi_h^{N-1,k},\sigma_h^{N,k});\]
    \State \textbf{Step 3--\textbf{Adjoint problem}}: \[\left(\left(p_h^{n,k}\right)_{n=1}^{N-1},\left(q_h^{n,k}\right)_{n=1}^{N-1},\left(r_h^{n,k}\right)_{n=1}^{N-1}\right)=Q_2^h\left(p_h^{N,k},q_h^{N,k},r_h^{N,k},\left(\phi_h^{n,k}\right)_{n=1}^N,\left(\mu_h^{n,k}\right)_{n=1}^N,\left(\sigma_h^{n,k}\right)_{n=1}^N\right);\]
    \State \textbf{Step 4--\textbf{Gradient projection with line search}}:
    Calculate 
    \[\tilde{\phi}_h^{0,k}=\mathcal{P}_{k,\phi}^{\rho}(\phi_h^{0,k-1}) \;\; \text{and set} \;\; v_{k,\phi}:=\tilde{\phi}_h^{0,k}-{\phi}_h^{0,k-1};\]
    \[\tilde{\sigma}_h^{0,k}=\mathcal{P}_{k,\sigma}^{\rho}(\sigma_h^{0,k-1}) \;\; \text{and set} \;\; v_{k,\sigma}:=\tilde{\sigma}_h^{0,k}-{\sigma}_h^{0,k-1};\]
    \If {$\norm{v_{k,\phi}}_H+\norm{v_{k,\sigma}}_H< \text{tol}_v$}
    \State $\bar{\phi}_h^0\leftarrow \tilde{\phi}_h^{0,k}$, \, $\bar{\sigma}_h^0\leftarrow \tilde{\sigma}_h^{0,k}$;
    \State break;
\EndIf 
\State Determine the minimal $m_k$ such that
\begin{multline}
\mathcal{J}\left(\phi_h^{N,k},\phi_h^{0,k-1}+\rho^{m_k}v_{k,\phi},\sigma_h^{N,k},\sigma_h^{0,k-1}+\rho^{m_k}v_{k,\sigma} \right)\leq \mathcal{J}\left(\phi_h^{N,k},\phi_h^{0,k-1},\sigma_h^{N,k},\sigma_h^{0,k-1}\right)\\
+\rho^{m_k}\iota (\Pi_{\phi}(\phi_h^{0,k-1}),v_{k,\phi})+\rho^{m_k}\iota (\Pi_{\sigma}(\sigma_h^{0,k-1}),v_{k,\sigma});
\end{multline}
\State Update $\phi_h^{0,k}=\phi_h^{0,k-1}+\rho^{m_k}v_{k,\phi},\, \sigma_h^{0,k}=\sigma_h^{0,k-1}+\rho^{m_k}v_{k,\sigma}$;
\EndFor
\State Set $\bar{\phi}_h^{0}=\phi_h^{0,k}$, $\bar{\sigma}_h^{0}=\sigma_h^{0,k}$.
\end{algorithmic}
\end{algorithm}
\newpage

\section{Simulation study and results}
\label{sec:simulations}

In this section, we report numerical results regarding the application of Algorithm \ref{alg:opt} to representative test cases in two and three space dimensions. We consider a first Test Case $1$ in two spatial dimensions where we reconstruct both the initial conditions $\phi_0$ and $\sigma_0$, given target distributions for $\phi(T)$ and $\sigma(T)$, showing qualitatively the good performance of the optimisation algorithm \ref{alg:opt}. In subsequent Test Cases $2$-$6$ we consider the problem of approximating only the initial condition $\phi_0$ for the tumour phase-field, given a target distribution for $\phi(T)$, in both two and three spatial dimensions. As detailed in the Introduction, this is the relevant case in real applications, in which we typically only have measurements for $\phi(T)$. For what concerns $\sigma_0$, we essentially fix it to be a stationary solution of equation \eqref{eq:sigmaad} with chemotactic and source terms depending on $\phi_0$. 
Such a procedure is still physically reasonable, since nutrient diffusion generally happens on a much faster time scale than tumour growth. For Test Cases $2$-$6$ we report the qualitative results proving the good performance of the optimisation algorithm \ref{alg:opt} in reconstructing the Ground Truth initial condition, together with further detailed results investigating its convergence properties. We moreover explore the behaviour of the algorithm by varying in different ways all the following parameters: the Tikhonov regularisation parameter, the final time horizon, the position of the initial guess and the spatial dimension ($2D$ or $3D$).

In all test cases, we consider the following values for the model parameters, which are chosen in the biologically feasible ranges reported in \cite{Agosti1}: 
\begin{center}
\begin{tabular}{ c c c }
 Parameter & Value & Units \\ 
 \hline
 $P_0$ & 0.1 & $day^{-1}$ \\  
 $\tilde{\delta}$ & 0.001 & adimensional\\
 $D_{\phi}$ & 0.00053 & $\frac{mm^3}{N\,day}$\\
 $D_{\sigma}$ & $0.001$ & $\frac{mm^3}{N\,day}$\\
 $\Gamma$ & 2.5 & $N/mm$\\
 $\tilde{\chi}$ & 0.5 & adimensional\\
 $\eps$ & 0.025 & $\left(\frac{N}{mm}\right)^{1/2}mm$\\
 $c$ & 0.02 & $day^{-1}$ \\ 
 $\kappa$ & 0.12 & $day^{-1}$ \\ 
 \hline
\end{tabular}
\end{center}
The parameter values have been reported with units corresponding to two space dimensions.
In all test cases, we consider $\rho=0.9$ as the starting guess for the learning parameter. 
In Test Case $1$ we take $\lambda_1=\lambda_2=1$ (adimensional) and $\alpha_1=\alpha_2=0.01$ in \eqref{eqn:jdisc}. We consider $\alpha_1=\tilde{\alpha}_1/\Gamma$, with $\tilde{\alpha}_1$ adimensional. In this way all the contributions in \eqref{eqn:jdisc} have the same units of $mm^2$, hence are comparable.
In Test Cases $2-6$ we consider situations where we only have measurements for $\phi(T)$, hence we take $\lambda_2=\alpha_2=0$ in \eqref{eqn:jdisc}. In the latter test cases we put $\lambda_1=1$ (adimensional) and $\alpha_1:=\alpha=\tilde{\alpha}/\Gamma$, with $\tilde{\alpha}$ adimensional, varying between different test cases. Also, we consider $\Delta t=0.01$ days. The value of the final time $T$ varies throughout the test cases. We observe that a patient diagnosed with a tumour undergoes control routines, consisting of MRI acquisitions, with a frequency which typically varies from $1$ month to $3$ months after diagnosis (see e.g. the clinical test cases concerning brain tumour considered in \cite{Agosti1}). Hence, thinking of future applications, we will vary the values of $T$ up to the value of $3$ months, to test the efficiency and stability of the numerical algorithm in these time ranges. Anyhow, we observe that we could afford the resolution of the optimisation problems also with higher values of $T$.

\subsection{Test case $1$}
\begin{figure}[t!]
\includegraphics[width=0.95\linewidth]
{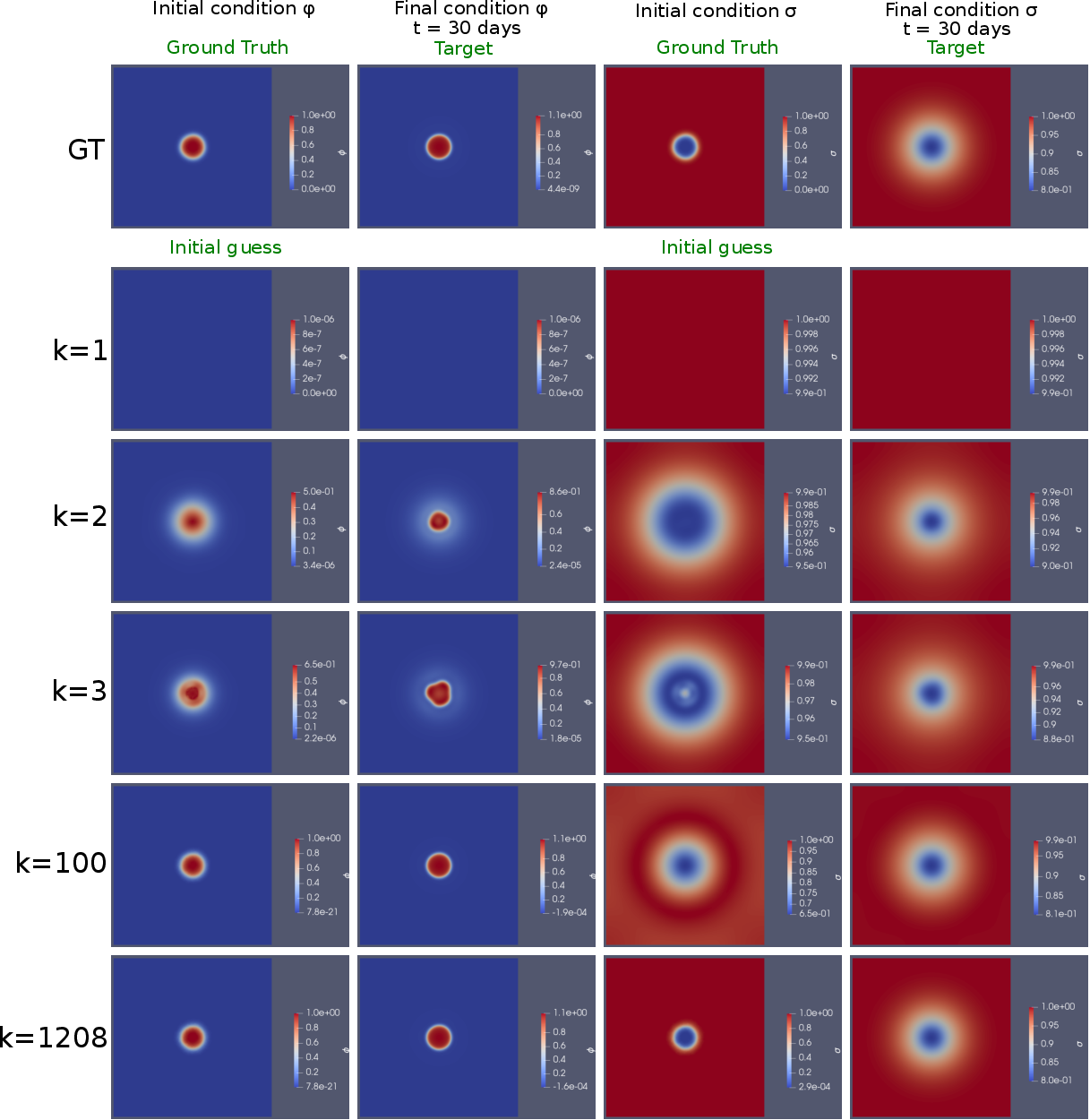}
\centering
\caption{Plots of the initial condition and the final condition at $t=30$ days, for both the variables $\phi$ (I-II columns) and $\sigma$ (III-IV columns), for the ground truth (GT) and for different values of the iteration step $k$, in the case of the initial guess $(\phi_0=0,\sigma_0=1)$.}
\label{fig:01}
\end{figure}

\begin{figure}[t!]
\includegraphics[width=0.95\linewidth]
{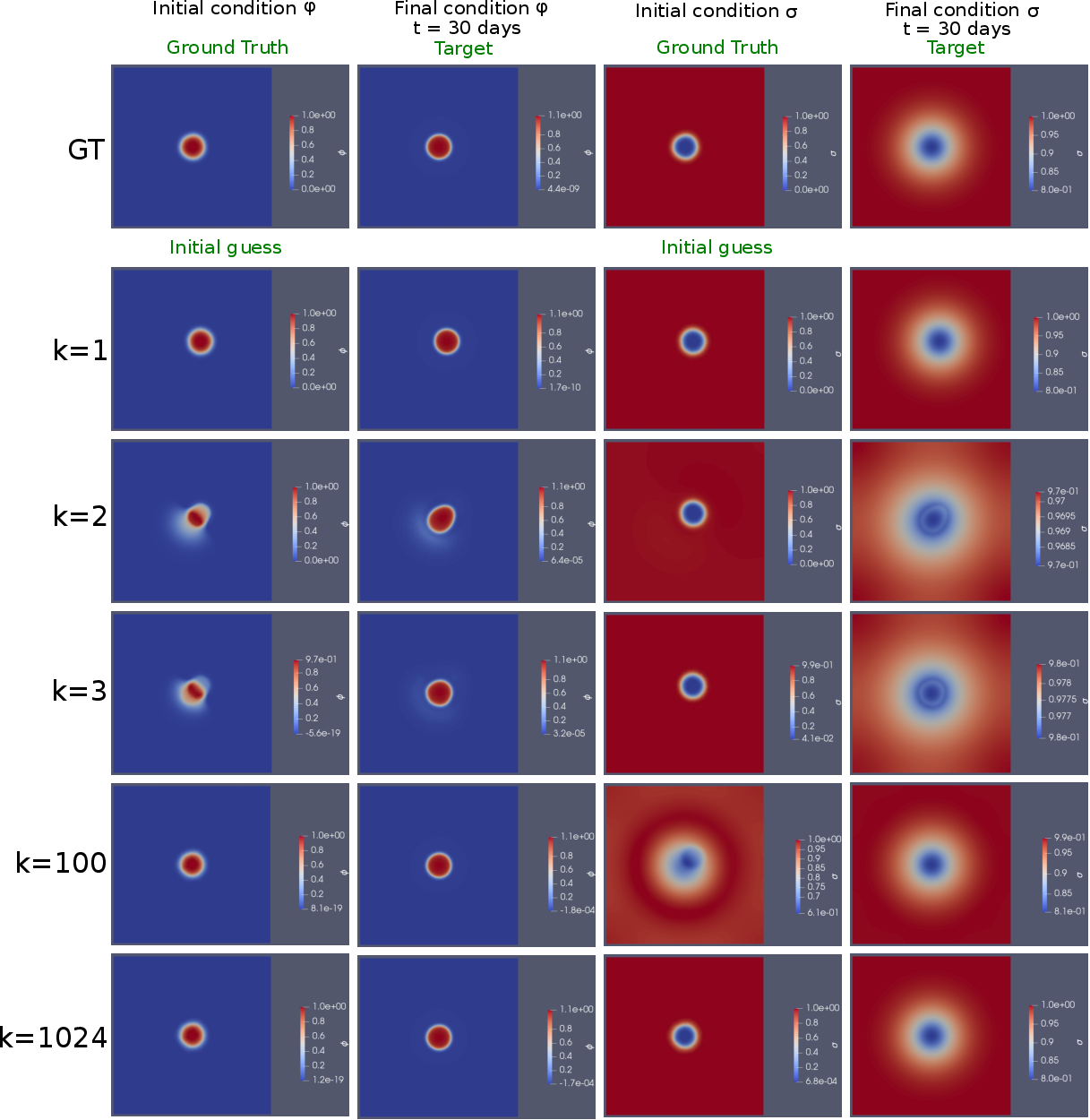}
\centering
\caption{Plots of the initial condition and the final condition at $t=30$ days, for both the variables $\phi$ (I-II columns) and $\sigma$ (III-IV columns), for the ground truth (GT) and for different values of the iteration step $k$, in the case in which the initial guess is the characteristic function of a circle shifted from the ground truth configuration.}
\label{fig:02}
\end{figure}

In Test Case $1$ we firstly show the qualitative results of two numerical tests in two spatial dimensions in which we have measurements of both $\phi(T)$ and $\sigma(T)$, varying the initial guesses.
We consider a domain $\Omega=[-5 mm, 5 mm]^2$, and a ground truth initial condition $\phi_0$ given by a smoothed characteristic function of a circular domain of radius $0.6$ centred in the origin, and $\sigma_0=1-\phi_0$. The targets $\phi_{\text{meas}}$ and $\sigma_{\text{meas}}$ are obtained as the final conditions, corresponding to the evolution of the ground truth, after $T=30$ days. The numerical simulations employed to generate the ground truths and the targets are obtained with a base mesh of size $80\times 80$, with further refinement as detailed in the previous sections. 
\noindent
We compare the results obtained by initiating the Algorithm \ref{alg:opt} firstly with a null initial guess $\phi_0=0$ for the tumour concentration and a saturated initial guess $\sigma_0=1$ for the nutrient concentration, (Figure \ref{fig:01}), and secondly with an initial guess where $\phi_0$ is given by the characteristic function of a circle shifted from the centred position of the Ground Truth configuration, and $\sigma_0=1-\phi_0$, (Figure \ref{fig:02}), with a base mesh of size $64\times 64$ in both situations. In figures \ref{fig:01} and \ref{fig:02} we report the plots of the initial conditions and the final conditions at $t=30$ days, for both the variables $\phi$ and $\sigma$, for the ground truth (GT) and for different values of the iteration step $k$, corresponding to the initial guess $(\phi_0=0,\sigma_0=1)$ and to the shifted circle configuration respectively.


We observe that the ground truth configurations for both $\phi_0$ and $\sigma_0$ are correctly reconstructed in both situations. Due to the difference in the time scales between the dynamics of $\phi$ and $\sigma$, we particularly observe that the ground truth reconstruction for $\sigma_0$ is much slower than that for $\phi_0$. Indeed, since the time scale of the nutrient dynamics is orders of magnitude smaller than that of the tumour cells dynamics, the final condition $\sigma(T)$ at each iterative step $k$ is close to a stationary solution of \eqref{eq:sigmaad} with the source term depending on $\phi(T)$; hence, as far as $\phi_h^{0,k}$ approximates the Ground Truth, the iterative increments in the gradient projected method to update $\sigma_h^{0,k}$ become smaller even if $\sigma_h^{0,k}$ is far from the Ground Truth, making the convergence for the $\sigma$ variable much slower than for the variable $\phi$.

\subsection{Test case $2$}
\begin{figure}[t!]
\includegraphics[width=1.0\linewidth]
{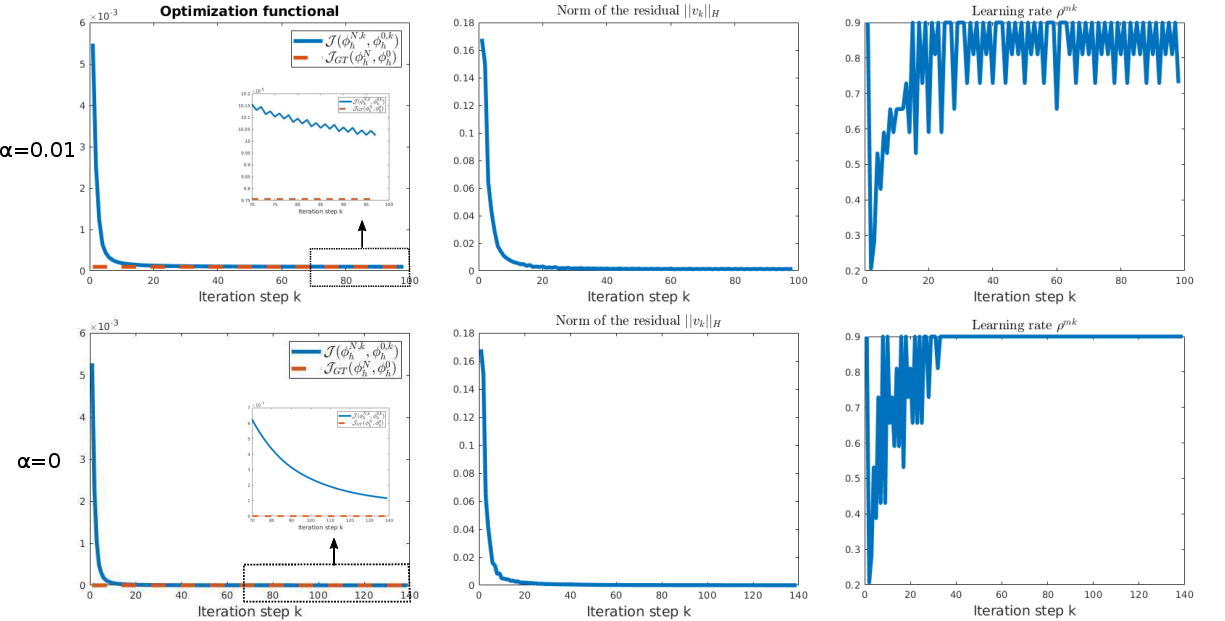}
\centering
\caption{Plots of the optimisation functional $\mathcal{J}$, of the norm $\norm{v_k}_H$ and of the learning rate $\rho^{m_k}$ vs the iteration steps $k$ for the cases ${\alpha}=0.01$ and ${\alpha}=0$.}
\label{fig:0}
\end{figure}

\begin{figure}[h!]
\includegraphics[width=0.95\linewidth]
{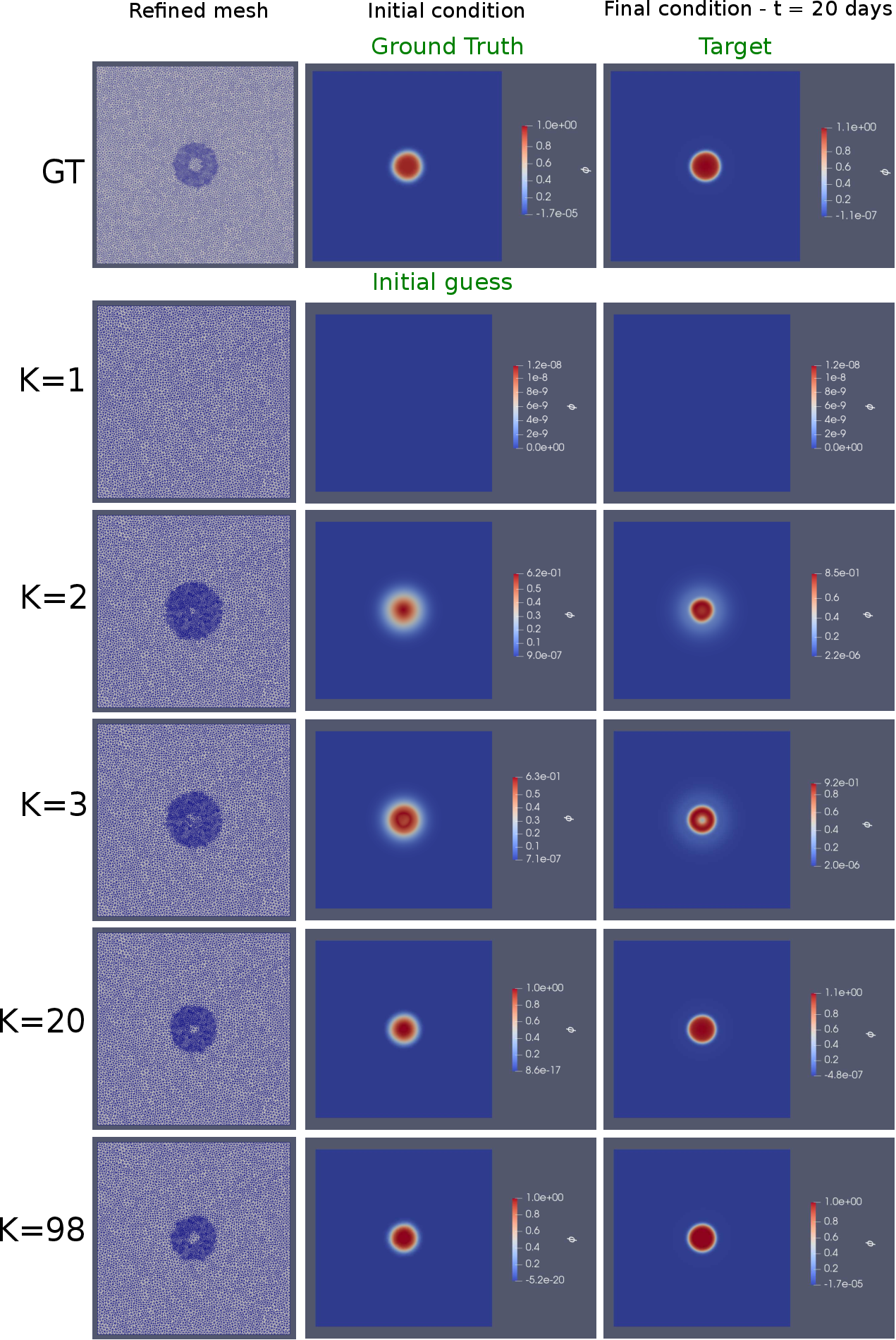}
\centering
\caption{Plots of the refined mesh, the initial condition and the final condition at $t=20$ days for the ground truth (GT) and for different values of the iteration step $k$, in the case ${\alpha}=0.01$.}
\label{fig:1}
\end{figure}

\begin{figure}[hp!]
\includegraphics[width=0.9\linewidth]
{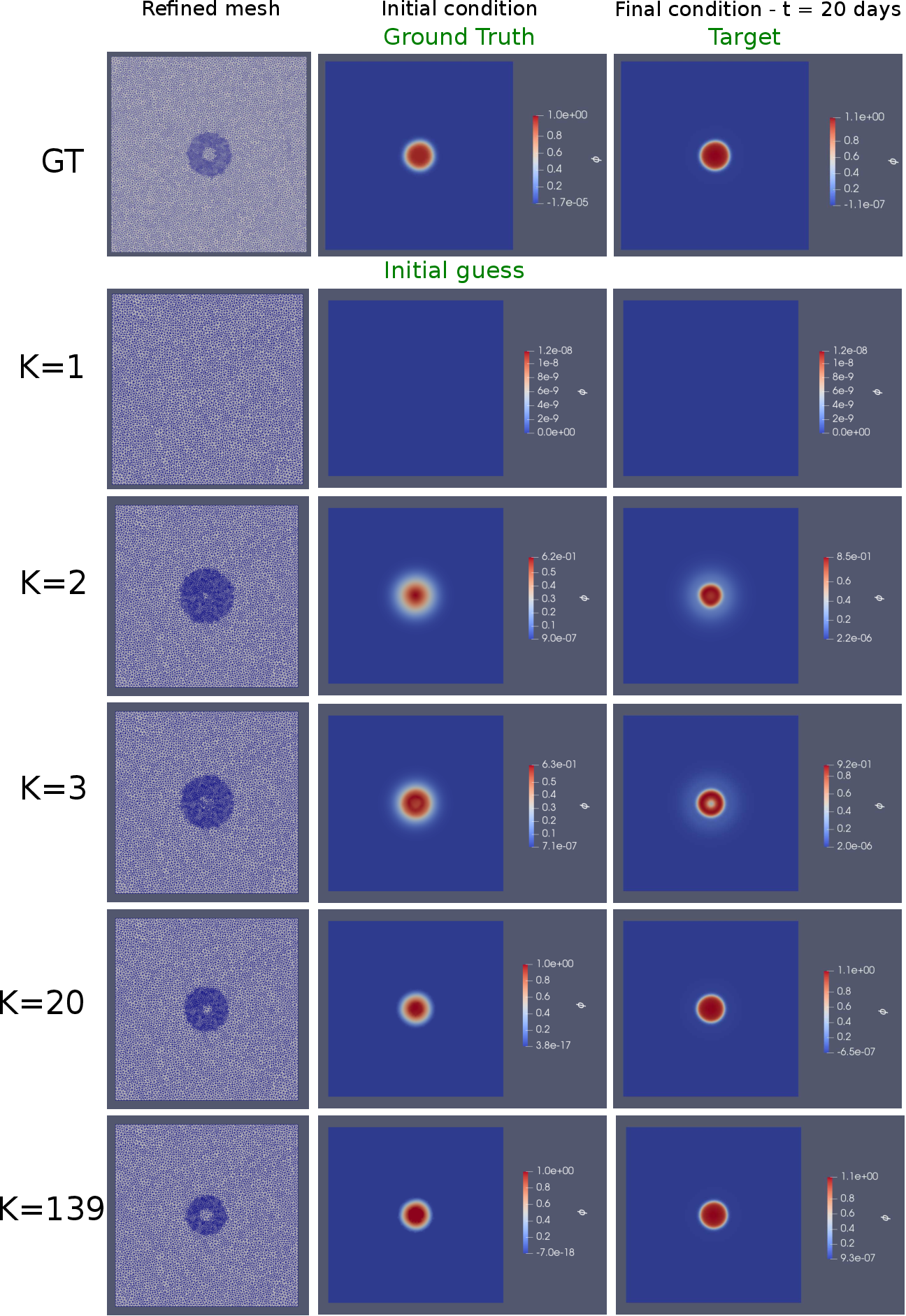}
\centering
\caption{Plots of the refined mesh, the initial condition and the final condition at $t=20$ days for the ground truth (GT) and for different values of the iteration step $k$, in the case ${\alpha}=0$.}
\label{fig:2}
\end{figure}

We now specialise the test cases to situations where we have only a measurement of $\phi(T)$ (which is relevant for applications), hence setting $\lambda_2=\alpha_2=0$ in \eqref{oc:problem}.
As a first test case in this direction, we consider, as in Test Case $1$, a domain $\Omega=[-5 mm, 5 mm]^2$, and a ground truth initial condition given by a smoothed characteristic function of a circular domain of radius $0.6$ centred in the origin. The target $\phi_{\text{meas}}$ is obtained as the final condition, corresponding to the evolution of the ground truth, after $T=20$ days. The numerical simulations employed to generate the ground truth and the target are obtained with a base mesh of size $80\times 80$, with further refinement as detailed in the previous section. We then initiate the Algorithm \ref{alg:opt} with a null initial guess, and a base mesh of size $64\times 64$. We compare the results obtained by choosing ${\alpha}=0.01$ and ${\alpha}=0$, in order to assess the effects of the Tikhonov regularisation on the numerical results. In figure \ref{fig:0} we plot the optimisation functional $\mathcal{J}$, the norm $\norm{v_k}_H$ and the learning rate $\rho^{m_k}$ versus the iteration steps $k$ for the cases ${\alpha}=0.01$ and ${\alpha}=0$.
We observe that the numerical algorithm converges to a local stationary state, represented by the condition that $\norm {v_k}_H \to 0$ as $k\to 0$, with a focal convergence in the case $\alpha=0.01$ and a nodal convergence in the case $\alpha=0$. In figures \ref{fig:1} and \ref{fig:2} we report the plots of the refined mesh, the initial condition and the final condition at $t=20$ days for the ground truth (GT) and for different values of the iteration step $k$, for the cases ${\alpha}=0.01$ and ${\alpha}=0$ respectively.

We observe that the ground truth configuration is correctly reconstructed in both situations, with only slight differences at late iteration steps between the two. Hence, we conclude that Algorithm \ref{alg:opt} works also in the absence of the Tikhonov regularisation.
This is essentially due to the fact that, by Remark \ref{rmk:lipstab} we expect Lipschitz stability for the inverse problem in finite-dimensional spaces, thus Landweber-type algorithms are locally convergent. Notice, indeed, that in the absence of a regularisation (i.e. $\alpha = 0$) our method is practically equivalent to a Landweber scheme. We further explore the effect of the Tikhonov regularisation in the next experiments.

\subsection{Test case $3$}
\begin{figure}[hp!]
\includegraphics[width=1.0\linewidth]
{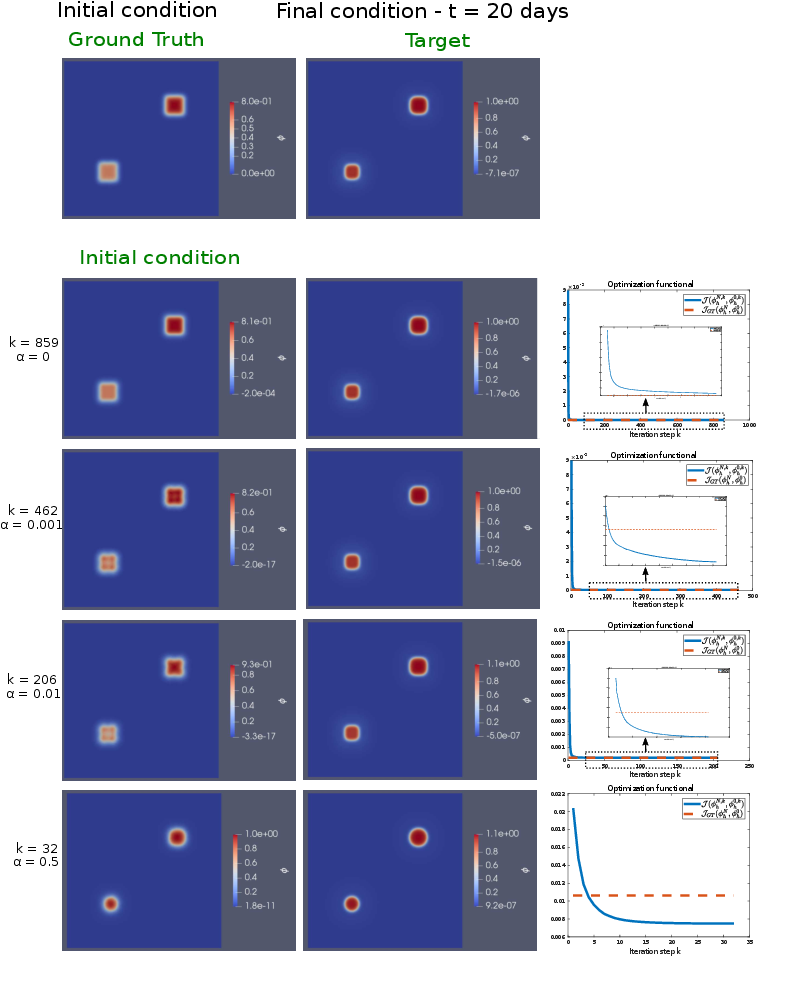}
\centering
\caption{I-II columns: Plots of the ground truth and the target configurations, compared with the reconstructed initial condition and the final condition at $T=20$ days, for values of the iteration step $k$ attained at convergence, in the cases ${\alpha}=0,0.001,0.01,0.5$, given initial conditions with squared topology and out of equilibria homogeneous values inside the square domains and a null initial guess. III column: Plot of the optimisation functional $\mathcal{J}$ versus the iteration steps $k$ for the cases ${\alpha}=0$, ${\alpha}=0.001$, ${\alpha}=0.01$ and ${\alpha}=0.5$.}
\label{fig:squares}
\end{figure}

In order to better highlight the effects of Tikhonov regularisation for increasing values of the regularisation parameter ${\alpha}$, as a third test case we consider a ground truth initial condition given by two smoothed characteristic functions of square domains of side $0.25$,  displaced along one diagonal of the domain $\Omega=[-5 mm, 5 mm]^2$ away from the origin, multiplied by $0.6$ and $0.8$ respectively. Hence, inside the two square domains, the phase field variable does not attain a minimum of the potential $F$ at the initial time. Also, the square topology is not a minimiser for the Cahn--Hilliard functional, employed as a regularisation term in \eqref{eqn:jdisc}. The target $\phi_{\text{meas}}$ is obtained as the final condition, corresponding to the evolution of the ground truth, after $T=20$ days. All the discretisation parameters are chosen as in Test Case $1$. As in the Test Case $2$, we initiate the Algorithm \ref{alg:opt} with a null initial guess. Here we compare the results obtained by choosing ${\alpha}=0$, ${\alpha}=0.001$, ${\alpha}=0.01$ and ${\alpha}=0.5$. In figure \ref{fig:squares} we report the plots of the initial condition and the final condition at $t=20$ days for the ground truth (GT) and for the value of the iteration step $k$ attained at convergence, for the cases ${\alpha}=0$, ${\alpha}=0.001$, ${\alpha}=0.01$ and ${\alpha}=0.5$, together with the plot of the optimisation functional $\mathcal{J}$ versus the iteration steps $k$ for the aforementioned cases.
Observing Figure \ref{fig:squares} we highlight the following facts:
\begin{itemize}
\item The results from the algorithm with Tikhonov regularisation ${\alpha}>0$ converge, as ${\alpha}\to0$, to the results obtained with a Landweber-like scheme (i.e. taking ${\alpha}=0$ in our algorithm).
\item Increasing the values of ${\alpha}$, the second term, i.e. the regularisation component, in the optimisation functional \eqref{eqn:jdisc} dominates with respect to the first term, i.e. the data-driven component. In the case ${\alpha}=0.5$, we observe that the optimisation functional $\mathcal{J}$ decreases further from the value attained at the ground truth, reaching local minimisers of the Cahn--Hilliard functional, i.e. circular domains with $\phi=1$ inside and $\phi=0$ outside of the domains. This latter configuration does not coincide with the ground truth, hence the $L^2$ distance between the final condition and the target is sensibly different from zero; anyhow, its Cahn--Hilliard energy is much lower than the one associated with square domains. Also in the cases $\alpha=0.001$ and $\alpha=0.01$ we can observe that the optimisation functional $\mathcal{J}$ decreases slightly further from the value attained at the ground truth.
\item The optimisation algorithm with higher values of ${\alpha}$ converges faster than with lower values of ${\alpha}$.
\end{itemize}

\subsection{Test case $4$}
\begin{figure}[t!]
\includegraphics[width=1.0\linewidth]
{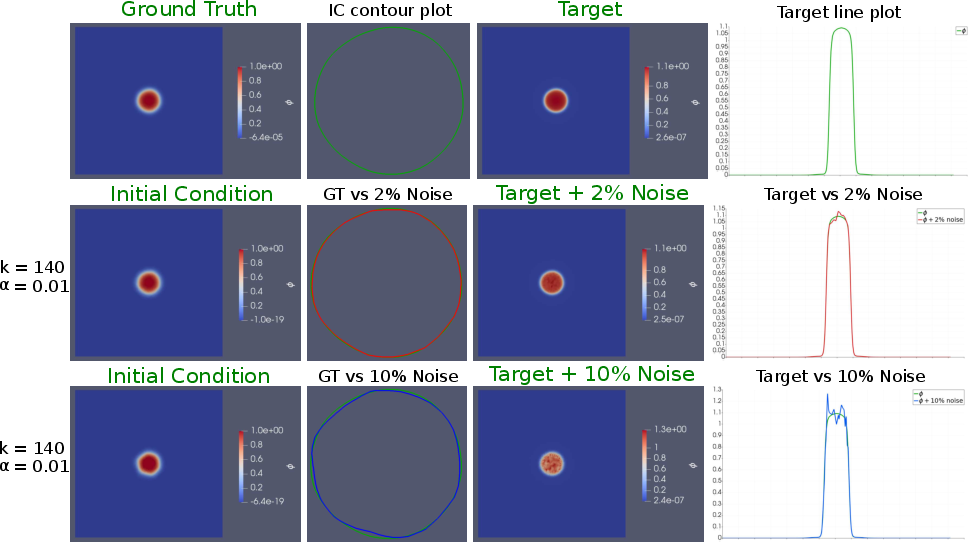}
\centering
\caption{I-II columns: plots of the ground truth and of the reconstructed initial condition at the iteration step $k=140$, together with a comparison of the contour plots for the level set $\phi=0.5$. III-IV columns: plots of the target, obtained at $T=20$ days, and of the noisy target with $2\%$ and $10\%$ level of noise, together with a comparison of the corresponding line plots along a diagonal of the square domain.}
\label{fig:noise}
\end{figure}

As a fourth test case, we consider the situation in which some noise is added to the target configuration, mimicking real measurements which are typically noisy, and a Tikhonov regularisation with $\alpha=0.01$ is considered in the optimisation functional \eqref{eqn:jdisc}. The domain is $\Omega=[-5 mm, 5 mm]^2$, and the same ground truth initial condition as Test Case $2$ is adopted.
The target $\phi_{\text{meas}}$ is obtained as the final condition, corresponding to the evolution of the ground truth, after $T=20$ days. We then add to the target a normal distributed noise with zero mean and variance equal to one, comparing two different noise levels: $2\%$ and $10\%$ of the target signal.
All the model and discretisation parameters are chosen as in Test Case $1$. As in the Test Case $2$, we initiate the Algorithm \ref{alg:opt} with a null initial guess. In figure \ref{fig:noise} we compare the results corresponding to the two noise levels. 
We observe from figure \ref{fig:noise} that the optimisation algorithm with a Tikhonov regularisation is stable against the presence of noise in the measured data, i.e. it is able to correctly reconstruct the ground truth configuration even in the presence of noisy data, with small deviations as the level of noise increases.

\subsection{Test case $5$}
\begin{figure}[hp!]
\includegraphics[width=0.9\linewidth]
{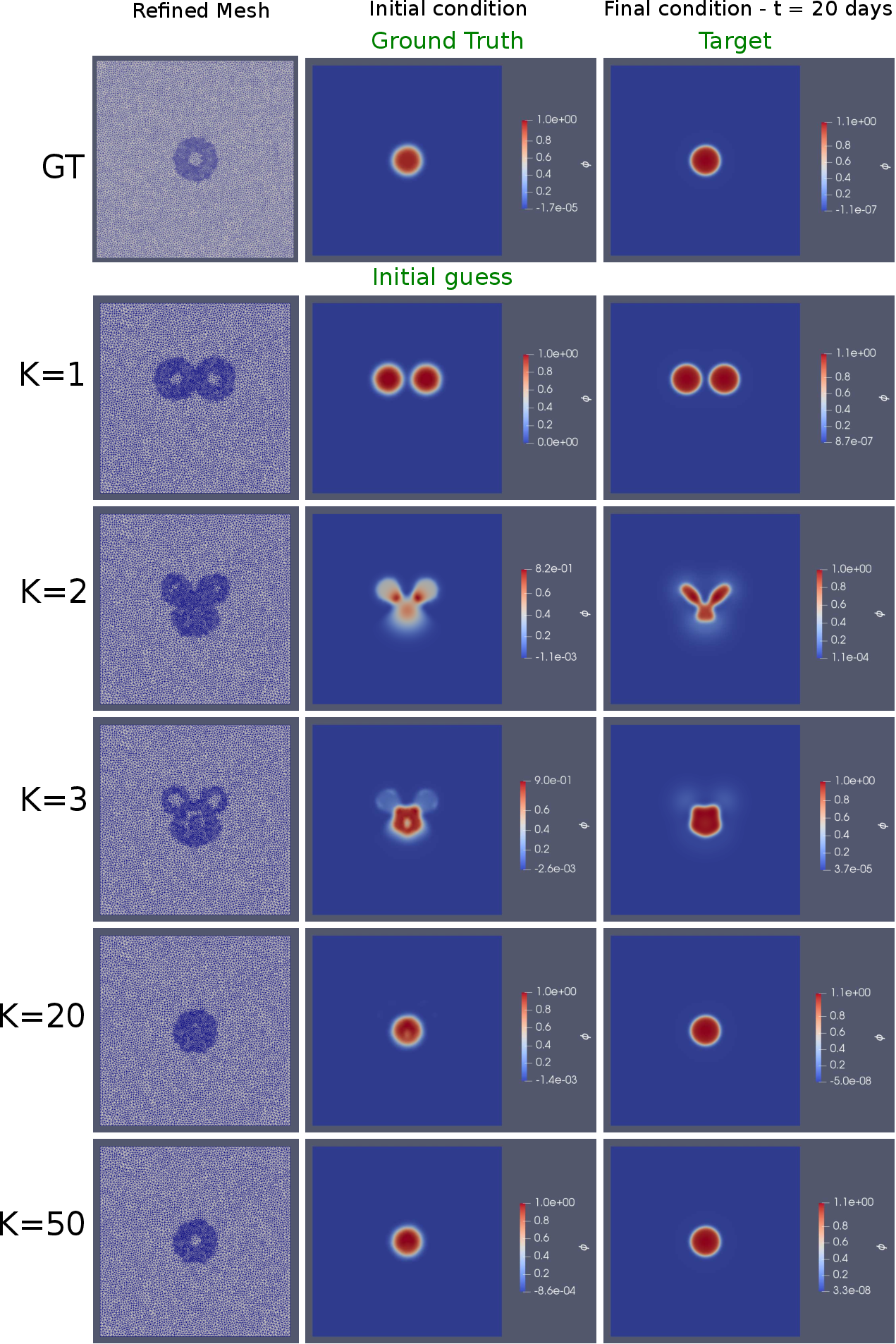}
\centering
\caption{Plots of the refined mesh, the initial condition and the final condition at $t=20$ days for the ground truth (GT) and for different values of the iteration step $k$, in the case ${\alpha}=0.01$.}
\label{fig:6}
\end{figure}

\begin{figure}[hp!]
\includegraphics[width=0.9\linewidth]
{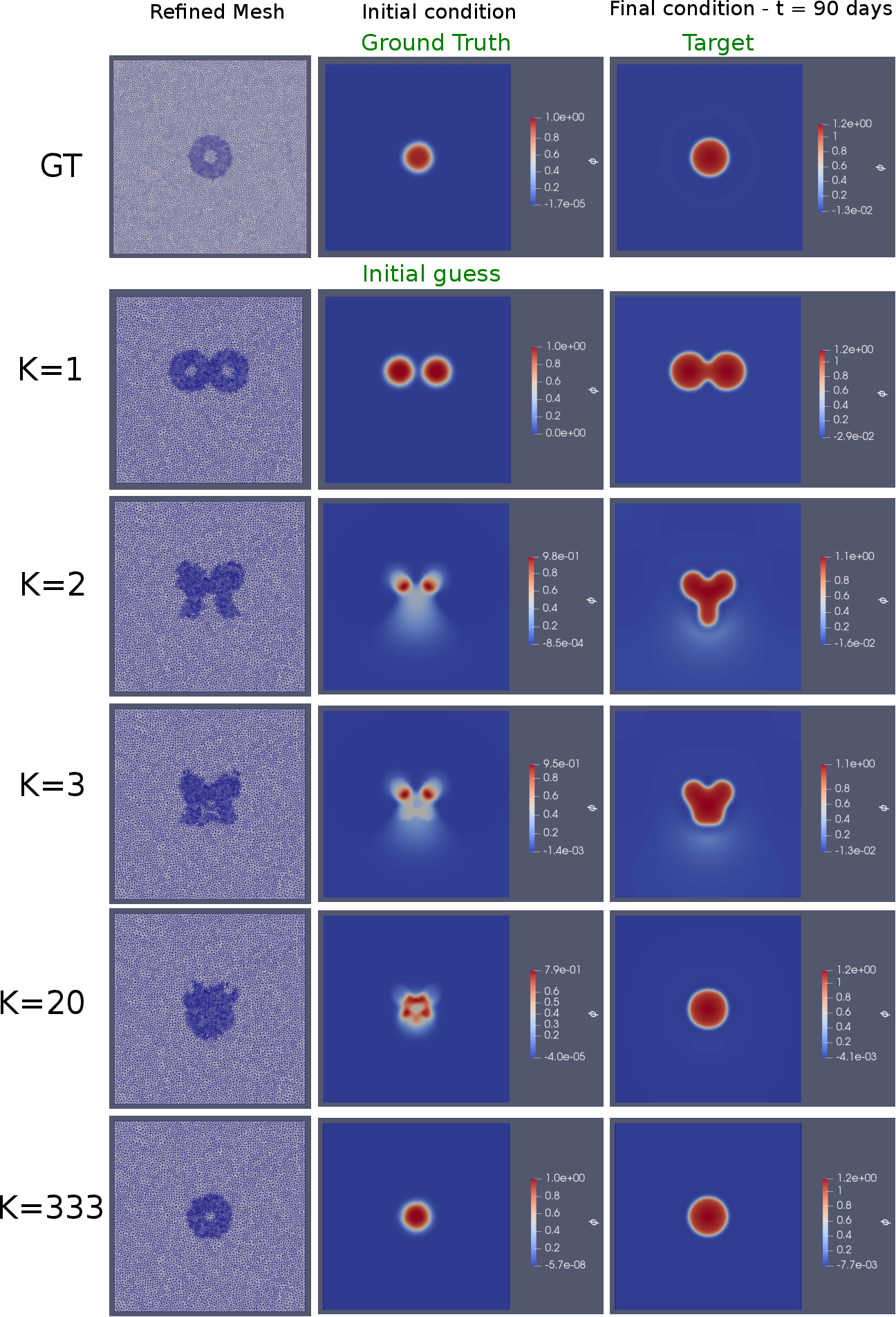}
\centering
\caption{Plots of the refined mesh, the initial condition and the final condition at $t=90$ days for the ground truth (GT) and for different values of the iteration step $k$, in the case ${\alpha}=0.01$.}
\label{fig:7}
\end{figure}

As a fifth test case, we consider a domain $\Omega=[-5 mm, 5 mm]^2$, and the same ground truth initial condition of Test Case $2$. We compare two situations where the target $\phi_{\text{meas}}$ is obtained as the final condition, corresponding to the evolution of the ground truth, after $T=20$ days and $T=90$ days, in order to assess the influence of the time interval on the numerical results. The numerical simulations employed to generate the ground truth and the target are obtained with a base mesh of size $80\times 80$, with further refinement as detailed in the previous section. We then initiate the Algorithm \ref{alg:opt} with an initial guess which is given by two circles shifted from the centred position of the ground truth configuration, in order to visualise the effects of a different topology for the initial guess with respect to the ground truth, and a base mesh of size $64\times 64$. We choose ${\alpha}=0.01$. In figures \ref{fig:6} and \ref{fig:7} we report the plots of the refined mesh, the initial condition and the final condition at $t=20$ days and $t=90$ days respectively, for the ground truth (GT) and for different values of the iteration step $k$.

We observe that the ground truth configuration is correctly reconstructed in both situations, with a much higher number of iterations required for convergence in the case with a longer time evolution.
This is expected as backward problems for parabolic equations become more and more ill-posed as the final time $T$ grows larger.

\subsection{Test case $6$}
\begin{figure}[hp!]
\includegraphics[width=0.8\linewidth]
{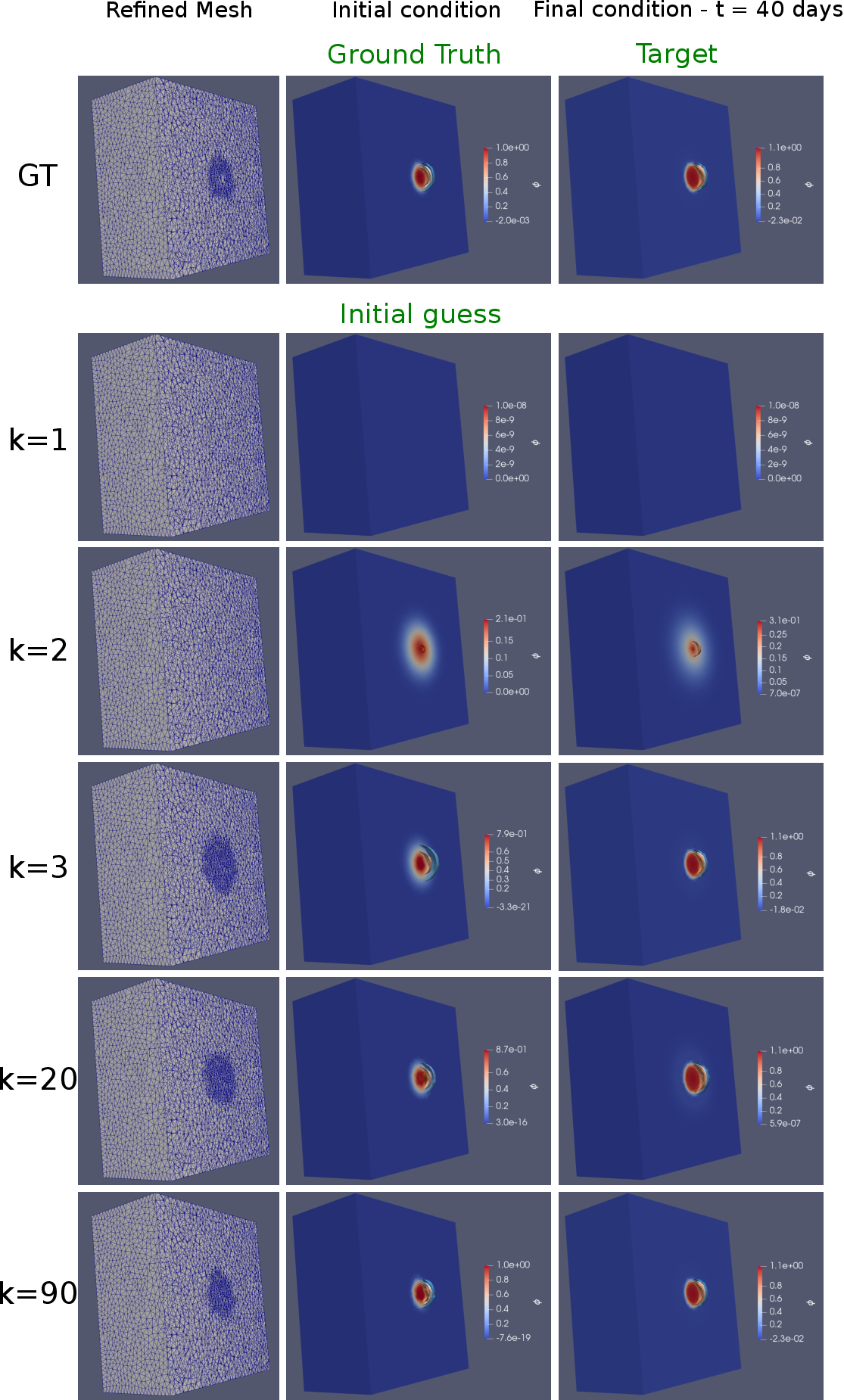}
\centering
\caption{Slice plots of the refined mesh, the initial condition and the final condition at $t=40$ days, together with the plots of level-set surfaces corresponding to the values $\phi=0.2,0.5,0.8$, for the ground truth (GT) and for different values of the iteration step $k$, in the case ${\alpha}=0.01$. }
\label{fig:8}
\end{figure}

We consider a final test case in three space dimensions, with domain $\Omega=[-5 mm, 5 mm]^3$, and a ground truth initial condition given by a smoothed characteristic function of a circular domain of radius $0.6$ centred in the origin. The target $\phi_{\text{meas}}$ is obtained as the final condition, corresponding to the evolution of the ground truth, after $T=40$ days. The numerical simulations employed to generate the ground truth and the target are obtained with a base mesh of size $80\times 80\times 80$, with further refinement as detailed in the previous section. We then initiate the Algorithm \ref{alg:opt} with a null initial guess, and a base mesh of size $60\times 60\times 60$. We choose ${\alpha}=0.01$. In Figure \ref{fig:8} we report the slice plots of the refined mesh, the initial condition and the final condition at $t=40$ days for the ground truth (GT), together with the plots of level-set surfaces corresponding to the values $\phi=0.2,0.5,0.8$, for different values of the iteration step $k$.
Also in the three-dimensional case, we observe that the ground truth configuration is correctly reconstructed, in a comparable number of iterations to the two-dimensional case.

\appendix

\section{Appendix}

\subsection{Proof of Theorem \ref{thm:weaksols}}

The argument is heavily inspired by those of \cite[Theorem 1]{FGR2015_TumGrowth} and \cite[Theorem 3.1]{GY2020}.
For this reason, we proceed formally and mostly highlight the differences due to the presence of the additional source terms.
We just note that, to make the argument rigorous, one should employ a Faedo--Galerkin discretisation with discrete spaces generated by the eigenvectors of the operator $\mathcal{N}$.
At the same time, if $s \in (4,6)$ in hypothesis \ref{ass:fconv}, one should approximate the potential $F$ with a sequence of approximating potentials having lower growth, as done in \cite[Lemma 2]{FGR2015_TumGrowth}. 
Therefore, we will just assume $2 \le s \le 4$ here and refer to the proofs of \cite[Theorem 1]{FGR2015_TumGrowth} and \cite[Theorem 3.1]{GY2020} for the more general case $2 \le s \ls 6$.

For the main energy estimate, we insert $\zeta = \mu$ in \eqref{varform:phi}, $\zeta = - \partial_t \phi$ in \eqref{varform:mu}, $\zeta = \sigma + \chi (1 - \phi)$ in \eqref{varform:sigma} and sum them all up to obtain
\begin{equation}
    \label{eq:weak:enerest}
    \begin{split}
        & \ddt \left( \mezzo \norm{\nabla \phi}^2_H + \mezzo \norm{\sigma}^2_H + \int_\Omega F(\phi) \, \de x + (\sigma, \chi(1 - \phi))_H \right) \\
        & \qquad + \norm{\nabla \mu}^2_H + \norm{\nabla (\sigma + \chi (1-\phi))}^2_H + \norm{\sqrt{P(\phi)}(\sigma + \chi (1-\phi) - \mu)}^2_H \\
        & \quad = - (\hh(\phi) c, \mu)_H + (1 - \sigma, \sigma)_H.
    \end{split}
\end{equation}
Now, for any $t \in (0,T)$, we estimate the energy term 
\[ E(t) = \int_\Omega \mezzo \abs{\nabla \phi}^2 + F(\phi) + \mezzo \abs{\sigma}^2 + \chi \sigma (1 - \phi) \, \de x \]
from below by using hypothesis \ref{ass:fbelow} and Cauchy--Schwarz and Young's inequalities. Indeed, we have that
\begin{align*}
    E(t) & \ge \mezzo \norm{\nabla \phi}^2_H + \mezzo \int_\Omega F(\phi) \, \de x + \mezzo \int_\Omega F(\phi) \, \de x + \mezzo \norm{\sigma}^2_H - \chi \norm{\sigma}_H \norm{1 - \phi}_H \\
    & \ge \mezzo \norm{\nabla \phi}^2_H + \mezzo \int_\Omega F(\phi) \, \de x + \frac{c_1}{2} \norm{\phi}^2_H - \frac{c_2}{2} \abs{\Omega} + \mezzo \norm{\sigma}^2_H - \alpha \norm{\sigma}^2_H - \frac{\chi^2}{4 \alpha} \norm{\phi}^2_H - C \\
    & \ge \mezzo \norm{\nabla \phi}^2_H + \mezzo \int_\Omega F(\phi) \, \de x + \left( \frac{c_1}{2} - \frac{\chi^2}{4 \alpha} \right) \norm{\phi}^2_H + \left( \mezzo - \alpha \right) \norm{\sigma}^2_H - C,
\end{align*}
where this holds for any $\alpha \gs 0$. 
However, given that $c_1 \gs \chi^2$ by \ref{ass:fbelow}, we can choose $\alpha = \mezzo - \delta$ for some sufficiently small $\delta \in (0, \mezzo)$ to get:
\[ E(t) \ge \mezzo \norm{\nabla \phi}^2_H + \mezzo \int_\Omega F(\phi) \, \de x + \underbrace{\left( \frac{c_1}{2} - \frac{\chi^2}{2  - 4\delta} \right)}_{:= \gamma \gs 0} \norm{\phi}^2_H + \delta \norm{\sigma}^2_H - C. \]
Next, we test \eqref{varform:mu} by $\zeta = 1$, which is possible within Galerkin's discretisation since the first eigenfunction of $\mathcal{N}$ is constant, and by hypothesis \ref{ass:fder} and H\"older's inequality we get an estimate on the mean value $\mu_\Omega$ of $\mu$, namely
\begin{equation}
\label{eq:weak:mean_mu}
    \abs{\Omega} \abs{\mu_\Omega} = \abs*{ \int_\Omega \mu \, \de x} \le \int_\Omega \abs{F'(\phi)} + \chi \abs{\sigma} \, \de x \le c_3 \int_\Omega F(\phi) \, \de x + C \norm{\sigma}^2_H + C \abs{\Omega}.
\end{equation}
At this point, we can estimate the terms on the right-hand side of \eqref{eq:weak:enerest} in the following way:
\begin{align*}
    & \abs{(\hh(\phi) c, \mu)_H} \le \hh_\infty c_\infty \int_\Omega \abs{\mu - \mu_\Omega} + \abs{\mu_\Omega} \, \de x \le C \int_\Omega \abs{\nabla \mu} \, \de x + \abs{\Omega} \abs{\mu_\Omega} \\
    & \qquad \le \mezzo \norm{\nabla \mu}^2_H + C \int_\Omega F(\phi) \, \de x + C \norm{\sigma}^2_H + C, \\
    & (1 - \sigma, \sigma)_H \le C \norm{\sigma}^2_H + C,
\end{align*}
where we used Poincaré--Wirtinger's inequality for the first term.
Then, by putting all together and integrating on $(0,t)$, for any $t \in (0,T)$, we deduce that 
\begin{align*}
    &\mezzo \norm{\nabla \phi}^2_H + \gamma \norm{\phi}^2_H + \delta \norm{\sigma}^2_H + \mezzo \int_\Omega F(\phi) \, \de x \\
    & \qquad + \mezzo \int_0^t \norm{\nabla \mu}^2_H \, \de s + \int_0^t \norm{\nabla (\sigma + \chi (1-\phi))}^2_H \, \de s + \int_0^t \norm{\sqrt{P(\phi)}(\sigma + \chi (1-\phi) - \mu)}^2_H \, \de s \\
    & \quad \le \max\left\{ \mezzo, \gamma \right\} \norm{\phi_0}^2_V + \delta \norm{\sigma_0}^2_H + C \int_0^T \left( \int_\Omega F(\phi) \, \de x \right) \, \de t + C \int_0^T \norm{\sigma}^2_H \, \de t + C.
\end{align*}
Hence, by applying Gronwall's lemma, we obtain the following uniform estimates:
\begin{equation}
\label{eq:weak:unifbounds}
    \begin{split}
        & \norm{\phi}_{\LT \infty V} \le C, \quad \norm{\sigma}_{\LT \infty H} \le C \\ 
        & \norm{\nabla \mu}_{\LT 2 H} \le C, \quad \norm{\nabla (\sigma + \chi (1 - \phi))}_{\LT 2 H} \le C \\
        & \norm{F(\phi)}_{\LT \infty {\Lx1}} \le C, \quad \norm{\sqrt{P(\phi)}(\sigma + \chi (1-\phi) - \mu)}_{\LT 2 H} \le C.
    \end{split}
\end{equation}
Moreover, by \eqref{eq:weak:mean_mu} and \eqref{eq:weak:unifbounds}, by means of Poincaré--Wirtinger's inequality, we also infer that
\begin{equation}
    \label{eq:weak:mu_l2v}
    \norm{\mu}_{\LT 2 V} \le C.
\end{equation}
We now observe that, by using H\"older's inequality and Sobolev's embeddings, for any $\zeta \in V$ we have that
\begin{align*}
    \abs{ (P(\phi) (\sigma + \chi (1 - \phi) - \mu), \zeta) }_H & \le \norm{P(\phi)}_{\Lx{3/2}} \norm{\sigma + \chi (1 - \phi) - \mu}_{\Lx6} \norm{\zeta}_{\Lx6} \\
    & \le C \norm{P(\phi)}_{\Lx{3/2}} \norm{\sigma + \chi (1 - \phi) - \mu}_V \norm{\zeta}_V.
\end{align*}
Next, since $\phi \in \LT \infty V \hookrightarrow \LT \infty {\Lx6}$ by \eqref{eq:weak:unifbounds}, by hypothesis \ref{ass:p} it follows that
\[ \norm{P(\phi)}_{\Lx{3/2}} \le \left( \int_\Omega c_5 (1 + \abs{\phi}^q)^{3/2} \, \de x \right)^{2/3} \le \left( \int_\Omega C ( 1 + \abs{\phi}^{3q/2} ) \, \de x \right)^{2/3} \in \Lt\infty, \]
considering that $1 \le q \le 4$.
Hence, since $\sigma + \chi(1-\phi) - \mu \in \LT 2 V$ by \eqref{eq:weak:unifbounds}, 
\[ \norm{P(\phi)(\sigma + \chi (1 - \phi) - \mu)}_{\LT 2 {V^*}} \le C, \]
and, subsequently, by comparison in \eqref{varform:phi} and \eqref{varform:sigma}, one can easily see that 
\begin{equation}
    \label{eq:weak:h1vstar}
    \norm{\phi}_{\HT 1 {V^*}} + \norm{\sigma}_{\HT 1 {V^*}} \le C.
\end{equation}
Note that the additional terms $-\hh(\phi)c$ and $(1-\sigma)$ do not make the analysis more difficult, since they are both easily bounded in $\LT 2 H$.
Finally, we just need to recover the additional $\LT 2 {\Hx3}$ regularity for $\phi$.
To do this, we consider equation \eqref{eq:mu}, which can be seen as an elliptic equation for $\phi$, namely
\begin{equation}
    \label{eq:weak:phi_elliptic}
    - \Delta \phi = \mu - F'(\phi) - \chi \sigma,
\end{equation}
together with the boundary condition $\partial_{\n} \phi = 0$, which is enclosed in the corresponding variational formulation \eqref{varform:mu}.
First, observe that by hypothesis \ref{ass:fconv} 
\[
    \int_\Omega \abs{F'(\phi)}^2 \, \de x \le \int_\Omega C ( 1 + \abs{\phi}^{2(s - 1)} ) \, \de x,
\]
which is uniformly bounded in $\Lt\infty$ if $s \le 4$, since $\phi \in \LT \infty {\Lx6}$. 
Then, the right-hand side of \eqref{eq:weak:phi_elliptic} is uniformly bounded in $\LT 2 H$, therefore by elliptic regularity we can infer that 
\begin{equation}
    \label{eq:weak:phil2w}
    \norm{\phi}_{\LT 2 W} \le C.
\end{equation}
Next, by means of Gagliardo--Nirenberg's interpolation inequality with suitable choices of the parameters, one can see that the following embeddings hold:
\begin{equation}
    \label{eq:weak:interpolation}
    \begin{split}
    & \nabla \phi \in \LT \infty H \cap \LT 2 V \hookrightarrow L^{10/3}(Q_T) \\  
    & \phi \in \LT \infty V \cap \LT 2 W \hookrightarrow L^{10}(Q_T).
    \end{split}
\end{equation}
Indeed, one needs to use $p=10/3$, $j=0$, $N=3$, $r=2$, $m=1$, $\alpha = 3/5$ and $q = 2$ for the first embedding, while $p=10$, $j=0$, $N=3$, $r=2$, $m=2$, $\alpha = 1/5$ and $q = 6$ for the second one.
With such information, by using also H\"older's inequality in $Q_T$ with $1/2 = 1/5 + 3/10$ and hypothesis \ref{ass:fconv}, one can also estimate
\begin{align*}
    \int_0^T \int_\Omega \abs{\nabla F'(\phi)}^2 \, \de x \, \de t & \le \int_0^T \int_\Omega \abs{F''(\phi) \nabla \phi}^2 \, \de x \, \de t \\
    & \le \norm{\nabla \phi}_{\Lqt{10/3}} \norm{F''(\phi)}_{\Lqt5} \\
    & \le \norm{\nabla \phi}_{\Lqt{10/3}} \int_0^T \int_\Omega C ( 1 + \abs{\phi}^{5(s-2)} ) \, \de x \, \de t \le C, 
\end{align*}
since \eqref{eq:weak:interpolation} hold and $s \le 4$.
Then, the right-hand side of \eqref{eq:weak:phi_elliptic} is now also bounded in $\LT 2 V$, hence elliptic regularity theory entails that 
\begin{equation}
    \label{eq:weak:l2h3}
    \norm{\phi}_{\LT 2 {\Hx3}} \le C.
\end{equation}

At this point, one could pass to the limit in the Galerkin discretisation scheme as in \cite[Theorem 1]{FGR2015_TumGrowth} and \cite[Theorem 3.1]{GY2020} and recover all the estimates \eqref{eq:weak:unifbounds}, \eqref{eq:weak:mu_l2v}, \eqref{eq:weak:h1vstar} and \eqref{eq:weak:l2h3} in the limit, as well as the validity of the system in weak form. 
We recall that the proof above works if $s \le 4$ in hypothesis \ref{ass:fconv}; if $s \in (4,6)$ instead one can argue as in \cite[Theorem 1]{FGR2015_TumGrowth}.
This concludes the proof of Theorem \ref{thm:weaksols}.

\section*{Acknowledgements} A. Agosti, C. Cavaterra, M. Fornoni and E. Rocca have been partially supported by the MIUR-PRIN Grant 2020F3NCPX ``Mathematics for industry 4.0 (Math4I4)''. 
C. Cavaterra has been partially supported by the MIUR-PRIN Grant 2022  
	``Partial differential equations and related geometric-functional inequalities''. 
 A. Agosti and E. Rocca also acknowledge the support of Next Generation EU Project No.P2022Z7ZAJ (A unitary mathematical framework for modelling muscular dystrophies).
The research of E. Beretta has been supported by the Project AD364 - Fund of NYU Abu Dhabi.
 The research of C. Cavaterra is part of the activities of ``Dipartimento di Eccellenza 2023-2027'' of Universit\`a degli Studi di Milano.
A Agosti, C. Cavaterra, M. Fornoni and E. Rocca are members of  
	GNAMPA (Gruppo Nazionale per l'Analisi Matematica, la Probabilit\`a e le loro Applicazioni)
	of INdAM (Istituto Nazionale di Alta Matematica).
 A. Agosti and M. Fornoni also acknowledge some support from the GNAMPA of INdAM through the GNAMPA project CUP E53C23001670001.



\footnotesize

\end{document}